\documentclass[authoryear,11pt]{elsarticle}
\usepackage{cancel}
\usepackage{caption}
\usepackage{color}
\usepackage{lscape}
\usepackage{afterpage}
\usepackage{pstricks}
\usepackage{pst-plot}
\usepackage{longtable}
\usepackage{dcolumn}
\usepackage{pst-node}
\usepackage{amsmath}
\usepackage{amssymb}
\usepackage{amsthm}
\usepackage{multirow}
\usepackage{colortab}
\usepackage{color}
\usepackage{rotating}
\usepackage{array}
\usepackage{textcomp}
\usepackage{pstricks, pst-node}
\usepackage{pst-all}
\usepackage{algorithm,algorithmic}
\usepackage{url}
\usepackage{soul}
\usepackage{hyperref}
\psset{arrows=->, labelsep=3pt, mnode=circle}
\usepackage{lineno}
\usepackage{booktabs}
\usepackage{tabularx}
\usepackage{subcaption}
\usepackage{eurosym}
\usepackage[normalem]{ulem}
\usepackage{tikz}
\usetikzlibrary{decorations.pathreplacing,calligraphy}

\newfont{\rams}{msbm10 scaled\magstep1}
\newcommand{\rea}{\mbox{\rams \symbol{'122}}}

\newcommand{\nat}{\mbox{\rams \symbol{'116}}}

\newenvironment{resumeT}{\begin{list}{}{\setlength{\rightmargin}{\leftmargin}}\item[]
{\centering {\bf \it~~~}
\par}\item[]\ignorespaces}{\unskip\end{list}}

\newtheorem{example}{Example}[section]

\usepackage{setspace}
\begin{document}
%\begin{center}
%    \large{\textbf{Title page for the paper}}\\[8mm]
%    \Large{\textbf{Deck of Cards method for Hierarchical,\\
%    Robust and Stochastic Ordinal Regression}}\\[8mm]
%\end{center}
%\textbf{First author}\\
%Name: Salvatore Corrente\\
%Affiliation: Department of Economics and Business, University of Catania\\
%Address: Corso Italia 55, 95129, Catania, Italy\\
%Email: salvatore.corrente@unict.it\\[4mm]
%\textbf{Second author}\\
%Name: Salvatore Greco\\
%Affiliation: Department of Economics and Business, University of Catania\\
%Address: Corso Italia 55, 95129, Catania, Italy\\
%Email: salgreco@unict.it\\[4mm]
%\textbf{Third author}\\
%Name: Silvano Zappalà\\
%Affiliation: Department of Economics and Business, University of Catania\\
%Address: Corso Italia 55, 95129, Catania, Italy\\
%Email: silvano.zappala@phd.unict.it\\[4mm]
%\textbf{Corresponding author}\\
%Name: Silvano Zappalà\\
%Affiliation: Department of Economics and Business, University of Catania\\
%Address: Corso Italia 55, 95129, Catania, Italy\\
%Email: silvano.zappala@phd.unict.it\\
%Phone number: +393272461669
%\newpage
\title{Deck of Cards method for Hierarchical, \\ Robust and Stochastic Ordinal Regression}

\author[Eco]{\rm Salvatore Corrente}
\ead{salvatore.corrente@unict.it}
\author[Eco]{\rm Salvatore Greco}
\ead{salgreco@unict.it}
\author[Eco]{\rm Silvano Zappalà}
\ead{silvano.zappala@phd.unict.it}

\address[Eco]{Department of Economics and Business, University of Catania, Corso Italia, 55, 95129  Catania, Italy}

\date{}
\maketitle
\vspace{-1cm}

%%%%%%%%%%%%%%%%%%%%

\begin{resumeT}
{\large {\bf Abstract:}}  
We consider the recently introduced application of the Deck of Cards Method (DCM) to ordinal regression proposing two extensions related to two main research trends in Multiple Criteria Decision Aiding, namely scaling and ordinal regression generalizations. On the one hand, procedures, different from DCM (e.g. AHP, BWM, MACBETH) to collect and elaborate Decision Maker's (DM's) preference information are considered to define an overall evaluation of reference alternatives. On the other hand, Robust Ordinal Regression and Stochastic Multicriteria Acceptability Analysis are used to offer the DM more detailed and realistic decision-support outcomes. More precisely, we take into account preference imprecision and indetermination through a set of admissible comprehensive evaluations of alternatives provided by the whole set of value functions compatible with DM's preference information rather than the univocal assessment obtained from a single value function. In addition, we also consider alternatives evaluated on a set of criteria hierarchically structured. The methodology we propose allows the DM to provide precise or imprecise information at different levels of the hierarchy of criteria. Like scaling procedures, the compatible value function we consider can be of a different nature, such as weighted sum, linear or general monotone value function, or Choquet integral. Consequently, the approach we propose is versatile and well-equipped to be adapted to DM's characteristics and requirements. The applicability of the proposed methodology is shown by a didactic example based on a large ongoing research project in which Italian regions are evaluated on criteria representing Circular Economy, Innovation-Driven Development and Smart Specialization Strategies.
%\vspace{1cm}
%%%%%%%%%%%%%%%%%%%%
%%%%%%% Keywords %%%%%%%
%%%%%%%%%%%%%%%%%%%%

\vspace{0,3cm}
\end{resumeT}

\pagenumbering{arabic}

\section{Introduction}
%\textcolor{brown}{Vediamo come va.}
Complex decisions require to take into consideration a plurality of points of view. For example, decisions related to circular economics \citep{stahel2016circular} have to be based on environmental aspects, economic aspects such as gross domestic product or employment rate, and also sustainability \citep{elliott2012introduction} or smart specialization \citep{Smartspecialization} aspects. Evaluation and comparison of different alternatives in complex decision problems –  in our example, feasible economic policies and strategies – ask for adequately articulated models embedded in advanced decision support procedures. They allow to aggregate partial evaluations with respect to the many considered criteria into a global evaluation. In addition, observe that in such complex problems, many heterogeneous stakeholders, as well as several experts in different domains, are involved together with a plurality of policy-makers. The results supplied by the formal decision model adopted for these types of complex problems heavily depend on the adopted parameters, so it is necessary to verify the stability of the obtained recommendation at the variation of the considered parameters. Just to give an example, if the trade-offs between the different criteria are represented by the weights assigned to them, it is reasonable to check if and how the comparisons between the alternatives change with a variation of these weights. Moreover,  all the above-mentioned actors in the decision process – experts, stakeholders and policy-makers – in general, have no expertise in the methodology and techniques for decision analysis. Consequently, the decision support procedure for complex decisions needs:
\begin{itemize}
\item to permit handling several heterogeneous criteria;
\item to collect the preference information from the different actors with a procedure as simple as possible; 
\item to collect rich preference information permitting to define in the most precise way the parameters of the decision models required to deal with the decision problem at hand;
\item to provide information about the stability of the comparisons between alternatives when the parameters of the decision model are perturbed. 
\end{itemize}
Multiple Criteria Decision Aiding (MCDA) (for an updated collection of state-of-the-art surveys see \citealt{GrecoEhrgottFigueira2016}, while for a review of the evolution of MCDA over the past 50 years, with a discussion on the perspectives and a future research agenda, see \citealt{greco2024fifty}) has developed a certain number of concepts, methods, techniques and procedures that can deal with the above requirements. Among them, we consider the ordinal regression \citep{JacquetLagrezeSiskos1982}. It asks the Decision Maker (DM) to compare pairwise a subset of the considered alternatives, called reference alternatives, in terms of preference. These comparisons define a preference model, very often a value function, aiming to represent in the most faithful way the DM’s preferences. As, in general, there is a plurality of value functions compatible with the DM’s preferences, a theoretical development of the ordinal regression, the Robust Ordinal Regression (ROR) \citep{GrecoMousseauSlowinski2008}, proposes to considers the whole set of decision models compatible with DM’s preferences answering to robustness concerns. 
%In this way, perturbation of the decision model preference parameters are taken into consideration. 
With respect to the richness of the preference information provided by the DM, recently it has been proposed to extend the ordinal regression paradigm considering preference information related not only to the mere preference between reference alternatives but to its intensity too \citep{barbati2023new}. This means the DM can specify that, among three alternatives $a$, $b$, and $c$, not only $a$ is preferred to $b$ and $b$ to $c$, but also that $a$ is preferred to $b$ more strongly than $b$ to $c$. To enable the DM to present this type of information in a clear and understandable way, the Deck of Cards Method (DCM) \citep{FigueiraRoy2002,AbastanteEtAl2022a} has been applied. To this aim, the DM is initially provided with a card for each reference alternative and a certain number of blank cards. Firstly, the DM rank orders the cards according to their preferences on the related alternatives. Secondly, the DM expresses the intensity of these preferences including a certain number of blank cards between two consecutive alternatives so that the greater the number of blank cards, the greater the preference between the alternatives. Once this preference information has been provided by the DM and a value function (for example, a weighted sum, an additive value function or a Choquet integral) has been selected to restore it, the instance of the value function that better represents these preferences can be obtained solving a linear programming problem. In particular, it is solved by minimizing the sum of the deviations between the values given to the reference alternatives by the DCM and the value assigned by the induced value function that is then used to give a value to all considered alternatives (not only reference ones). 
%\textcolor{blue}{In the DCM the DMs rank alternatives and indicate preference intensity by spacing blank cards between them. The optimal value function is then computed using linear programming.} 
This approach has been called Deck-of-Cards-based Ordinal Regression (DOR) \citep{barbati2023new}. It has the advantage of collecting very rich preference information, that is, not only information about the preference ranking of reference alternatives as in ordinal regression but also information related to the intensity of preference between them. Even if DOR considers richer information than ordinal regression, it keeps its main weak points. In particular, 
%DOR does not take into account the possibility that there exists a plurality of value functions that represent the DM's preference with similar faithfulness. 
DOR does not consider the possibility that multiple value functions may represent the DM's preferences with comparable accuracy. Moreover, DOR does not admit that the DM could provide imprecise preference information, for example, because they cannot exactly say how many blank cards have to be included between two alternatives. Another aspect not considered in DOR is that the set of criteria may be structured in a hierarchical way. Consequently, in this paper, we propose some extensions of DOR permitting to handle all the above problems and, more precisely:
\begin{itemize}
\item to define overall evaluations of reference alternatives by procedures collecting and elaborating DM's preference information different from DCM, such as AHP \citep{Saaty1977},  BWM \citep{Rezaei2015} and MACBETH \citep{BanaCostaVansnick1994}, which, for different reasons (for example because already known by the DM or the analyst) can be more appropriate in the specific decision problem at hand; 
\item to consider the whole set of value functions compatible with the DM's preferences for DOR method, applying the ROR \citep{GrecoMousseauSlowinski2008}. It builds a necessary and a possible preference relation holding in case an alternative is preferred to another for all or for at least one compatible value  function;
\item to consider the whole set of value functions compatible with the DM's preferences for DOR method, applying the Stochastic Multicriteria Acceptability Analysis (SMAA) \citep{LahdelmaHokkanenSalminen1998,PelissariEtAl2020}. Based on a sample of compatible value functions, it computes the frequency with which an alternative is ranked in a certain position or the frequency with which an alternative is preferred to another one;
\item to generalize the DOR method by taking into account the possibility that the DM provides imprecise information about the number of blank cards between reference alternatives;
%, giving only a minimal number and a maximal number of blank cards between them, or only one of these bounds or, even, no one of them;
\item to permit the DOR approach to take into account problems presenting a set of criteria hierarchically structured. This allows the DM to give preference information and to obtain final recommendations on the considered problem both at global and partial levels. To this aim, the Multiple Criteria Hierarchy Process (MCHP) \citep{CorrenteGrecoSlowinski2012} will be adapted to this context.    
\end{itemize}
The proposed extensions make DOR method very flexible and adaptable to real-world applications, taking into account two main research trends in MCDA: on the one hand, the scaling procedures such AHP, BWM and MACBETH, and, on the other hand, the extensions of ordinal regression, that is, ROR and Stochastic Ordinal Regression \citep{KadzinskiTervonen2013b}. 

The paper is structured in the following way. In the next section, we introduce the basic framework recalling the basic principles of DOR methodology. Section \ref{originalmodel} extends DOR to take into account the previously mentioned issues discussing in detail the possibility of applying scaling procedures different from the DCM. The DOR's extensions to ROR, SMAA and MCHP are provided in Sections \ref{ROR}, \ref{SMAASection} and \ref{MCHPSection}, respectively. How to include imprecise information provided by the DM in DOR is described in Section \ref{ImpreciseInformation}.
%The different preference models to be used to represent the DM's preferences are shown in Section \ref{PreferenceModels}.
A didactic example, based on a large ongoing research project, is presented in Section \ref{ExampleSection} to illustrate the application of the new proposal. Finally, Section \ref{ConclusionSection} collects conclusions and future directions of research.

%%%%%%%%%%%%%%%%%%%%%%%%%%%%%%%%%%%%%%%%%%%%%%%%%%%%%%%%%%%%%
\section{The Deck of Cards based Ordinal Regression}%%%%%%%%%
\label{originalmodel}
%%%%%%%%%%%%%%%%%%%%%%%%%%%%%%%%%%%%%%%%%%%%%%%%%%%%%%%%%%%%%
Let us denote by $A=\{a_1,a_2,\ldots\}$ a set of alternatives evaluated on a coherent family of criteria \citep{Roy1996} $G = \{g_1, \ldots, g_n\}$. We shall assume that each criterion $g_i\in G$ is a real-valued function $g_i:A\rightarrow\rea$ and, consequently, $g_i(a_j)\in\rea$ is the evaluation of $a_j\in A$ on criterion $g_i$. For each $g_i\in G$, $X_i =\{x_{i}^{1}, \ldots, x_{i}^{m_i}\}$ is the whole set of evaluations that can be taken on $g_i$ and they are such that $x_i^{1}<\cdots<x_i^{m_i}$. In the following, for the sake of simplicity and without loss of generality, let us assume that all criteria have an increasing direction of preference (the greater $g_i(a_j)$, the better $a_j$ is on $g_i$). Thus, ($x_1^{m_1}, \ldots, x_n^{m_n}$) represents the ideal alternative, achieving the highest possible evaluation across all criteria, whereas ($x_1^{1}, \ldots, x_n^{1}$) represents the anti-ideal alternative, with the lowest possible evaluation on all criteria. 

Looking at the evaluations of the alternatives on the considered criteria, the only objective information that can be obtained is the dominance relation for which $a_{j_1}\in A$ dominates $a_{j_2}\in A$ iff $a_{j_1}$ is at least as good as $a_{j_2}$ for all criteria and better for at least one of them ($g_i(a_{j_1}) \geqslant g_i(a_{j_2}),$ for all $i=1,\ldots,n$, and there exists at least one $g_i \in G$ such that $g_i(a_{j_1}) > g_i(a_{j_2})$). Even if it is objective that if $a_{j_1}$ dominates $a_{j_2}$, then $a_{j_1}$ is at least as good as $a_{j_2}$, the dominance relation is quite poor since, comparing two alternatives, in general, one is better than the other for some criteria but worse for the others. For this reason, to provide a final recommendation on the alternatives at hand, one needs to aggregate the evaluations taken by any alternative $a \in A$ on the criteria from $G$. In the field of MCDA, three main approaches have been proposed for aggregation of criteria: 
\begin{enumerate}
	\item using a value function $U: \mathbb{R}^n\rightarrow \mathbb{R}$ non decreasing in each of its arguments, such that, for any  $a_{j_1}, a_{j_2} \in A$, if $U(g_1(a_{j_1}),\ldots,g_n(a_{j_1})) \geqslant U(g_1(a_{j_2}),\ldots,g_n(a_{j_2}))$, then $a_{j_1}$ is comprehensively at least as good as $a_{j_2}$ \citep{KeeneyRaiffa1976};
	\item using an outranking relation $S$ on $A$ defined on the basis of evaluations taken by alternatives from $A$ on criteria from $G$,  such that for any  $a_{j_1}, a_{j_2} \in A, a_{j_1} S a_{j_2}$ means $a_{j_1}$ is comprehensively at least as good as $a_{j_2}$; the relation $S$ is reflexive, but, in general, neither transitive nor complete \citep{Roy1996};
	\item using a set of decision rules expressed in natural language, such as, for example, if ``$a$ is fairly preferred to $a^{\prime}$ on criterion $g_{i_1}$ and extremely preferred on criterion $g_{i_2}$, then $a$ is comprehensively at least as good as $a^{\prime}$''; the decision rules are induced from some examples of decisions provided by the DM \citep{GrecoMatarazzoSlowinski2001}.
\end{enumerate}
In this paper, we consider the aggregation of criteria through a value function $U$ and we take into consideration a recently proposed methodology called DOR (Deck of Cards based Ordinal Regression) \citep{barbati2023new}. It collects the preference information from the DM using the Deck of Cards method (DCM) \citep{FigueiraRoy2002,AbastanteEtAl2022a} and defines the parameters of the value function $U$ using an ordinal regression approach \citep{JacquetLagreze1982} by means of a mathematical programming problem that allows to represent the DM's preferences. 

Considering a set of reference alternatives $A^R\subseteq A$, in the following, we shall briefly review the main steps of the DOR method:
\begin{enumerate}
\item The DM has to rank-order alternatives in $A^{R}$ from the least to the most preferred in sets $L_1,L_2,\ldots,L_s\subseteq A^R$. Alternatives in $L_{h+1}$ are preferred to alternatives in $L_{h}$ for all $h=1,\ldots,s-1$ and alternatives in $L_h$ are indifferent among them for all $h=1,\ldots,s$;
\item The DM can include a certain number of blank cards $e_h$ between sets $L_{h}$ and $L_{h+1}$. The greater $e_h$, the more alternatives in $L_{h+1}$ are preferred to alternatives in $L_h$. Observe that the absence of blank cards between $L_{h}$ and $L_{h+1}$ does not mean that the alternatives in $L_h$ are indifferent to the alternatives in $L_{h+1}$ but that the difference between the value of the alternatives in $L_{h+1}$ and the value of the alternatives in $L_h$ is minimal; 
\item Following \cite{AbastanteEtAl2022a}, the DM has to provide the number of blank cards $e_0$ between the ``fictitious zero alternative" and the set of the least preferred alternatives, that is, $L_1$. In this context, the ``fictitious zero alternative" is a fictitious alternative $a_0$ having a null value, that is, $U(a_0)=0$. To make coherent the notation introduced above, let us assume $L_0=\{a_0\}$;
\item Each alternative $a \in A^R$ is assigned a value $\nu(a)$ such that  if $a\in L_{h+1}$ and $a^{\prime}\in L_h$, $\nu\left(a\right)=\nu\left(a^{\prime}\right)+(e_h+1)$, $h=1,\ldots,s-1$ and if $a\in L_1$, then $\nu(a)=e_0+1$. Consequently, in general, for $a\in L_h$ we have: $\nu(a)=\displaystyle\sum_{p=0}^{h-1} (e_{p}+1);$
\item The parameters of the value function $U$ are determined in a way that for all $a \in A^R$, $U(a)$ deviates as less as possible from $k\cdot\nu(a)$ with $k$ a scalarizing constant that can be interpreted as the value of a single blank card. More precisely, for all $a \in A^R$, one considers the positive and negative deviations $\sigma^+(a)$ and $\sigma^-(a)$ between $U(a)=U(g_1(a),\ldots,g_n(a))$ and $k\cdot\nu(a)$, that is $U(a)-\sigma^+(a)+\sigma^-(a)=k\cdot\nu(a)$. The sum of all deviations $\sigma^+(a)$ and $\sigma^-(a), a\in A^R$, is then minimized solving the following optimization problem having as unknown variables the constant $k$ and the deviations $\sigma^+(a)$ and $\sigma^-(a)$, in addition to the parameters of the value function $U$    
\begin{equation}\label{regressionmodel}
\begin{array}{l}
\;\;\overline{\sigma}=\min \displaystyle \sum_{a \in A^R} (\sigma^+(a)+\sigma^-(a)),\; \mbox{subject to},\\[1mm]
\left.
\begin{array}{l}
E^{Model}\\[2mm]
%U(x_1^{m_1}, \ldots, x_n^{m_n}) = 1, \\[2mm]
U(a) \geqslant U(a^\prime), \;\mbox{for all} \;a,a^\prime \in A^R\; \;\mbox{such that } \nu(a) \geqslant \nu(a^\prime),\\[2mm]
%U(a_{j_1})=U(a_{j_2}), \;\;\mbox{for all}\;\;a_{j_1},a_{j_2}\in L_j,\;\;\mbox{and for all}\;\;j=1,\ldots,s,\\[2mm]
U(a)-\sigma^{+}(a) +\sigma^{-}(a)=k\cdot \nu(a)\;\; \mbox{for all }a \in A^R,\\[2mm]
k \geqslant 0,\\[2mm]
\sigma^{+}(a) \geqslant 0, \;\sigma^{-}(a) \geqslant 0 \mbox{ for all } a \in A^R, \\[2mm]
\end{array}
\right\}E^{DM}
\end{array}
\end{equation}
\noindent where:
\begin{itemize}
\item $E^{Model}$ is the set of technical constraints related to the considered function $U$: for example, if $U$ is formulated as a \textit{weighted sum}, we have:
$$
\left.
\begin{array}{l}
U(a)=\displaystyle\sum_{i=1}^{n}w_i\cdot g_i(a) \mbox{  for all } a \in A,\\[1mm]
w_i\geqslant 0,\;\mbox{for all}\;i=1,\ldots,n, \\[1mm]
\displaystyle\sum_{i=1}^{n}w_i=1,\\[1mm]
\end{array}
\right\}E^{Model}_{WS}
$$
%}
\noindent while, if $U$ is formulated as a \textit{general additive value function}, we have:
$$
\left.
\begin{array}{l}
U(a)=\displaystyle\sum_{i=1}^{n} u_i(g_i(a))\;\mbox{for all}\; a\in A,\\[1mm]
u_i\left(x_i^{f-1}\right) \leqslant u_i\left(x_i^{f}\right),\; \mbox{for all}\;i= 1, \ldots, n\; \mbox{and for all}\; f = 2,\ldots,m_i,\\[2mm]
u_i\left(x_i^1\right) = 0, \; \mbox{for all}\; i = 1, \ldots, n,\\[2mm]
\displaystyle\sum_{i=1}^{n} u_i\left(x_i^{m_i}\right) = 1.
\end{array}
\right\}E^{Model}_{GA}
$$
%: for example, if $U$ is formulated as a weighted sum, we have
%$$
%\left.
%\begin{array}{l}
%U(a)=\displaystyle\sum_{i=1}^{n}w_i\cdot g_i(a) \mbox{  for all } a \in A^R\\%[1mm]
%w_i\geqslant 0,\;\mbox{for all}\;i=1,\ldots,n \\[1mm]
%\displaystyle\sum_{i=1}^{n}w_i=1.\\[1mm]
%end{array}
%\right\}E^{Model}_{WS} 
%$$
More details on different formulation of value function $U$ and, consequently, of the set of constraints $E^{Model}$ will be given in the \ref{PreferenceModels};
%\item $U(x_1^{m_1}, \ldots, x_n^{m_n}) = 1$ defines an upper bound for the value of all the alternatives in $A$\footnote{Since $\left(x_1^{m_1}, \ldots, x_n^{m_n}\right)$ is the performances vector of the \textit{ideal} alternative having the best performance on all the criteria, and all alternatives $\left(g_1(a),\ldots,g_n(a)\right)\in X_1\times \cdots \times X_n$ are weakly dominated by $\left(x_1^{m_1}, \ldots, x_n^{m_n}\right)$, the monotonicity constraints common to all considered preference models imply that $U\left(g_1(a),\ldots,g_n(a)\right)\leqslant U(x_1^{m_1}, \ldots, x_n^{m_n}) = 1$ for all $a\in A$.},
\item $U(a) \geqslant U(a^\prime)$ for all $a,a^\prime \in A^R$ such that $\nu(a) \geqslant \nu(a^{\prime})$, imposes that the ranking of alternatives from $A^R$ provided by the value function $U$ is concordant with the one provided by the DM,
%\item $U(a_{j_1}) = U(a_{j_2})$ imposes that alternatives $a_{j_1}$ and $a_{j_2}$ have the same value since they belong to the same subset $(L_j)$,
\item $U(a)-\sigma^{+}(a) +\sigma^{-}(a)=k\cdot \nu(a)$ imposes that the value of $a \in A^R$ is proportional to $\nu(a)$,
\item $k\geqslant0$ imposes the non-negativity of the scalarizing constant $k$,
\item $\sigma^{+}(a) \geqslant 0,\; \sigma^{-}(a) \geqslant 0, a \in A^R$ impose the non-negativity of the positive and negative deviations.
\end{itemize}
\end{enumerate}
The following cases can occur: 
\begin{description}
\item[1.] $\overline{\sigma}=0$ and $k>0$: in this case, the function $U$ obtained as solution of the LP problem (\ref{regressionmodel}) represents the preferences of the DM without any error;
\item[2.] $\overline{\sigma}=0$ and $k=0$: in this case, the function $U$ obtained as solution of the LP problem (\ref{regressionmodel}) does not represent the preferences of the DM since the scalarizing constant $k$ is zero and, therefore, all alternatives in $A$ have a null evaluation being indifferent among them;
\item[3.] $\overline{\sigma}>0$ and $k>0$: in this case, the function is able to represent the preferences of the DM with the minimal error $\overline{\sigma}$; 
\item[4.] $\overline{\sigma}>0$ and $k=0$: in this case, the function $U$ obtained as solution of the LP problem (\ref{regressionmodel}) does not represent the preferences of the DM since, as in above point 2., the scalarizing constant $k$ is null.
\end{description}

In the following, we shall present a small didactic example showing the DOR application.
\begin{example}
\label{didacticexample1}
Consider four regions ${ {\tt R_1}}, { {\tt R_2}}, { {\tt R_3}}, { {\tt R_4}}$ evaluated on a $0-100$ scale with respect to the three criteria of Circular Economy ($g_1$), Innovation-Driven Development ($g_2$),  and Smart Specialization ($g_3$), as shown in Table \ref{tab:evaluations}.
\begin{table}[!h]
\begin{center}
\caption{Evaluations of regions with respect to considered criteria }\label{tab:evaluations}
		\begin{tabular}{cccc}
	\hline
  Regions & Circular Economy: $g_1$ & Innovation-Driven Development: $g_2$ & Smart Specialization: $g_3$\\
	\hline
${ {\tt R_1}}$ & 90 &	100 &	80 \\ 
  ${ {\tt R_2}}$& 100 &	70 &	40 \\
${ {\tt R_3}}$& 30 &	50 &	60 \\
${ {\tt R_4}}$ & 20 &	40 &	40 \\
		\hline
		\end{tabular}
\end{center}
\end{table}
Using the DCM and taking into consideration a ``zero Region'' ${\tt R_0}$ as a reference of a null value level, the four regions ${\tt R_1}, {\tt R_2}, {\tt R_3}$ and ${\tt R_4}$ are ordered from worst to best so that $L_{0}=\{{\tt R_0}\}$, $L_{1}=\{{\tt R_4}\}$, $L_{2}=\{{\tt R_3}\}$, $L_{3}=\{{\tt R_2}\}$ and $L_{4}=\{{\tt R_1}\}$. Moreover, let $e_s$ be the number of blank cards (written between square brackets $[ \; ]$) between $L_s$ and $L_{s+1}$, $s=0,\ldots,3$, so that the complete preference information is the following:
$$ 
{{\tt R_0}} \;\; [5] \;\; {{\tt R_4}} \;\; [2] \;\;{{\tt R_3}} \;\; [5] \;\;{{\tt R_2}} \;\; [2] \;\;{{\tt R_1}}.
$$
By applying the DCM, we assign the following value to each project:
\begin{itemize}
	\item $\nu({\tt R_0}=[0,0,0])=0$,
	\item $\nu({\tt R_4}=[20,40,40])=\nu({\tt R_0})+e_0+1=6$,
	\item $\nu({\tt R_3}=[30,50,60])=\nu({\tt R_4})+e_1+1=9$,
	\item $\nu({\tt R_2}=[100,70,40])=\nu({\tt R_3})+e_2+1=15$,
	\item $\nu({\tt R_{1}}=[90,100,80])=\nu({\tt R_2})+e_3+1=18$.
	\end{itemize}
Considering a value function $U(\cdot)$ expressed in terms of a weighted sum,  the following LP problem has to be solved for the unknown variables $w_1, w_2, w_3, \sigma^+({\tt R_{i}}),\sigma^-({\tt R_i}), i=1,\ldots,4$, and $k$: 
\begin{equation}\label{regressionmodel_WS}
\begin{array}{l}
\;\;\;\;\overline{\sigma}=\min\displaystyle\sum_{i=1}^4\left(\sigma^+({\tt R_i})+\sigma^-({\tt R_i})\right),\; \mbox{subject to},\\[1,5mm]
\left.
\begin{array}{l}
\left.
\begin{array}{l}
U({\tt R_i})=w_1\cdot g_1({\tt R_i})+w_2\cdot g_2({\tt R_i})+w_3\cdot g_3({\tt R_i}) \;\; i=1,\ldots,4,\\[1,5mm]
w_1\geqslant 0,w_2\geqslant 0,w_3\geqslant 0, \\[1,5mm]
w_1+w_2+w_3=1,\\[1,5mm]
\end{array}
\right\}E^{Model}_{WS}\\[1,5mm]
\;\;U({\tt R_1})\geqslant U({\tt R_2})\geqslant U({\tt R_3})\geqslant U({\tt R_4}), \\[1,5mm]
\;\;U({\tt R_i})-\sigma^+({\tt R_i})+\sigma^-({\tt R_i})=k\cdot \nu({\tt R_i}),\;\; i=1,\ldots,4,\\[1,5mm]
%\displaystyle\sum_{j=1}^{m}w_j=1,\\[1mm]
\;\;k\geqslant 0,\\[1,5mm]
\;\;\sigma^+({\tt R_i}) \geqslant 0, \sigma^-({\tt R_i})\geqslant 0, \;\; i=1,\ldots,4.\\[1mm]
\end{array}
\right\}E^{DM}
\end{array}
\end{equation}
The solution of the LP problem (\ref{regressionmodel_WS}) yields a sum of errors $\overline{\sigma}=1.176$, a scaling constant $k=4.902$ and the following weights for the considered criteria $w_1= 0.471, w_2=0.176, w_3=0.353$. Using them, one can obtain the values listed in Table \ref{FirstExample}.  $\blacksquare$
%}
\begin{table}[!h]
\begin{center}
\caption{Scores assigned to regions by the value function $U(\cdot)$ obtained solving the LP problem (\ref{regressionmodel_WS})}\label{FirstExample}
\begin{tabular}{cccccc}
\hline
  Regions & $U(\cdot)$ & $\nu(\cdot)$ & $k\cdot\nu(\cdot)$& $\sigma^+(\cdot)$ & $\sigma^-(\cdot)$ \\
	\hline
${ {\tt R_1}}$ & 88.24 &	18 &	88.24 & 0&0\\ 
  ${ {\tt R_2}}$& 73.53 &	15 &	73.53 &0&0\\
${ {\tt R_3}}$& 44.12 &	9 &	44.12 &0&0\\
${ {\tt R_4}}$ & 30.59 &	6 &	29.21 &1.18&0\\
		\hline
		\end{tabular}
\end{center}
\end{table}
\end{example}
As previously explained, $k$ represents the value of a blank card and, consequently, it has to be greater than zero to avoid that all alternatives from $A$ are indifferent between them. Therefore, in order to ensure that this happens, we solve in an iterative way the following LP problem:
\begin{equation}
\label{CompatibleMaximizingk}
\begin{array}{l}
\;\;k^*=\max k\;\mbox{subject to},\\[1mm]
\left.
\begin{array}{l}
E^{DM},\\[1mm]
\displaystyle \sum_{a \in A^R} (\sigma^+(a)+\sigma^-(a))\leqslant \overline{\sigma}+\eta\left(\overline{\sigma}\right),
\end{array}
\right\}E^{DM{'}}
\end{array}
\end{equation}
\noindent where $\eta\left(\overline{\sigma}\right)$ is an admitted deterioration error with respect to the optimal value $\overline{\sigma}$ obtained solving (\ref{regressionmodel}). 
If $k^* = 0$, to keep the deterioration error as low as possible, $\eta\left(\overline{\sigma}\right)$ should be increased in increments of 0.01, starting from 0, until $k^* > 0$.
\\
In the following, we shall provide a small didactic example to show the necessity of maximizing the $k$ value.

\begin{example}\label{didacticexample2}
%Before demonstrating how to apply our proposal in the next section, we will first justify the 2-step procedure outlined in Section \ref{originalmodel} using a brief example. \\
Consider three alternatives evaluated on two criteria with increasing direction of preference, as shown in Table \ref{EvaluationsExample}:\\
\begin{minipage}[t]{0.45\textwidth}
\centering
\captionof{table}{Alternatives' evaluations on two criteria}
\label{EvaluationsExample}
%\centering
\begin{tabular}{ccc}
\toprule
Alternative&$g_1$&$g_2$\\
\midrule
$a$&$0.3$&$0.7$\\
$b$&$0.4$&$0.6$\\
$c$&$0.8$&$1$\\
\bottomrule
\end{tabular}
\end{minipage}
%\hspace{10mm}
\begin{minipage}[t]{0.45\textwidth}
\centering
\captionof{table}{DM's preference information}
\label{ExamplePreferences}
\begin{tabular}{cc}
\toprule
Alternative&Level\\
\midrule
$a$&$\nu(a) = 100$\\
&$e_1 = 29$\\
$b$&$\nu(b) = 70$\\
&$e_0 = 69$\\
\bottomrule
\end{tabular}
\end{minipage}\\[2mm]
Suppose the DM provides preference information for alternatives $a$ and $b$ as shown in Table \ref{ExamplePreferences}. Specifically, $a$ is preferred to $b$, putting 69 blank cards between $b$ and the zero level and 29 blank cards between $a$ and $b$. \\
%\begin{center}
%\captionof{table}{DM's preference information}
%\label{ExamplePreferences}
%\begin{tabular}{cc}
%\toprule
%Alternative&Level\\
%\midrule
%$a$&$\nu(a) = 100$\\
%&$e_1 = 29$\\
%$b$&$\nu(b) = 70$\\
%&$e_0 = 69$\\
%\bottomrule
%\end{tabular}
%\end{center}
%Assuming the DM's preference model follows a general additive value function, to check for an instance of the assumed preference model compatible with the DM's preferences one has to solve the following LP problem:
%\begin{center}
%\begin{equation}\label{ExampleGA}
%\begin{array}{l}
%\;\;\;\;\overline{\sigma} = \min \left \{ \sigma^{+}(a) + \sigma^{-}(a) + \sigma^{+}(b) + \sigma^{-} (b) \right \},\;\;\mbox{subject to},\\[3mm]
%\left.
%\begin{array}{l}
%\left.
%\begin{array}{l}
%U(a)= u_1(0.3) + u_2(0.7) \\[1,5mm]
%U(b) = u_1(0.4) + u_2(0.6)\\[1,5mm]
%U(c) = u_1(0.8) + u_2(1)\\[1,5mm]
%u_1(0.3) \leqslant u_1(0.4) \leqslant u_1(0.8)\\[1,5mm]
%u_2(0.6) \leqslant u_2(0.7) \leqslant u_2(1)\\[1,5mm]
%u_1(0.3) = 0\\[1,5mm]
%u_2(0.6) = 0\\[1,5mm]
%u_1(0.8) + u_2(1) = 1\\[1,5mm]
%\end{array}
%\right \} E^{Model}_{GA}\\[1,5mm]
%\;\;U(a) \geqslant U(b)\\[1,5mm]
%\;\;U(a) - \sigma^{+}(a) + \sigma^{-}(a) = k\cdot 100\\[1,5mm]
%\;\;U(b) - \sigma^{+}(b) + \sigma^{-}(b) = k\cdot 70\\[1,5mm]
%\;\;k \geqslant 0\\[1,5mm]
%\;\;\sigma^{+}(a),\; \sigma^{-}(a), \; \sigma^{+}(b), \;\sigma^{-}(b) \geqslant 0.
%\end{array}
%\right \} E^{DM}
%\end{array}
%\end{equation}
%\end{center}
Assuming the DM's preference model follows a general additive value function, the LP problem (\ref{regressionmodel}) has to be solved having as unknown variables $\sigma^{+}(a)$, $\sigma^{-}(a)$, $\sigma^{+}(b)$, $\sigma^{-}(b)$, $u_1\left(0.3\right)$, $u_1\left(0.4\right)$, $u_1\left(0.8\right)$, $u_2\left(0.6\right)$, $u_2\left(0.7\right)$, $u_2\left(1\right)$ and $k$ (see the LP problem (\ref{ExampleGA}) in \ref{sec:appendixB} for the extended formulation of the problem to be solved in this case). Solving this problem, one finds the solution shown in Table \ref{ExampleLP1}. 
\begin{center}
\captionof{table}{Marginal value functions obtained solving the LP problem (\ref{ExampleGA})}
\label{ExampleLP1}
\begin{tabular}{cccccccc}
\toprule
$u_1(0.3)$& $u_1(0.4)$&$u_1(0.8)$& $u_2(0.6)$&$u_2(0.7)$& $u_2(1)$ &$k$&$\overline{\sigma}$\\
\midrule
$0$&$0$&$0$&$0$&$0$&$1$&$0$&$0$ \\
\bottomrule
\end{tabular}
\end{center}
%Here, $k=0$ means that the blank card has null value. In order to maximize it, one has to solve the following LP problem:
%\begin{center}
    %\begin{equation}
    %\label{secondstepexample}
        %\begin{array}{l}
        %\;\;\ k^{*} = \max k,\; \text{subject to},\\[1,5mm]
        %\left.
        %\begin{array}{l}
             %E^{DM}\\[1,5mm]
             %\sigma^{+}(a)+\sigma^{-}(a)+\sigma^{+}(b)+\sigma^{-}(b) \leqslant \overline{\sigma} + \eta(\overline{\sigma}).
        %\end{array}
        %\right \} E^{DM^{'}}
        %\end{array}
    %\end{equation}
%\end{center}
Here, $k=0$ means that the blank card has null value. In order to maximize it, one has to solve the LP problem (\ref{CompatibleMaximizingk}) (see the LP problem (\ref{secondstepexample}) in \ref{sec:appendixB} for the extended formulation of the problem to be solved). Setting $\eta(\overline{\sigma}) = 0$ in problem (\ref{secondstepexample}), one gets the marginal values shown in Table \ref{ExampleLP2}:
\begin{center}
\captionof{table}{Marginal value functions obtained solving the LP problem (\ref{secondstepexample})}
\label{ExampleLP2}
    \setlength{\tabcolsep}{0.8\tabcolsep}
    \begin{tabular}{cccccccc}
    \toprule
    $u_1(0.3)$& $u_1(0.4)$&$u_1(0.8)$& $u_2(0.6)$&$u_2(0.7)$& $u_2(1)$ &$k$&$\overline{\sigma}$\\
    \midrule
    $0$& $0.4118$ &$0.4118$ &$0$&$0.5882$&$0.5882$&$0.0059$&$0$\\
    \bottomrule
    \end{tabular}
\end{center}
This value function accurately reflects the DM's preferences with a positive value for the blank card. $\blacksquare$

\end{example}
Some remarks are now at order:
\begin{enumerate}
\item LP problem (\ref{regressionmodel}) assumes that overall evaluations $\nu(\cdot)$ and $U(\cdot)$ are expressed on a ratio scale, that is, $U(a)/U(b)=\nu(a)/\nu(b)$ for all $a,b \in A$. If instead evaluations $\nu(\cdot)$ and $U(\cdot)$ are to be expressed on an interval scale, that is, $[U(a)-U(b)] / [U(c) - U(d)] = [\nu(a)-\nu(b)] / [\nu(c) - \nu(d)]$, for all $a,b,c,d \in A$ \citep{Stevens1946}, LP problem (\ref{regressionmodel}) has to be properly reformulated. More precisely, remembering that  the admissible transformations for a ratio scale are the multiplications by a positive constant, that is $U(\cdot)=k\cdot \nu(\cdot), k>~0$,  while the admissible transformations for an interval scale are the positive affine transformations, that is $U(\cdot)=k\cdot \nu(\cdot)+h$, $k>0, h\in \mathbb{R}$, we can conclude that  to handle interval scales the regression model  (\ref{regressionmodel}) has to be reformulated replacing the constraint
\begin{equation}\label{ratio}
U(a)-\sigma^{+}(a) +\sigma^{-}(a)=k\cdot \nu(a)\;\; \mbox{for all }a \in A^R,
\end{equation}
with the constraint
\begin{equation}\label{interval}
U(a)-\sigma^{+}(a) +\sigma^{-}(a)=k\cdot \nu(a)+h\;\; \mbox{for all }a \in A^R.
\end{equation}
Observe that if $\nu(\cdot)$ and $U(\cdot)$ are expressed on an interval scale, the information related to the number of blank cards $e_0$ between the worst alternative and the zero alternative is redundant, because it is not related to any difference $\nu(a)-\nu(b)$, $a,b \in A^R$. For the sake of simplicity, if not explicitly mentioned, we refer to the case in which $\nu(\cdot)$ and $U(\cdot)$ are expressed on a ratio scale. Anyway, what we discuss with respect to the ratio scale can be straightforwardly extended to the interval scale. \\[0,1cm]
\indent\; In the following, we shall continue the Example \ref{didacticexample1} showing how the described procedure changes once $U(\cdot)$ and $\nu(\cdot)$ are expressed on an interval scale.
\addtocounter{example}{-2}
\begin{example}[\textbf{continuation}]
Assuming $\nu(\cdot)$ and $U(\cdot)$ are expressed on an interval scale, the solution of the LP problem (\ref{regressionmodel_WS}) in which constraints 
$$
U({\tt R_i})-\sigma^{+}({\tt R_i})+\sigma^{-}({\tt R_i})=k\cdot\nu({\tt R_i}), \;\;i=1,\ldots,4
$$
\noindent are replaced by constraints
$$
U({\tt R_i})-\sigma^{+}({\tt R_i})+\sigma^{-}({\tt R_i})=k\cdot\nu({\tt R_i})+h, \;\;i=1,\ldots,4
$$
yields a null sum of errors $\overline{\sigma}$, $k=4.722$, $h=1.667$ and the following weights for the considered criteria: $w_1= 0.5$, $w_2=0.08$, $w_3=0.42$. Using them, one can obtain the values listed in Table \ref{Tab:regressionmodel_WS_interval}. 
\begin{table}[!h]
\begin{center}
\caption{Scores assigned to regions by the value function $U(\cdot)$ expressed on an interval scale and obtained solving the LP problem (\ref{regressionmodel_WS})}\label{Tab:regressionmodel_WS_interval}
		\begin{tabular}{cccccc}
	\hline
  Regions & $U(\cdot)$ & $\nu(\cdot)$ & $k\cdot\nu(\cdot)+h$& $\sigma^+(\cdot)$ & $\sigma^-(\cdot)$ \\
	\hline
${ {\tt R_1}}$ & 88.67 &	18 &	88.67 & 0&0\\ 
  ${ {\tt R_2}}$& 72.50 &	15 &	72.50 &0&0\\
${ {\tt R_3}}$& 44.17 &	9 &	44.17 &0&0\\
${ {\tt R_4}}$ & 30.00 &	6 &	30.00 &0&0\\
		\hline
		\end{tabular}
\end{center}
\end{table}    $\blacksquare$
\end{example}
\item Until now, we considered the DCM to assign an overall evaluation $\nu(a)$ to the reference alternatives $a \in A^R$. However, for several reasons, the Decision Maker (DM) may encounter difficulties in using the DCM. For instance, when a high number of blank cards is required to evaluate the alternatives, the method may be perceived as essentially asking the DM to directly assess their values. In some cases, this perception could make the approach seem unrealistic, contradicting the original purpose of the DCM, introduced by \cite{simos1990evaluer} to reduce the DM’s cognitive burden by using a limited number of cards, thus making it easier to manipulate a ratio scale.
 However, as already noted in \cite{barbati2023new}, depending on the specific problem, the DM's prior experience, and their individual predisposition, alternative procedures can be used to define overall evaluations $\nu(a)$. Among the simplest of these scaling methods, let us remember the direct rating and point allocations \citep{doyle1997judging}, as well as SMART and its extensions \citep{Edwards1977,EdwardsBarron1994}. Anyway, the most well-known among the methods proposed in literature to assign a value to alternatives on the basis of preference information provided by the DM are AHP \citep{Saaty1977},  BWM \citep{Rezaei2015} and MACBETH \citep{BanaCostaVansnick1994}.
%\item Until now, we considered the DCM to assign an overall evaluation $\nu(a)$ to the reference alternatives $a \in A^R$. However, on the basis of the previous experience and of the DM's predisposition, as already noted in \cite{barbati2023new}, one can use other procedures for defining overall evaluations $\nu(a)$. Among the simplest of these scaling methods, let us remember the direct rating and point allocations \citep{doyle1997judging}, as well as SMART and its extensions \citep{Edwards1977,EdwardsBarron1994}. Anyway, the most well-known among the methods proposed in literature to assign a value to alternatives on the basis of preference information provided by the DM are AHP \citep{Saaty1977},  BWM \citep{Rezaei2015} and MACBETH \citep{BanaCostaVansnick1994}. 
All these methods are based on pairwise comparisons of alternatives from $A^R$. For example, AHP provides evaluations $\nu(a), a \in A^R$, on a ratio scale, 
because it maintains the ratio $\nu(a)/\nu(b)=\nu'(a)/\nu'(b)$ between the evaluations of any $a, b \in A$, 
taking into consideration the DM's  pairwise comparisons $c_{a,b}$ of reference alternatives $a,b\in A^R$ in the classical  1-9 Saaty scale. Assuming that the DM retains $b$ not preferred to $a$, the points in the scale have the following interpretation: 
\begin{description}
\item[1-] $a$ and $b$ are equally preferable, denoted by $c_{a,b}=1$, translated to $\nu(a)=\nu(b)$,
\item[3-] $a$ is moderately preferred to $b$, denoted by $c_{a,b}=3$, translated to $\nu(a)=3 \cdot \nu(b)$,
\item[5-] $a$ is strongly preferred to $b$, denoted by $c_{a,b}=5$, translated to $\nu(a)=5 \cdot \nu(b)$,
\item[7-] $a$ is very strongly preferred to $b$, denoted by $c_{a,b}=7$, translated to $\nu(a)=7 \cdot \nu(b)$,
\item[9-] $a$ is extremely preferred to $b$, denoted by $c_{a,b}=9$, translated to $\nu(a)=9 \cdot \nu(b)$.
\end{description}
Values 2, 4, 6 and 8 denote a hesitation between 1-3, 3-5, 5-7 and 7-9, respectively. AHP considers the matrix $C=[c_{a,b}]_{a,b \in A^R}$, with $c_{a,b}=\frac{1}{c_{b,a}}$ for all $a,b\in A^{R}$. In general, there is no set of evaluations $\nu(a),  a\in A^R$, satisfying all equations $\nu(a)=c_{a,b} \cdot \nu(b)$, $a,b\in A^R$. AHP computes the values $\nu(a), a\in A^R$ as the components of the eigenvector $\hat{\nu}=[\hat{\nu}(a)]_{a \in A^R}$ associated to the maximum eigenvalue of $C$ $(\lambda_{max})$ normalized so that $\displaystyle\sum_{a \in A^{R}} \hat\nu(a)=1$. Formally we have
\begin{equation}
\label{AHP}
\begin{cases}
C\cdot \hat{\nu}=\lambda_{max}\cdot \hat{\nu}\\
\displaystyle{\sum_{a \in A} \hat{\nu}(a)=1}.
\end{cases}\,
\end{equation}
Procedures different from the eigenvector have been proposed to obtain the values $\nu(a), a\in A^{R}$, from the judgments in the matrix $C=[c_{a,b}]_{a,b\in A^{R}}$ and, among them, we remember the logarithmic least squares \citep{CrawfordWilliams1985}, %corresponding to the geometric mean
the least squares \citep{Jensen1984}, the weighted least squares \citep{ChuKalabaSpingarn1979,Blankmeyer1987}, the logarithmic least absolute values \citep{CookKress1988} and the geometric least square \citep{IsleiLockett1988}. \\[0,1cm] 
Some considerations related to decision psychology and to the great number of pairwise comparison judgments $c_{a,b}, a,b \in A^R$ required by AHP to fill the matrix $C$ suggested to consider only pairwise comparisons with respect to the best and the worst alternatives originating the BWM \citep{Rezaei2015}. A reduced number of pairwise comparison judgments is taken into account also in two recently proposed procedures \citep{AbastanteEtAl2019,CorrenteGrecoRezaei2024} that correct the DM's direct rating of alternatives through the overall evaluations provided by AHP or BWM to a limited number of reference alternatives. \\[0,1cm]
Requiring preference information similar to AHP, MACBETH provides evaluations $\nu(a), a \in A^R$, on an interval scale. In particular, assuming that the DM retains $b$ not preferred to $a$, MACBETH considers pairwise comparisons $c^M_{a,b}, a,b \in A^R$, on a 0-6 scale, having the following interpretation:
\begin{description}
\item[$0_M$-] $a$ and $b$ are equally preferable, denoted as $c^M_{a,b}=0$,
\item[$1_M$-] $a$ is very weakly preferred to $b$, denoted as $c^M_{a,b}=1$,
\item[$2_M$-] $a$ is weakly preferred to $b$, denoted as $c^M_{a,b}=2$, 
\item[$3_M$-] $a$ is  moderately preferred to $b$, denoted as $c^M_{a,b}=3$,
\item[$4_M$-] $a$ is  strongly preferred to $b$, denoted as $c^M_{a,b}=4$,
\item[$5_M$-] $a$ is  very strongly preferred to $b$, denoted as $c^M_{a,b}=5$, 
\item[$6_M$-] $a$ is extremely preferred to $b$, denoted as $c^M_{a,b}=6$.
\end{description}
Below, we consider the following MACBETH-like LP problem to define the overall evaluations $\nu(a), a \in A^{R}$: 
\begin{small}
\begin{subequations}
%\hspace{-10mm}
\label{MACBETH}
\renewcommand{\theequation}{\theparentequation.\arabic {equation}}
\begin{tikzpicture}[overlay, remember picture,thick]
\draw [decorate,decoration={brace,amplitude=8pt,mirror},line width=1.5pt] 
(13.3, -5.5) -- (13.3,-1.2) node[right, midway]{$\;\;\;E^{}_{Macb}$};
\end{tikzpicture}
\begin{align}
%\begin{array}{l}
&\hspace{-20mm}\max\gamma,\;\mbox{subject to},\notag\\[1mm]
%\left.
%\begin{array}{l}
&\hspace{-20mm}\nu(a^*)=100, \; \nu(a_*)=0\label{MACBETH_1}\\[1mm]
&\left.
\begin{array}{l}
\hspace{-23mm}\left[\nu(a)-\sigma^+(a)+\sigma^-(a)\right]-\left[\nu(b)-\sigma^+(b)+\sigma^-(b)\right]=0 \;\; \mbox{if }c^M_{a,b}=0,\\[1mm]
\hspace{-23mm}\left[\nu(a)-\sigma^+(a)+\sigma^-(a)\right]-\left[\nu(b)-\sigma^+(b)+\sigma^-(b)\right]=\delta_{e},\;\mbox{if}\;c_{a,b}^{M}=e,\; e=1,\ldots,6,\\[1mm]
\end{array}
\right\}a,b \in A^R\label{MACBETH_2} \\[1mm]
&\hspace{-20mm}\delta_{e+1} -\delta_e \geqslant \gamma, e=1,\ldots,5,\label{MACBETH_3}\\[1mm]
&\hspace{-20mm}\displaystyle{\sum_{a \in A^R}\left(\sigma^+(a)+\sigma^-(a)\right)\leqslant \overline{\sigma}},\label{MACBETH_4}\\[1mm]
&\hspace{-20mm}\sigma^+(a) \geqslant 0, \sigma^-(a)\geqslant 0.
\label{MACBETH_5}
%\end{array}
%\right\}
%\end{array}
\end{align}
\end{subequations}
\end{small}
The value of $\overline{\sigma}$ is obtained solving the following auxiliary optimization problem 
\begin{equation}
\hspace*{-8mm}
\label{MACBETH_AUS}
\;\;\;\;\displaystyle\overline{\sigma}=\max\sum_{a \in A^R}\left(\sigma^+(a)+\sigma^-(a)\right),\;\; \mbox{subject to}\;\; E'_{Macb}
\end{equation}
\noindent where $E'_{Macb}$ is obtained by $E_{Macb}$ replacing constraint (\ref{MACBETH_4}) with $\gamma\geqslant\varepsilon$; $\varepsilon$ is a small positive number defined to ensure that $\gamma$ assumes a strictly positive value; $a_{*}$ and $a^{*}$ are the worst and the best alternatives in $A^{R}$, respectively.
\addtocounter{example}{-1}
\begin{example}[\textbf{continuation}]
Applying AHP, let us assume the pairwise comparisons $c_{\tt R_i,\tt R_j}, i,j=1,2,3,4$, shown in the following matrix: 
\begin{equation*}
C= 
\bordermatrix{ & \tt R_1 & \tt R_2 & \tt R_3 & \tt R_4 \cr
\tt R_1 & 1 & 1 & 2 & 3 \cr
\tt R_2 & 1 & 1 & 2 & 2 \cr
\tt R_3 & 1/2 & 1/2 & 1 & 1 \cr
\tt R_4 & 1/3 & 1/2 & 1 & 1 }
\end{equation*}
Solving problem (\ref{AHP}) one obtains the normalized eigenvector $\hat\nu$ of $C$ providing the following overall evaluations for the regions  $\tt R_1, \tt R_2, \tt R_3$ and $\tt R_4$: $\hat\nu({\tt R_1})=0.36, \hat\nu({\tt R_2})=0.33, \hat\nu({\tt R_3})=0.16, \hat\nu({\tt R_4})=0.15$. As AHP expresses overall evaluations $\nu({\tt R_i}), i=1,2,3,4$, on a ratio scale, replacing them in the LP problem (\ref{regressionmodel_WS}) one gets the total error $\displaystyle\overline{\sigma}=\displaystyle\sum_{i=1}^4\left(\sigma^+({\tt R_i})+\sigma^-({\tt R_i})\right)=8.25$ and the results listed in Table \ref{Tab:regressionmodel_WS} where $k=248.87$ and the weights for the considered criteria are $w_1= 0.57, w_2=0.23, w_3=0.2$. 
\begin{table}[!h]
\begin{center}
\caption{Scores assigned to regions by the value function $U(\cdot)$ expressed on a ratio scale and obtained solving the LP problem (\ref{regressionmodel_WS}) with respect to AHP overall evaluations $\nu({\tt R_i}), i=1,2,3,4$.}\label{Tab:regressionmodel_WS}
		\begin{tabular}{cccccc}
	\hline
  Regions & $U(\cdot)$ & $\nu(\cdot)$ & $k\cdot\nu(\cdot)$& $\sigma^+(\cdot)$ & $\sigma^-(\cdot)$ \\
	\hline
${ {\tt R_1}}$ & 90.31 &	0.36 &	90.31 & 0&0\\ 
  ${ {\tt R_2}}$& 81.15 &	0.32 &	81.15 &0&0\\
${ {\tt R_3}}$& 40.58 &	0.16 & 40.58& 0&0\\
${ {\tt R_4}}$ & 28.59 &	0.15 &36.83& 0 &8.25\\
		\hline
		\end{tabular}
\end{center}
\end{table}
\\
Considering MACBETH and assuming that the pairwise comparisons $c^M_{\tt R_i,\tt R_j}, i,j=1,2,3,4$ are collected as follows
\begin{equation*}
C^M= 
\bordermatrix{ & \tt R_1 & \tt R_2 & \tt R_3 & \tt R_4 \cr
\tt R_1 & 0 & 1 & 4 & 5 \cr
\tt R_2 &  & 0 & 3 & 4 \cr
\tt R_3 &  &  & 0 & 1 \cr
\tt R_4 &  &  &  & 1 }
\end{equation*}
solving sequentially problems (\ref{MACBETH_AUS}) and (\ref{MACBETH}) one obtains the following overall evaluations for regions  ${\tt R_1},{\tt R_2}, {\tt R_3}$ and ${\tt R_4}$: $\nu({\tt R_1})=100, \nu({\tt R_2})=80, \nu({\tt R_3})=20, \nu({\tt R_4})=0$. 
As MACBETH expresses overall evaluations $\nu({\tt R_i}), i=1,2,3,4$, on an interval scale, using them in the LP problem (\ref{regressionmodel_WS}) where eqs. (\ref{ratio}) are replaced by eqs. (\ref{interval}) one gets a null total error $\left(\displaystyle\sum_{i=1}^4\left(\sigma^+({\tt R_i})+\sigma^-({\tt R_i})\right)=0\right)$ and the results listed in Table \ref{Tab:regressionmodel_WS_MACBETH} where $k=0.607$, $h=30$, and the weights for the considered criteria are $w_1= 0.5, w_2=0.29, w_3=0.21$. $\blacksquare$
\begin{table}[!h]
\begin{center}
\caption{Scores assigned to regions by the value function $U(\cdot)$ obtained solving the reformulation of LP problem (\ref{regressionmodel_WS}) for interval scales with respect to MACBETH overall evaluations $\nu({\tt R_i}), i=1,2,3,4$}\label{Tab:regressionmodel_WS_MACBETH}
		\begin{tabular}{cccccc}
	\hline
  Regions & $U(\cdot)$ & $\nu(\cdot)$ & $k\cdot\nu(\cdot)+h$& $\sigma^+(\cdot)$ & $\sigma^-(\cdot)$ \\
	\hline
${ {\tt R_1}}$ & 90.71 &	100 &	90.71 & 0&0\\ 
  ${ {\tt R_2}}$& 78.57 &	80 &	78.57 &0&0\\
${ {\tt R_3}}$& 42.14 &	20 & 42.14& 0&0\\
${ {\tt R_4}}$ & 30 &	0 &30.00& 0 &0\\
		\hline
		\end{tabular}
\end{center}
\end{table}
\end{example}
\item The number of cards $e_h, h=0,\ldots,s-1$ of the DCM can be interpreted in a more ordinal form in the sense that instead of $U(a)-U(b)=k\cdot e_h$ for all $a\in L_{h+1}, b\in L_h$, we could consider the much weaker relation for which if $e_h \geqslant e_{h'}, h,h'=0,\ldots,s$, then $U(a)-U(b) \geqslant U(c)- U(d)$, for all $a\in L_{h+1}$, $b\in L_h$, $c\in L_{h'+1}$ and $d\in L_{h'}$. Accepting this meaning of the blank cards, the regression problem (\ref{regressionmodel}) can be reformulated as follows:
\begin{subequations}\label{regressionmodel_weaker}\renewcommand{\theequation}{\theparentequation.\arabic {equation}}
%\begin{array}{l}
\begin{align}
&\hspace{-30mm}\max \gamma \mbox{ subject to},\notag\\[1mm]
%\left.
%\begin{array}{l}
&\hspace{-30mm}E^{Model}\label{eq:regressionmodel_weaker_1}\\[2mm]
%U(x_1^{m_1}, \ldots, x_n^{m_n}) = 1, \\[2mm]
&\hspace{-30mm}U(a) \geqslant U(a^\prime) \;\mbox{for all}\;a \in L_{h+1},a^\prime \in L_h, h=0,\ldots,s-1, \label{eq:regressionmodel_weaker_2}\\[2mm]
&\hspace{-30mm}U'(a)=U(a)-\sigma^{+}(a) +\sigma^{-}(a)\;\; \mbox{for all }a \in A^R,\label{eq:regressionmodel_weaker_3}\\[2mm]
&\hspace{-30mm}U'(a)-U'(b)=\delta_h\;\; \mbox{for all }a \in L_{h+1} \mbox{ and }b \in L_{h}, h=0,\ldots,s-1,\label{eq:regressionmodel_weaker_4}\\[2mm]
&\hspace{-30mm}\delta_h - \delta_{h'} \geqslant \gamma \mbox{ if } e_h> e_{h'}, \;h,h'=0,\ldots,s-1,\label{eq:regressionmodel_weaker_5}\\[2mm]
&\hspace{-30mm}\displaystyle \sum_{a \in A^R} (\sigma^+(a)+\sigma^-(a)) \leqslant \overline{\sigma}\label{eq:regressionmodel_weaker_6}\\[2mm]
&\hspace{-30mm}\sigma^{+}(a) \geqslant 0, \;\sigma^{-}(a) \geqslant 0 \mbox{ for all } a \in A^R.\label{eq:regressionmodel_weaker_7}
%\end{array}
%\right\}
%\end{array}
\end{align}
\begin{tikzpicture}[overlay, remember picture,thick]
    \draw [decorate,decoration={brace,amplitude=8pt,mirror},line width=1.5pt] 
    (13.5, 0.8) -- (13.5,7) node[right, midway] {$\;\;\;E^{}_{Ordinal}$};
\end{tikzpicture}
\end{subequations}
The value of $\overline{\sigma}$ is obtained solving the following auxiliary LP problem
\begin{equation}\label{regressionmodel_weaker_AUS}
%\begin{array}{l}
\;\;\overline{\sigma}=\min \displaystyle \sum_{a \in A^R} (\sigma^+(a)+\sigma^-(a)),\;\; \mbox{subject to} \;\;E^{\prime}_{Ordinal}
%\left.
%E^{\prime}_{Ordinal},
%\right\}
%\end{array}
\end{equation}
where $E^{\prime}_{Ordinal}$ is obtained from $E_{Ordinal}$ replacing constraints (\ref{eq:regressionmodel_weaker_5}) and (\ref{eq:regressionmodel_weaker_6}) with the following one 
$$
\delta_h \geqslant \delta_{h'}  \mbox{ if } e_h   \geqslant e_{h'}, \;\;h,h'=0,\ldots,s-1.
$$ 
%\begin{equation}\label{regressionmodel_weaker_AUS}
%\begin{array}{l}
%\;\;\overline{\sigma}=\min \displaystyle \sum_{a \in A^R} (\sigma^+(a)+\sigma^-(a)),\; \mbox{subject to},\\[1mm]
%\left.
%\begin{array}{l}
%E^{Model}\\[2mm]
%U(a) \geqslant U(a^\prime), \;a \in L_{h+1},a^\prime \in L_h, h=0,\ldots,s-1,\\[2mm]
%U'(a)=U(a)-\sigma^{+}(a) +\sigma^{-}(a)\;\; \mbox{for all }a \in A^R,\\[2mm]
%U'(a)-U'(b)=\delta_h\;\; \mbox{for all }a \in L_{h+1} \mbox{ and }b \in L_{h}, h=0,\ldots,s-1,\\[2mm]
%\delta_h \geqslant \delta_{h'}  \mbox{ if } e_h   \geqslant e_{h'}, \;\;h,h'=0,\ldots,s-1,\\[2mm]
%\sigma^{+}(a) \geqslant 0, \;\sigma^{-}(a) \geqslant 0 \mbox{ for all } a \in A^R. \\[2mm]
%\end{array}
%\right\}
%\end{array}
%\end{equation}
Observe that the DM's preference information could be collected also in terms of ratio instead of difference between overall evaluations of alternatives, that is,  if $e_h \geqslant e_{h'}, h,h'=0,\ldots,s$, then $U(a)/U(b) \geqslant U(c)/U(d)$, for all $a\in L_{h+1}$, $b\in L_h$, $c\in L_{h'+1},$ $d\in L_{h'}$. This interpretation of the meaning of blank cards requires the following reformulation of the regression problem:
\begin{subequations}\label{regressionmodel_weaker_ratio}\renewcommand{\theequation}{\theparentequation.\arabic {equation}}
%\begin{array}{l}
\begin{align}
&\hspace{-30mm}\max \gamma,\; \mbox{subject to},\notag\\[1mm]
%\left.
%\begin{array}{l}
&\hspace{-30mm}E^{Model}\label{regressionmodel_weaker_ratio_1}\\[2mm]
&\hspace{-30mm}U(a) \geqslant U(a^\prime), \;\mbox{for all}\;a \in L_{h+1},a^\prime \in L_h, h=1,\ldots,s-1,\label{regressionmodel_weaker_ratio_2} \\[2mm]
&\hspace{-30mm}U'(a)=U(a)-\sigma^{+}(a) +\sigma^{-}(a)\;\; \mbox{for all }a \in A^R,\label{regressionmodel_weaker_ratio_3}\\[2mm]
&\hspace{-30mm}U'(a)/U'(b)=\varphi_h\;\; \mbox{for all }a \in L_{h+1} \mbox{ and }b \in L_{h}, h=1,\ldots,s-1,\label{regressionmodel_weaker_ratio_4}\\[2mm]
&\hspace{-30mm}\varphi_h / \varphi_{h'} \geqslant \gamma \mbox{ if } e_h>e_{h'}, \;h, h'=1,\ldots,s-1,\label{regressionmodel_weaker_ratio_5}\\[2mm]
&\hspace{-30mm}\displaystyle \sum_{a \in A^R} \left(\sigma^+(a)+\sigma^-(a)\right) \leqslant \overline{\sigma}\label{regressionmodel_weaker_ratio_7}\\[1mm]
&\hspace{-30mm}\varphi_h \geqslant 0, h=1,\ldots,s-1,\label{regressionmodel_weaker_ratio_6}\\[2mm]
&\hspace{-30mm}\sigma^{+}(a) \geqslant 0, \;\sigma^{-}(a) \geqslant 0 \mbox{ for all } a \in A^R.\label{regressionmodel_weaker_ratio_8}
%\end{array}
%\right\}
%\end{array}
\end{align}
\begin{tikzpicture}[overlay, remember picture,thick]
    \draw [decorate,decoration={brace,amplitude=8pt,mirror},line width=1.5pt] 
    (13.5, 0.8) -- (13.5,8) node[right, midway] {$\;\;\;E^{}_{Ratio}$};
\end{tikzpicture}
\end{subequations}
The value of  $\overline{\sigma}$ is obtained solving the following auxiliary optimization problem
\begin{equation}\label{regressionmodel_weaker_ratio_AUS}
\begin{array}{l}
\;\;\overline{\sigma}=\min \displaystyle \sum_{a \in A^R} (\sigma^+(a)+\sigma^-(a)),\;\; \mbox{subject to},\;\; E^{\prime}_{Ratio},
\end{array}
\end{equation}
where $E^{\prime}_{Ratio}$ is obtained from $E_{Ratio}$ replacing constraints (\ref{regressionmodel_weaker_ratio_5}) and (\ref{regressionmodel_weaker_ratio_7}) with the following one
$$
\varphi_h \geqslant \varphi_{h'}  \mbox{ if } e_h   \geqslant e_{h'}, h, h'=1,\ldots,s-1.
$$
%\begin{equation}\label{regressionmodel_weaker_ratio_AUS}
%\begin{array}{l}
%\;\;\overline{\sigma}=\min \displaystyle \sum_{a \in A^R} (\sigma^+(a)+\sigma^-(a)),\; \mbox{subject to},\\[1mm]
%\left.
%\begin{array}{l}
%E^{Model}\\[2mm]
%U(a) \geqslant U(a^\prime), \;\mbox{for all}\;a \in L_{h+1},a^\prime \in L_h, h=1,\ldots,s-1, \\[2mm]
%U'(a)=U(a)-\sigma^{+}(a) +\sigma^{-}(a)\;\; \mbox{for all }a \in A^R,\\[2mm]
%U'(a)/U'(b)=\varphi_h\;\; \mbox{for all }a \in L_{h+1} \mbox{ and }b \in L_{h}, h=1,\ldots,s-1,\\[2mm]
%\varphi_h \geqslant \varphi_{h'}  \mbox{ if } e_h   \geqslant e_{h'}, h, h'=1,\ldots,s-1,\\[2mm]
%\varphi_h \geqslant 0, h=1,\ldots,s-1,\\[2mm]
%\sigma^{+}(a) \geqslant 0, \;\sigma^{-}(a) \geqslant 0 \mbox{ for all } a \in A^R. \\[2mm]
%\end{array}
%\right\}
%\end{array}
%\end{equation}
Even if the regression problem (\ref{regressionmodel_weaker_ratio}) is generally not linear, it can be handled by means of the many available nonlinear programming solvers. In the following, we shall refer to (\ref{regressionmodel_weaker}) and (\ref{regressionmodel_weaker_ratio}) as difference-based and ratio-based regression models, respectively. Observe that in regression problems (\ref{regressionmodel_weaker}) and (\ref{regressionmodel_weaker_ratio}) one can consider information about intensity of preferences expressed in terms of the DCM, 1-9 AHP scale or 0-6 MACBETH scale. However, we could consider information expressed in terms of intensity of preference of a merely ordinal type such as ``$a$ is preferred to $b$ at least as strongly as $c$ is preferred to $d$", without considering the predetermined scale expressed in terms of cards, 1-9 AHP scale or 0-6 MACBETH scale. This type of preference information has been discussed in ordinal regression and ROR in \cite{FigueiraGrecoSlowinski2009}.
\addtocounter{example}{-1}
\begin{example}[\textbf{continuation}]
Applying the difference-based LP regression model (\ref{regressionmodel_weaker}) to the preference information over regions ${\tt R_i}$, $i=1,\ldots,4$, collected through the DCM and considering a value function expressed in terms of a weighted sum, we get $\delta_1=\delta_3=30$ and $\delta_2=\delta_4=13.75$. Reminding that $e_1=e_3=5$ and $e_2=e_4=2$, this means that five and two blank cards have a value of 30 and 13.75, respectively. Solving the difference-based LP regression model (\ref{regressionmodel_weaker}) one gets a null total sum of the errors ($\overline{\sigma}=0$) and the following weights of criteria $w_1= 0.5, w_2=0.125, w_3=0.375$, so that, the overall evaluations of regions are obtained as listed in Table \ref{Tab:regressionmodel_WS_weaker}.\\
\begin{minipage}[t]{0.45\textwidth}
\begin{center}
\captionof{table}{Scores assigned to regions by the value function $U(\cdot)$ obtained solving the difference-based LP problem (\ref{regressionmodel_weaker})}\label{Tab:regressionmodel_WS_weaker}
\begin{tabular}{ccccc}
\toprule
Regions & $U(\cdot)$ & $U'(\cdot)$ &  $\sigma^+(\cdot)$ & $\sigma^-(\cdot)$ \\
\midrule
${ {\tt R_1}}$ & 87.5 &		87.5 & 0&0\\ 
  ${ {\tt R_2}}$& 73.75  &	73.75 &0&0\\
${ {\tt R_3}}$& 43.75 &		43.75 &0&0\\
${ {\tt R_4}}$ & 30.00 &		30.00 &0&0\\
\bottomrule
\end{tabular}
\end{center}    
\end{minipage}
\hfill
\begin{minipage}[t]{0.45\textwidth}
\begin{center}
\captionof{table}{Scores assigned to regions by the value function $U(\cdot)$ obtained solving the ratio-based regr. problem (\ref{regressionmodel_weaker_ratio})}\label{Tab:regressionmodel_WS_weaker_ratio}
\begin{tabular}{ccccc}
\midrule
Regions & $U(\cdot)$ & $U'(\cdot)$ &  $\sigma^+(\cdot)$ & $\sigma^-(\cdot)$\\
\midrule
${ {\tt R_1}}$ & 94.37 &		94.37 & 0&0\\ 
  ${ {\tt R_2}}$& 71.03  &	71.03 &0&0\\
${ {\tt R_3}}$& 47.55 &		47.55 &0&0\\
${ {\tt R_4}}$ & 35.79 &		35.79 &0&0\\
\bottomrule
\end{tabular}
\end{center}    
\end{minipage}\\[2mm]
\noindent Applying the ratio-based model (\ref{regressionmodel_weaker_ratio}), one gets $\varphi_1=\varphi_3=1.49$ and $\delta_2=\delta_4=1.33$ with $\gamma=1.12$. Reminding that $e_1=e_3=5$ and $e_2=e_4=2$, this means that five and two blank cards correspond to a ratio 1.33 and 1.12, respectively. Solving the ratio-based regression model (\ref{regressionmodel_weaker_ratio}) one gets a null total sum of the errors ($\overline{\sigma}=0$) and the following weights of criteria $w_1= 0.21, w_2=0.61, w_3=0.18$, so that, the overall evaluations of regions are obtained as listed in Table  \ref{Tab:regressionmodel_WS_weaker_ratio}. $\blacksquare$
\end{example}
\end{enumerate}
In the following, we will refer to the standard regression problem (\ref{regressionmodel}) based on the DCM. Of course, any extension we propose in the next sections can be generalized to the other regression models discussed above.

%%%%%%%%%%%%%%%%%%%%%%%%%%%%%%%
\section{A more robust and richer DOR}%%%
\label{sec:exentions}
%%%%%%%%%%%%%%%%%%%%%%%%%%%%%%%
In this section, we shall extend the DOR approach described in the previous section aiming to provide more robust recommendations on the considered problem and to take into account hierarchical structures of criteria. 
%%%%%%%%%%%%%%%%%%%%%%%%%%%%%%%%%%%%%%%%%%%%
\subsection{Robust Ordinal Regression}
\label{ROR}
%%%%%%%%%%%%%%%%%%%%%%%%%%%%%%%%%%%%%%%%%%%%
Solving the LP problems (\ref{regressionmodel}) and (\ref{CompatibleMaximizingk}) one finds a value function $\overline{U}$ that is compatible or deviates as little as possible from the DM's preferences. However, taking into account robustness concerns, it is reasonable to consider other value functions close to $\overline{U}$ obtained by some perturbation of its parameters. In this perspective, we shall consider \textit{compatible with the DM's preferences} any value function $U$ satisfying the constraints $E^{DM{'}}$ (with $k>0$), so that, it is reasonable taking into account all of them. With this aim we apply the ROR \citep{GrecoMousseauSlowinski2008,CorrenteEtAl2013},  defining a necessary and a possible preference relation on the set of alternatives $A$. Given $a_{j_1},a_{j_2}\in A$, on the one hand, $a_{j_1}$ is necessarily (weakly) preferred to $a_{j_2}$, denoted by $a_{j_1}\succsim^{N} a_{j_2}$, if $a_{j_1}$ is at least as good as $a_{j_2}$ for all compatible value functions, while, on the other hand, $a_{j_1}$ is possibly (weakly) preferred to $a_{j_2}$, denoted by $a_{j_1}\succsim^{P} a_{j_2}$, if $a_{j_1}$ is at least as good as $a_{j_2}$ for at least one compatible value function. \\
To compute the necessary and possible preference relations, the following LP problems need to be solved for each pair of alternatives $(a_{j_1},a_{j_2})\in A\times A$:\\
\begin{minipage}{0.45\textwidth}
\begin{equation}
\label{possiblepreferences}
\begin{array}{l}
\;\;\varepsilon^{P}(a_{j_1},a_{j_2})=\max\varepsilon, \;\;\mbox{subject to} \\[2mm]
\left.
\begin{array}{l}
U(a_{j_1})\geqslant U(a_{j_2}), \\[2mm]
E^{DM{''}} \\[2mm]
\end{array}
\right\}E^{P}(a_{j_1},a_{j_2})
\end{array}
\end{equation}
\end{minipage}
\hfill
\begin{minipage}{0.45\textwidth}
\begin{equation}
\label{necessarypreferences}
\begin{array}{l}
\;\;\varepsilon^{N}(a_{j_1},a_{j_2})=\max\varepsilon, \;\;\mbox{subject to} \\[2mm]
\left.
\begin{array}{l}
U(a_{j_2})\geqslant U(a_{j_1})+\varepsilon, \\[2mm]
E^{DM{''}}. \\[2mm]
\end{array}
\right\}E^{N}(a_{j_1},a_{j_2})
\end{array}
\end{equation}
\end{minipage}\\[2mm]
\noindent where $E^{DM^{''}}$ is obtained by $E^{DM^{'}}$ replacing $k\geqslant0$ with $k\geqslant \varepsilon$. In particular:
\begin{description}
\item[Case P1)] $E^P(a_{j_1},a_{j_2})$ is infeasible or $\varepsilon^{P}(a_{j_1},a_{j_2})\leqslant 0$: then, there is not any compatible value function for which $U(a_{j_1})\geqslant U(a_{j_2})$. Consequently, $not(a_{j_1}\succsim^{P}a_{j_2})$ (implying that $a_{j_2}\succsim^Na_{j_1}$),
\item[Case P2)] $E^P(a_{j_1},a_{j_2})$ is feasible and $\varepsilon^{P}(a_{j_1},a_{j_2})>0$: then, there is at least one compatible value function such that $U(a_{j_1})\geqslant U(a_{j_2})$ and, therefore, $a_{j_1}\succsim^P a_{j_2}$;
\item[Case N1)] $E^{N}(a_{j_1},a_{j_2})$ is infeasible or $\varepsilon^{N}(a_{j_1},a_{j_2})\leqslant 0$: then, there is not any compatible value function for which $U(a_{j_2})>U(a_{j_1})$. Therefore, $a_{j_1}\succsim^{N}a_{j_2}$,
\item[Case N2)] $E^{N}(a_{j_1},a_{j_2})$ is feasible and $\varepsilon^{N}(a_{j_1},a_{j_2})>0$, then, there is at least one compatible value function for which $U(a_{j_2})>U(a_{j_1})$. Therefore, $not(a_{j_1}\succsim^{N}a_{j_2})$ (implying that $a_{j_2}\succsim^{P}a_{j_1}$).
\end{description}

%%%%%%%%%%%%%%%%%%%%%%%%%%%%%%%%%%%%%%%%%%%%%%%%%%%%%%%%%%%%%%%
\subsection{Stochastic Multicriteria Acceptability Analysis}
\label{SMAASection}
%%%%%%%%%%%%%%%%%%%%%%%%%%%%%%%%%%%%%%%%%%%%%%%%%%%%%%%%%%%%%%%
Even if in a different way in comparison with the ROR, the SMAA \citep{LahdelmaHokkanenSalminen1998,PelissariEtAl2020} provides recommendations on the alternatives at hand considering the whole space of compatible value functions.\\
SMAA gives information in statistical terms based on a set of compatible value functions sampled from the simplex defined by constraints in $E^{DM{''}}$. Because all the above constraints are linear, one can efficiently perform the sampling using the Hit-And-Run (HAR) method \citep{Smith1984,VanValkenhoefTervonenPostmus2014}. In our context, for each sampled value function, a ranking of the alternatives in $A$ can be obtained. Based on these rankings, SMAA provides the following indices:
\begin{itemize}
\item \textit{Rank Acceptability Index} (RAI), $b^{v}(a)$: it is the frequency with which $a\in A$ is in place $v$, $v=1,\ldots,|A|,$ in the considered rankings,
\item \textit{Pairwise Winning Index} (PWI), $p(a_{j_1},a_{j_2})$: it is the frequency with which $a_{j_1}$ is preferred to $a_{j_2}$. 
\end{itemize}

Even if RAIs provide more robust information on the considered problem, in general, they do not give a total ranking of the alternatives under consideration. To overcome this problem, two different procedures can be used: 
\begin{itemize}
\item \textit{Computing the expected ranking of each alternative:} Following \cite{LahdelmaSalminen2001}, each $a\in A$ can be associated to a value $ER(a)$ being the weighted average of its RAIs. Formally, $
ER(a)~=~\displaystyle\sum_{v=1}^{|A|}v\cdot b^v\left(a\right).$ On the basis of $ER(a)$, one can obtain a complete ranking of the alternatives in $A$ so that, for each $a,b \in A$, $a$ is ranked not worse than $b$ if $ER(a) \leqslant ER(b)$;
\item \textit{Computing the ranking obtained by the barycenter of the compatible space:} As mentioned earlier, the SMAA indices are calculated by considering the rankings of the alternatives, which are determined by using each value function sampled from the space defined by constraints in $E^{DM''}$. The value function obtained by averaging the sampled value functions has been proved to be able to well represent the preferences of the DM \citep{ArcidiaconoCorrenteGreco2023}. Therefore, one can use this compatible value function, which represents an approximation of the barycenter of the space defined by constraints in $E^{DM''}$, to rank the alternatives at hand. 
\end{itemize}   

%%%%%%%%%%%%%%%%%%%%%%%%%%%%%%%%%%%%%%%%%%%%%%%%%%%%
\subsection{Multiple Criteria Hierarchy Process}\label{MCHPSection}
%%%%%%%%%%%%%%%%%%%%%%%%%%%%%%%%%%%%%%%%%%%%%%%%%%%%
In real-world decision making problems criteria are generally structured in a hierarchical way. It is possible to consider a root criterion (the objective of the problem), and some macro-criteria descending from it until the elementary criteria placed at the bottom of the hierarchy and on which the alternatives are evaluated. The Multiple Criteria Hierarchy Process (MCHP \citealt{CorrenteGrecoSlowinski2012}) can then be used to deal with such problems. It permits to take into account the preferences of the DM at both partial and global levels as well as providing recommendations globally and partially.\\
According to the MCHP framework, in the following, $g_{\mathbf{r}}$ will be a generic criterion in the hierarchy, while by $g_{\mathbf{r}}=g_{\mathbf{0}}$ we refer to the whole set of criteria; $I_G$ is the set of the indices of all criteria in the hierarchy; $EL\subseteq I_G$ is the set of the indices of elementary criteria, while, $E(g_{\mathbf{r}})\subseteq EL$ is the set of indices of the elementary criteria descending from $g_{\mathbf{r}}$; the value of an alternative $a$ on criterion $g_{\mathbf{r}}$ will be denoted by $U_{\mathbf{r}}(a)$, while the global value of $a$ will be denoted by $U_{\mathbf{0}}(a)$. The value of $a$ on macro-criterion $g_{\mathbf{r}}$, that is, $U_{\mathbf{r}}(a)$ will depend only on its performance on the elementary criteria descending from it. Of course, the definition of $U_{\mathbf{r}}$ will change according to the type of function used to represent the DM's preferences (see \ref{PreferenceModels}).

Extending our proposal described in Section \ref{originalmodel} to the MCHP framework, the DM is therefore asked (but they are not obliged) to apply the DCM for each macro-criterion $g_{\mathbf{r}}$ following these steps:  
\begin{enumerate}
\item Rank-ordering the alternatives from the less preferred to the most preferred with respect to criterion $g_{\mathbf{r}}$ in sets $L_{(\mathbf{r},1)},L_{(\mathbf{r},2)},\ldots,L_{(\mathbf{r},s(\mathbf{r}))}$. That is, alternatives in $L_{(\mathbf{r},h+1)}$ are preferred to alternatives in $L_{(\mathbf{r},h)}$ on $g_{\mathbf{r}}$ for all $h=1,\ldots,s(\mathbf{r})-1$ and alternatives in $L_{(\mathbf{r},h)}$ are indifferent on $g_{\mathbf{r}}$, for all $h=1,\ldots,s(\mathbf{r})$;
\item Putting a certain number of blank cards $e_{(\mathbf{r},h)}$ between sets $L_{(\mathbf{r},h)}$ and $L_{(\mathbf{r},h+1)}$ to increase the difference between the value on $g_\mathbf{r}$ of the alternatives in $L_{(\mathbf{r},h+1)}$ and the value on $g_\mathbf{r}$ of the alternatives in $L_{(\mathbf{r},h)}$;
\item Providing the number of blank cards $e_{(\mathbf{r},0)}$ between the ``fictitious zero alternative on $g_{\mathbf{r}}$" and the alternatives in $L_{(\mathbf{r},1)}$. In this case fictitious zero alternative on $g_{\mathbf{r}}$ is a fictitious alternative $a_0$ having a null value on $g_{\mathbf{r}}$ ($U_{\mathbf{r}}(a_0)=0$);
\item Each alternative $a\in A^R$ is assigned a value $\nu_{\mathbf{r}}(a)$ for each macro-criterion $g_{\mathbf{r}}$ such that if $a\in L_{(\mathbf{r},h+1)}$ and $a^{\prime}\in L_{(\mathbf{r},h)}$, $\nu_{\mathbf{r}}(a)=\nu_{\mathbf{r}}(a^{\prime})+(e_{(\mathbf{r},h)}+1)$, $h=1,\ldots,s(\mathbf{r})-1$ and if $a\in L_{(\mathbf{r},1)}$ then $\nu_{\mathbf{r}}(a)=e_{(\mathbf{r},0)}+1.$ Consequently, for $a\in L_{(\mathbf{r},h)}$, we have: $
\nu_{\mathbf{r}}(a)=\displaystyle\sum_{p=0}^{h-1}\left(e_{(\mathbf{r},p)}+1\right);$
\item The parameters of the value function $U$ are determined in a way that for all $a \in A^R$ and all $g_{\mathbf{r}}$, $U_{\mathbf{r}}(a)$ deviates as less as possible from $k_{\mathbf{r}}\cdot\nu_{\mathbf{r}}(a)$ with $k_{\mathbf{r}}$ a scalarizing constant that can be interpreted as the value of a single blank card w.r.t. macro-criterion $g_{\mathbf{r}}$. More precisely, for all $a \in A^R$ and all macro-criteria $g_{\mathbf{r}}$, one considers the positive and negative deviations $\sigma_{\mathbf{r}}^+(a)$ and $\sigma_{\mathbf{r}}^-(a)$ between $U_{\mathbf{r}}(a)$ and $k_{\mathbf{r}}\cdot\nu_{\mathbf{r}}(a)$, that is $U_{\mathbf{r}}(a)-\sigma_{\mathbf{r}}^+(a)+\sigma_{\mathbf{r}}^-(a)=k_{\mathbf{r}}\cdot\nu_{\mathbf{r}}(a)$.  
\end{enumerate}
The preference information the DM provides on macro-criterion $g_{\mathbf{r}}$ through the DCM can be translated into the following set of constraints
$$
\left.
\begin{array}{l}
U_{\mathbf{r}}(a) \geqslant U_{\mathbf{r}}(a^{\prime}), \;\mbox{for all}\; a,a^{\prime}\in A^R,\;\mbox{such that}\;\nu_{\mathbf{r}}(a)\geqslant\nu_{\mathbf{r}}(a^{\prime}), \\[2mm]
U_{\mathbf{r}}(a)-\sigma_{\mathbf{r}}^{+}(a) +\sigma_{\mathbf{r}}^{-}(a)=k_{\mathbf{r}}\cdot \nu_{\mathbf{r}}(a)\;\; \mbox{for all} \;a\in A^R,\\[2mm]
k_{\mathbf{r}} \geqslant 0,\\[2mm]
\sigma_{\mathbf{r}}^{+}(a) \geqslant 0,\; \sigma_{\mathbf{r}}^{-}(a) \geqslant 0 \mbox{ for all}\;a\in A^R, \\[2mm]
\end{array}
\right\}E_{\mathbf{r}}
$$
\noindent where: 
\begin{itemize}
\item $U_{\mathbf{r}}(a)$ is the value assigned by $U$ to alternative $a$ w.r.t. macro-criterion $g_{\mathbf{r}}$. The formalization of $U_{\mathbf{r}}$ depends on the type of value function $U$ used to approximate the preferences of the DM (see Section \ref{PreferenceModels}), 
\item $k_{\mathbf{r}}$ is the value of a blank card on $g_{\mathbf{r}}$, 
\item $\sigma_{\mathbf{r}}^{+}(a)$ and $\sigma_{\mathbf{r}}^{-}(a)$ are over and under estimations related to $a$ and $g_{\mathbf{r}}$. 
\end{itemize}
To check a function $U$ compatible with the DM's preferences and presenting the minimum error, one has to solve the following problem: 
$$
\begin{array}{l}
\;\;\overline{\sigma}_{MCHP}=\min\displaystyle\sum_{\mathbf{r}\in I_G\setminus EL}\sum_{a\in A^R}\left(\sigma^{+}_{\mathbf{r}}(a)+\sigma^{-}_{\mathbf{r}}(a)\right), \;\;\mbox{subject to},\\[8mm]
\left.
\begin{array}{l}
E^{Model},\\[2mm]
%U_{\mathbf{0}}\left(x_1^{m_1},\ldots,x_{n}^{m_n}\right)=1,\\[2mm]
\displaystyle\cup_{\mathbf{r}\in I_G\setminus EL}E_{\mathbf{r}}.
\end{array}\right\}E^{DM}_{MCHP}
\end{array}
$$
To ensure that the value of each blank card $k_{\mathbf{r}}$ is greater than zero one has to solve iteratively the following LP problem:
\begin{equation}
\label{EpsilonMCHP}
\begin{array}{l}
\;\;\varepsilon_{MCHP}^*=\max\varepsilon,\;\;\mbox{subject to},\\[2mm]
\left.
\begin{array}{l}
E^{DM^{'}}_{MCHP},\\[2mm]
\displaystyle\sum_{\mathbf{r}\in I_G\setminus EL}\sum_{a\in A^R}\left(\sigma^{+}_{\mathbf{r}}(a)+\sigma^{-}_{\mathbf{r}}(a)\right)\leqslant\overline{\sigma}_{MCHP}+\eta\left(\overline{\sigma}_{MCHP}\right),
\end{array}
\right\}E^{DM{''}}_{MCHP}
\end{array}
\end{equation}
\noindent where $E^{DM^{'}}_{MCHP}$ is obtained by $E^{DM}_{MCHP}$ replacing the constraints $k_{\mathbf{r}}\geqslant 0$ (one in each $E_{\mathbf{r}}$) with $k_{\mathbf{r}}\geqslant\varepsilon$ and $\eta\left(\overline{\sigma}_{MCHP}\right)$ is an admitted deterioration error with respect to the optimal value obtained in the previous step, that is, $\overline{\sigma}_{MCHP}$. At the beginning, $\eta\left(\overline{\sigma}_{MCHP}\right)=0$. However, if $\varepsilon_{MCHP}^*=0$, then, one has to increase $\eta\left(\overline{\sigma}_{MCHP}\right)$ (as suggested in Section \ref{originalmodel}) until $\varepsilon_{MCHP}^*>0.$ \\
Let us observe that, differently from LP problem (\ref{CompatibleMaximizingk}) where we have maximized directly the value of the unique blank card $k$, here, we have replaced $k_{\mathbf{r}}\geqslant 0$ with $k_{\mathbf{r}}\geqslant \varepsilon$ and, then, we have maximized $\varepsilon$ to ensure that if $\varepsilon^{*}_{MCHP}>0$ then all $k_{\mathbf{r}}$ are greater than zero. \\
To conclude this section, let us observe the following: 
\begin{itemize}
\item The DM is not obliged to provide the preference information on all macro-criteria in the hierarchy but only on those they are more confident. Moreover, some information can also be imprecisely given (see Section \ref{ImpreciseInformation}),
\item ROR can be applied in the MCHP context. In particular, to check if $a_{j_1}\succsim_{\mathbf{r}}^{P}a_{j_2}$ ($a_{j_1}$ is possibly preferred to $a_{j_2}$ on $g_{\mathbf{r}}$) or $a_{j_1}\succsim_{\mathbf{r}}^{N}a_{j_2}$ ($a_{j_1}$ is necessarily preferred to $a_{j_2}$ on $g_{\mathbf{r}}$), one has to solve the following two LP problems: \\
\begin{minipage}{0.45\textwidth}
\begin{equation}
\label{PossibleMCHP}
\begin{array}{l}
\;\;\varepsilon_{\mathbf{r}}^{P}(a_{j_1},a_{j_2})=\max\varepsilon,\;\;\mbox{subject to} \\[2mm]
\left.
\begin{array}{l}
U_{\mathbf{r}}(a_{j_1})\geqslant U_{\mathbf{r}}(a_{j_2}), \\[2mm]
E_{MCHP}^{DM^{''}}. \\[2mm]
\end{array}
\right\}E_{\mathbf{r}}^{P}(a_{j_1},a_{j_2})
\end{array}
\end{equation}
\end{minipage}
\hfill
\begin{minipage}{0.45\textwidth}
\begin{equation}
\label{NecessaryMCHP}
\begin{array}{l}
\;\;\varepsilon_{\mathbf{r}}^{N}(a_{j_1},a_{j_2})=\max\varepsilon, \;\;\mbox{subject to} \\[2mm]
\left.
\begin{array}{l}
U_{\mathbf{r}}(a_{j_2})\geqslant U_{\mathbf{r}}(a_{j_1})+\varepsilon, \\[2mm]
E_{MCHP}^{DM^{''}}. \\[2mm]
\end{array}
\right\}E_{\mathbf{r}}^{N}(a_{j_1},a_{j_2})
\end{array}
\end{equation}
\end{minipage}\\[2mm]
As commented in Section \ref{ROR},
\begin{itemize}
\item $a_{j_1}\succsim^{P}_{\mathbf{r}}a_{j_2}$ iff $E_{\mathbf{r}}^{P}(a_{j_1},a_{j_2})$ is feasible and $\varepsilon_{\mathbf{r}}^{P}\left(a_{j_1},a_{j_2}\right)>0$;
\item $a_{j_1}\succsim^{N}_{\mathbf{r}}a_{j_2}$ iff $E_{\mathbf{r}}^{N}(a_{j_1},a_{j_2})$ is infeasible or $\varepsilon_{\mathbf{r}}^{N}\left(a_{j_1},a_{j_2}\right)\leqslant 0$;
\end{itemize}
\item SMAA can be applied in the MCHP context on the basis of a sampling of value functions from the simplex defined by constraints in $E^{DM^{''}}_{MCHP}$.
For each sampled value function, a ranking of the alternatives at hand can be done on each macro-criterion $g_{\mathbf{r}}$, $\mathbf{r}\in I_G\setminus EL$. Consequently, the RAIs and PWIs defined in Section \ref{SMAASection} can be computed for each macro-criterion $g_{\mathbf{r}}$:
\begin{itemize}
\item $b^{v}_{\mathbf{r}}(a)$: it is the frequency with which alternative $a\in A$ is in place $v$, $v=1,\ldots,|A|$, in the rankings produced on $g_{\mathbf{r}}$,
\item $p_{\mathbf{r}}(a_{j_1},a_{j_2})$: it is the frequency with which $a_{j_1}$ is preferred to $a_{j_2}$ on $g_{\mathbf{r}}$. 
\end{itemize}
Analogously, the expected ranking of each alternative $a$ on $g_{\mathbf{r}}$ ($ER_{\mathbf{r}}(a_j)$) as well as the ranking obtained using the approximation of the barycenter of the space defined by constraints in $E^{DM{''}}_{MCHP}$ can be computed. Formally, for each $g_{\mathbf{r}}$, $\mathbf{r}\in I_G\setminus EL$, $
ER_{\mathbf{r}}(a)=\displaystyle\sum_{v=1}^{|A|}v\cdot b_{\mathbf{r}}^{v}(a). 
$
\end{itemize}

%%%%%%%%%%%%%%%%%%%%%%%%%%%%%%%%%%%%%%%%%%%%
\section{Providing imprecise information}%%%
\label{ImpreciseInformation} %%%%%%%%%%%%%%%%%%%%%%%%%%%%%
%%%%%%%%%%%%%%%%%%%%%%%%%%%%%%%%%%%%%%%%%%%%
In this section, we shall extend the proposal described in Section \ref{originalmodel}. We shall consider the case in which the DM is not able to precisely provide the number of blank cards between two successive subsets of alternatives $L_{h}$ and $L_{h+1}$, with $h=1,\ldots,s-1.$  

%%%%%%%%%%%%%%%%%%%%%%%%%%%%%%%%%%%%
\subsection{Interval information}%%%
\label{IntervalInformation}
%%%%%%%%%%%%%%%%%%%%%%%%%%%%%%%%%%%%
Let us assume that the DM can define the minimum ($e_h^L$) and the maximum ($e_h^U$) number of blank cards to be included between sets $L_h$ and $L_{h+1}$. Formally, this means that the number of blank cards $e_h$ between the two subsets of alternatives is such that $e_h \in [e_h^L, e_h^U]$ for all $h = 0, \ldots, s-1$. %Denoting, again, by $a_j$ a generic alternative belonging to $L_j$
To infer a value function $U$ able to represent the DM's preference information expressed by the use of the ``imprecise" DCM, one has to solve the following LP problem:
\begin{equation*}
\label{zetamodel}
\begin{array}{l}
\;\;\;\;\overline{\sigma}_{Interval}=\min \displaystyle \sum_{a\in A^R} (\sigma^{+}(a)+\sigma^{-}(a)),\;\; \mbox{subject to}\\[3mm]
\left.
\begin{array}{l}
\;\;E^{Model}\\[2mm]
\;\;U(a) \geqslant U(a^{\prime}), \;\mbox{for all}\; a\in L_h, a^{\prime} \in L_{h'},\;\mbox{such that}\; h\geqslant h', h,h' = 1,\ldots,s \\[2mm]
\;\;U(a)-\sigma^{+}(a) +\sigma^{-}(a)=\widehat{\nu}(a)\; \mbox{for all}\;a\in A^R\\[2mm]
\;\;k \geqslant 0\\[2mm]
\;\;\sigma^{+}(a) \geqslant 0,\; \sigma^{-}(a) \geqslant 0 \;\mbox{for all}\; a\in A^R\\[2mm]
\left.
\begin{array}{l}
\widehat{\nu}(a) \geqslant \widehat{\nu}(a^{\prime}) + (e^{L}_{h} + 1)\cdot k, \\[2mm]
\widehat{\nu}(a) \leqslant \widehat{\nu}(a^{\prime}) + (e^{U}_{h} + 1)\cdot k, \\[2mm]
\end{array}
\right\} \;\mbox{for all}\; a\in L_{h+1}, a^{\prime} \in L_h, \;h = 0,\ldots,s-1,\\[3mm]
\end{array}
\right\}E^{DM}_{Interval}
\end{array}
\end{equation*}
\noindent where:
\begin{itemize}
\item $E^{Model}$ is the set of monotonicity and normalization constraints related to the considered preference function $U$ (see \ref{PreferenceModels});
\item $U(a) \geqslant U(a^{\prime})$ imposes that, following the ranking given by the DM, if $a\in L_h$ and $a^{\prime}\in L_{h^{\prime}}$ such that $h\geqslant h^{\prime}$, the value of $a$ is not lower than the value of $a^{\prime}$. Observe that this constraint implies that $U\left(a\right)=U\left(a^{\prime}\right)$ in case both of them belong to the same level $L_h$, with $h=1,\ldots,s$;
%\item $U\left(a_{j_1}\right)=U\left(a_{j_2}\right)$ imposes that alternatives $a_{j_1}$ and $a_{j_2}$ have the same value since they belong to the same subset of alternatives ($L_j$);
\item $U(a)-\sigma^{+}(a) +\sigma^{-}(a)=\widehat{\nu}(a)$ imposes that the value of $a$ (that is $U(a)$) is equal to $\widehat{\nu}(a)$ being implicitly given by the product between the value of a blank card $(k)$ and the number of levels up to alternative $a$ ($\nu(a)$), that is $\widehat{\nu}(a)=k\cdot \nu(a)$;
\item $k\geqslant 0$ is the value of a single blank card;
\item $\sigma^{+}(a) \geqslant 0, \sigma^{-}(a) \geqslant 0$ impose that the slack variables, representing overestimation and underestimation of $U(a)$, respectively, are not negative;
\item constraints $\widehat{\nu}(a) \geqslant \widehat{\nu}(a^{\prime}) + (e^{L}_{h} + 1)\cdot k$ and $\widehat{\nu}(a) \leqslant \widehat{\nu}(a^{\prime}) + (e^{U}_{h} + 1)\cdot k$ link the values assigned by the DCM to alternatives $a$ and $a^{\prime}$ taking into account the interval $[e^L_h,\;e^U_h]$ of the possible values of the blank card $e_h$ put between two contiguous levels $L_h$ and $L_{h+1}$, $h=0,\ldots,s-1$. Indeed, 
$$
e_h^{L}\leqslant e_h\leqslant e_h^U \Rightarrow \left(e_h^{L}+1\right)\leqslant \left(e_h+1\right)\leqslant \left(e_h^U+1\right) \Rightarrow \left(e_h^{L}+1\right)\cdot k\leqslant \left(e_h+1\right)\cdot k\leqslant \left(e_h^U+1\right)\cdot k\Rightarrow
$$
$$
\Rightarrow \widehat{\nu}\left(a^{\prime}\right)+\left(e_h^{L}+1\right)\cdot k\leqslant \underbrace{\widehat{\nu}\left(a^{\prime}\right)+\left(e_h+1\right)\cdot k}_{\widehat{\nu}\left(a\right)} \leqslant \widehat{\nu}\left(a^{\prime}\right)+\left(e_h^U+1\right)\cdot k.
$$
\end{itemize}
Following the description given in Section \ref{originalmodel}, the value of the blank card $k$ has to be greater than zero. To maximize it, one has to solve the following LP problem
\begin{equation}\label{IntervalLP}
\begin{array}{l}
\;\;\varepsilon_{Interval}^{*}=\max\varepsilon,\;\mbox{subject to},\\[2mm]
\left.
\begin{array}{l}
E^{DM'}_{Interval},\\[1mm]
\displaystyle \sum_{a \in A^R} (\sigma^{+}(a)+\sigma^{-}(a))\leqslant\overline{\sigma}_{Interval}+\eta\left(\overline{\sigma}_{Interval}\right)
\end{array}
\right\}
\end{array}
\end{equation}
\noindent where $E^{DM'}_{Interval}$ is obtained by $E^{DM}_{Interval}$ replacing the constraint $k\geqslant 0$ with $k\geqslant\varepsilon$ and $\eta\left(\overline{\sigma}_{Interval}\right)$ is an admitted deterioration error with respect to the optimal value obtained in the previous step, that is, $\overline{\sigma}_{Interval}$. At the beginning, $\eta\left(\overline{\sigma}_{Interval}\right)=0$. However, if $\varepsilon_{Interval}^*=0$, then, one has to increase $\eta\left(\overline{\sigma}_{Interval}\right)$ (as suggested in Section \ref{originalmodel}) until $\varepsilon_{Interval}^*>0.$ Notice that from a computational point of view, comparing the LP problems (\ref{IntervalLP}) and (\ref{CompatibleMaximizingk}), maximizing $\varepsilon$ in (\ref{IntervalLP}) is equivalent to maximize $k$ in (\ref{CompatibleMaximizingk}).

%%%%%%%%%%%%%%%%%%%%%%%%%%%%%%%%%%%%%%%%%%%%%%%%
\subsection{Imprecise or missing information}%%%
\label{ImprecisionSection}%%%%%%%%%%%%%%%%%%%%%%
%%%%%%%%%%%%%%%%%%%%%%%%%%%%%%%%%%%%%%%%%%%%%%%%
In some situations, the DM could not be able to specify the lower or the upper bounds of $e_h$, that is, $e_h^L$ or $e_h^U$. 
Summarizing, considering the number of blank cards $e_h$ between levels $L_h$ and $L_{h+1}$, with $h=1,\ldots,s-1$, five different cases can be observed and each of them is translated into a few linear equalities or inequalities:
\begin{description}
\item[A)] The DM can precisely assign a number to $e_h$, that is, $e_h\in\nat^{*}$. In this case, 
$$
\widehat{\nu}\left(a\right)=\widehat{\nu}(a^{\prime})+(e_h+1)\cdot k,\;\mbox{for all}\;a\in L_{h+1},a^{\prime}\in L_h.
$$
\item[B)] The DM can define the lower $e_h^L$ and the upper $e_h^U$ bound of $e_h$, that is, $e_{h}\in\left[e_h^{L},e_h^U\right]$. In this case:
$$
\left.
\begin{array}{l}
\widehat{\nu}(a)\geqslant\widehat{\nu}(a^{\prime}) + \left(e^{L}_{h} + 1\right)\cdot k, \\[2mm]
\widehat{\nu}(a)\leqslant\widehat{\nu}(a^{\prime}) + \left(e^{U}_{h} + 1\right)\cdot k 
\end{array}
\right\}\;\mbox{for all}\;a\in L_{h+1},a^{\prime}\in L_h.
$$
\item[C)] The DM can define only the lower bound $e_h^L$ of $e_h$, that is, $e_h\in\left[e_h^{L},?\right]$. In this case:
\begin{equation*}
\widehat{\nu}(a) \geqslant \widehat{\nu}(a^{\prime}) + \left(e^{L}_{h} + 1\right)\cdot k,\;\mbox{for all}\;a\in L_{h+1},a^{\prime}\in L_h. \\[1mm]
\end{equation*}
\item[D)] The DM can define only the upper bound $e_h^U$ of $e_h$, that is, $e_h\in\left[?,e_h^U\right]$. In this case:
$$
\left.
\begin{array}{l}
\widehat{\nu}(a) \geqslant \widehat{\nu}(a^{\prime})+k, \\[1mm]
\widehat{\nu}(a) \leqslant \widehat{\nu}(a^{\prime}) + \left(e^{U}_{h} + 1\right)\cdot k \\[2mm]
\end{array}
\right\}\;\mbox{for all}\;a\in L_{h+1},a^{\prime}\in L_h.
$$
\item[E)] The DM can define neither the lower bound $e_h^L$ nor the upper bound $e_h^U$ of $e_h$, that is, $e_h\in\left[?,?\right]$. In this case:
\begin{equation*}
\begin{array}{l}
\widehat{\nu}(a) \geqslant \widehat{\nu}(a^{\prime})+k,\;\mbox{for all}\;a\in L_{h+1},a^{\prime}\in L_h. \\[1mm]
\end{array}
\end{equation*}
\end{description}
Taking into account all these different cases, to check for a value function compatible with the possible imprecise preference information provided by the DM, one has to solve the following LP problem: 

\begin{equation*}
\label{zetamodelgeneral}
\begin{array}{l}
\;\;\;\;\;\;\overline{\sigma}_{Imprecise}=\min \displaystyle \sum_{a\in A^R} (\sigma^{+}(a)+\sigma^{-}(a)), \;\; \mbox{subject to}\\[7mm]
\left.
\begin{array}{l}
\;\;\;\;E^{Model}\\[1mm]
\;\;\;\;U(a) \geqslant U(a^{\prime}), \;\mbox{with}\; a\in L_h, a^{\prime}\in L_{h'},\;\mbox{such that}\; h\geqslant h^{\prime},\;h,h^{\prime} = 1, \ldots, s-1, \\[2mm]
\;\;\;\;U(a)-\sigma^{+}(a) +\sigma^{-}(a)=\widehat{\nu}(a)\;\; \mbox{for all}\; a \in A^R,\\[2mm]
%\;\;\;\;U(a_{j_1})=U(a_{j_2}), \;\;\mbox{for all}\;\;a_{j_1},a_{j_2}\in L_j,\;\;\mbox{and for all}\;\;j=1,\ldots,s,\\[2mm]
\;\;\;\;k \geqslant 0,\\[2mm]
\;\;\;\;\sigma^{+}(a) \geqslant 0, \;\sigma^{-}(a) \geqslant 0 \;\mbox{for all}\; a\in A^R, \\[2mm]
\left.
\begin{array}{lll}
\;\;\widehat{\nu}\left(a\right)=\widehat{\nu}(a^{\prime})+(e_h+1)\cdot k &\mbox{if}& e_h\in\nat^{*},\\[2mm]
\left.
\begin{array}{l}
\widehat{\nu}(a)\geqslant\widehat{\nu}(a^{\prime}) + \left(e^{L}_{h} + 1\right)\cdot k \\[3mm]
\widehat{\nu}(a)\leqslant\widehat{\nu}(a^{\prime}) + \left(e^{U}_{h} + 1\right)\cdot k 
\end{array}
\right\} &\mbox{if}&e_h\in\left[e_h^{L},e_h^{U}\right] \\[6mm]
\;\;\widehat{\nu}(a) \geqslant \widehat{\nu}(a^{\prime}) + \left(e^{L}_{h} + 1\right)\cdot k &\mbox{if}&e_h\in\left[e_h^L,?\right]\\[3mm]
\left.
\begin{array}{l}
\widehat{\nu}(a) \geqslant \widehat{\nu}(a^{\prime})+k \\[1mm]
\widehat{\nu}(a) \leqslant \widehat{\nu}(a^{\prime}) + \left(e^{U}_{h} + 1\right)\cdot k \\[2mm]
\end{array}
\right\}&\mbox{if}&e_h\in\left[?,e_h^U\right]\\[6mm]
\;\;\widehat{\nu}(a) \geqslant \widehat{\nu}(a^{\prime})+k & \mbox{if} & e_h\in\left[?,?\right]\\[3mm]
\end{array}
\right\}\mbox{for all}\;a\in L_{h+1},a^{\prime}\in L_h.\\[2mm]
\end{array}
\right\}E^{DM}_{Imprecise}
\end{array}
\end{equation*}
To ensure that the value of the blank card $k$ is greater than zero, one has to solve the following LP problem
\begin{equation*}
\label{EpsilonImprecise}
\begin{array}{l}
\;\;\varepsilon^*_{Imprecise}=\max\varepsilon,\;\mbox{subject to},\\[2mm]
\left.
\begin{array}{l}
E^{DM'}_{Imprecise},\\[1mm]
\displaystyle \sum_{a\in A^R} \left(\sigma^{+}(a)+\sigma^{-}(a)\right)\leqslant\overline{\sigma}_{Imprecise}+\eta\left(\overline{\sigma}_{Imprecise}\right)
\end{array}
\right\}E^{DM{''}}_{Imprecise}
\end{array}
\end{equation*}
\noindent where $E^{DM'}_{Imprecise}$ is obtained by $E^{DM}_{Imprecise}$ replacing the constraint $k\geqslant 0$ with $k\geqslant\varepsilon$ and $\eta\left(\overline{\sigma}_{Imprecise}\right)$ is an admitted deterioration error with respect to the optimal value obtained in the previous step, that is, $\overline{\sigma}_{Imprecise}$. At the beginning, $\eta\left(\overline{\sigma}_{Imprecise}\right)=0$. However, if $\varepsilon_{Imprecise}^*=0$, then, one has to increase $\eta\left(\overline{\sigma}_{Imprecise}\right)$ (as suggested in Section \ref{originalmodel}) until $\varepsilon_{Imprecise}^*>0.$   

%%%%%%%%%%%%%%%%%%%%%%%%%%%%%%%%%%%%%%%%%%%%%%%%
\subsection{ROR, SMAA and hierarchy of criteria in presence of imprecise information}\label{RORSMAAMCHPMissingInformationSection}%%%
%%%%%%%%%%%%%%%%%%%%%%%%%%%%%%%%%%%%%%%%%%%%%%%%
The imprecision information described in the previous section can also be taken into account in the ROR and SMAA as well as in case the problem presents criteria structured in a hierarchical way and dealt with by the MCHP:
\begin{itemize}
\item Regarding the ROR and, for each $(a_{j_1},a_{j_2})\in A\times A$, to check if $a_{j_1}\succsim^{P}a_{j_2}$ and if $a_{j_1}\succsim^N a_{j_2}$ one has to solve the LP problems (\ref{possiblepreferences}) and (\ref{necessarypreferences}), respectively, replacing $E^{DM''}$ with $E^{DM''}_{Imprecise}$,
\item Regarding the SMAA, the sampling of compatible value functions has to be done from the simplex defined by constraints in $E^{DM{''}}_{Imprecise}$,
\item Regarding the MCHP, reminding that for each macro-criterion $g_{\mathbf{r}}$, $\mathbf{r}\in I_G\setminus EL$, $e_{(\mathbf{r},h)}$ is the number of blank cards between the sets of reference alternatives $L_{(\mathbf{r},h)}$ and $L_{(\mathbf{r},h+1)}$ and that the number of subsets in which reference alternatives are ordered w.r.t. this criterion is denoted by $s({\mathbf{r}})$ (see Section \ref{MCHPSection}), we assume that the DM can provide information on the lower $e_{(\mathbf{r},h)}^{L}$ or upper $e_{(\mathbf{r},h)}^{U}$ bounds of $e_{(\mathbf{r},h)}$. Of course, all the imprecise information presented in Section \ref{ImprecisionSection} can also be taken into account in the MCHP with the corresponding constraints translating them. For each macro-criterion $g_{\mathbf{r}}$, one needs to define the following set of constraints 
\begin{small}
$$
\left.
\begin{array}{l}
\;\;\;\;U_{\mathbf{r}}(a) \geqslant U_{\mathbf{r}}(a^{\prime}), \;\mbox{for all}\; a\in L_h,a^{\prime}\in L_{h^{\prime}}\;\mbox{such that}\; h\geqslant h^{\prime},\;h,h^{\prime}=1,\ldots,s, \\[2mm]
\;\;\;\;U_{\mathbf{r}}(a)-\sigma_{\mathbf{r}}^{+}(a) +\sigma_{\mathbf{r}}^{-}(a)=\widehat{\nu}_{\mathbf{r}}(a)\;\; \mbox{for all}\;a\in A^R,\\[2mm]
%\;\;\;\;U_{\mathbf{r}}(a_{j_1})=U_{\mathbf{r}}(a_{j_2}), \;\;\mbox{for all}\;\;a_{j_1},a_{j_2}\in L_{(\mathbf{r},j)},\;\;\mbox{and for all}\;\;j=1,\ldots,s({\mathbf{r}}),\\[2mm]
\;\;\;\;k_{\mathbf{r}} \geqslant 0,\\[2mm]
\;\;\;\;\sigma_{\mathbf{r}}^{+}(a) \geqslant 0, \;\sigma_{\mathbf{r}}^{-}(a) \geqslant 0\; \mbox{for all}\; a\in A^R, \\[2mm]
\left.
\begin{array}{lll}
\;\;\widehat{\nu}_{\mathbf{r}}\left(a\right)=\widehat{\nu}_{\mathbf{r}}(a^{\prime})+\left(e_{(\mathbf{r},h)}+1\right)\cdot k_{\mathbf{r}}, &\mbox{if}& e_{(\mathbf{r},h)}\in\nat^{*},\\[2mm]
\left.
\begin{array}{l}
\widehat{\nu}_{\mathbf{r}}(a)\geqslant\widehat{\nu}_{\mathbf{r}}(a^{\prime}) + \left(e^{L}_{(\mathbf{r},h)} + 1\right)\cdot k_{\mathbf{r}} \\[3mm]
\widehat{\nu}_{\mathbf{r}}(a)\leqslant\widehat{\nu}_{\mathbf{r}}(a^{\prime}) + \left(e^{U}_{(\mathbf{r},h)} + 1\right)\cdot k_{\mathbf{r}} 
\end{array}
\right\} &\mbox{if}&e_{(\mathbf{r},h)}\in\left[e_{(\mathbf{r},h)}^{L},e_{(\mathbf{r},h)}^{U}\right] \\[6mm]
\;\;\widehat{\nu}_{\mathbf{r}}(a) \geqslant \widehat{\nu}_{\mathbf{r}}(a^{\prime}) + \left(e^{L}_{(\mathbf{r},h)} + 1\right)\cdot k_{\mathbf{r}} & \mbox{if}&e_{(\mathbf{r},h)}\in\left[e_{(\mathbf{r},h)}^L,?\right]\\[3mm]
\left.
\begin{array}{l}
\widehat{\nu}_{\mathbf{r}}(a) \geqslant \widehat{\nu}_{\mathbf{r}}(a^{\prime})+k_{\mathbf{r}}, \\[1mm]
\widehat{\nu}_{\mathbf{r}}(a) \leqslant \widehat{\nu}_{\mathbf{r}}(a^{\prime}) + \left(e^{U}_{(\mathbf{r},h)} + 1\right)\cdot k_{\mathbf{r}} \\[2mm]
\end{array}
\right\}&\mbox{if}&e_{(\mathbf{r},h)}\in\left[?,e_{(\mathbf{r},h)}^U\right]\\[6mm]
\;\;\widehat{\nu}_{\mathbf{r}}(a) \geqslant \widehat{\nu}_{\mathbf{r}}(a^{\prime})+k_{\mathbf{r}} & \mbox{if} & e_{(\mathbf{r},h)}\in\left[?,?\right]\\[3mm]
\end{array}
\right\}\mbox{for all}\;a,a^{\prime}\in A^R.\\[2mm]
\end{array}
\right\}E^{Imprecise}_{\mathbf{r}}
$$
\end{small}
The existence of a value function compatible with the imprecise information provided by the DM is checked solving the following LP problem: 
\begin{equation}\label{ImpreciseMCHP}
\begin{array}{l}
\;\;\overline{\sigma}_{MCHP}^{Imprecise}=\min\displaystyle\sum_{\mathbf{r}\in I_G\setminus EL}\sum_{a\in A^R}\left(\sigma^{+}_{\mathbf{r}}(a)+\sigma^{-}_{\mathbf{r}}(a)\right) \;\;\mbox{subject to},\\[4mm]
\left.
\begin{array}{l}
E^{Model},\\[2mm]
% U_{\mathbf{0}}\left(x_1^{m_1},\ldots,x_{n}^{m_n}\right)=1,\\[2mm]
\displaystyle\cup_{\mathbf{r}\in I_G\setminus EL}E_{\mathbf{r}}^{Imprecise}.
\end{array}\right\}E^{Imprecise}_{MCHP}
\end{array}
\end{equation}
To ensure that the value of each blank card $k_{\mathbf{r}}$ is greater than zero, one has to solve the following LP problem
\begin{equation}
\label{EpsilonMCHPImprecise}
\begin{array}{l}
\;\;\varepsilon_{MCHP}^{*,Imprecise}=\max\varepsilon,\;\;\mbox{subject to},\\[2mm]
\left.
\begin{array}{l}
E^{Imprecise^{'}}_{MCHP},\\[2mm]
\displaystyle\sum_{\mathbf{r}\in I_G\setminus EL}\sum_{a\in A^R}\left(\sigma^{+}_{\mathbf{r}}(a)+\sigma^{-}_{\mathbf{r}}(a)\right)\leqslant\overline{\sigma}_{MCHP}^{Imprecise}+\eta\left(\overline{\sigma}_{MCHP}^{Imprecise}\right)
\end{array}
\right\}E^{Imprecise^{''}}_{MCHP}
\end{array}
\end{equation}
\noindent where $E^{Imprecise^{'}}_{MCHP}$ is obtained by $E^{Imprecise}_{MCHP}$ replacing $k_{\mathbf{r}}\geqslant 0$ with $k_{\mathbf{r}}\geqslant\varepsilon$ in each set of constraints $E_{\mathbf{r}}^{Imprecise}$ and $\eta\left(\overline{\sigma}^{Imprecise}_{MCHP}\right)$ is an admitted deterioration error with respect to the optimal value obtained in the previous step, that is, $\overline{\sigma}^{Imprecise}_{MCHP}$. At the beginning, $\eta\left(\overline{\sigma}_{MCHP}^{Imprecise}\right)=0$. However, if $\varepsilon_{MCHP}^{*,Imprecise}=0$, then, one has to increase $\eta\left(\overline{\sigma}_{MCHP}^{Imprecise}\right)$ (as suggested in Section \ref{originalmodel}) until $\varepsilon_{MCHP}^{*,Imprecise}>~0.$ \\
To apply ROR and SMAA in this context is therefore necessary to proceed in the following way. For each macro-criterion $g_{\mathbf{r}}$, $\mathbf{r}\in I_G\setminus EL$,   
\begin{itemize}
\item for each $(a_{j_1},a_{j_2})\in A\times A$, to check if $a_{j_1}\succsim_{\mathbf{r}}^{P}a_{j_2}$ and if $a_{j_1}\succsim_{\mathbf{r}}^{N}a_{j_2}$ one has to solve the LP problems (\ref{PossibleMCHP}) and (\ref{NecessaryMCHP}), respectively, replacing $E^{DM^{''}}_{MCHP}$ with $E^{Imprecise^{''}}_{MCHP}$;
\item the sampling on which the SMAA indices computation is based, has to be performed from the simplex defined by constraints in $E^{Imprecise^{''}}_{MCHP}$.
\end{itemize}
\end{itemize}

%%%%%%%%%%%%%%%%%%%%%%%%%%%%%%%%%%%%%%%%%%%%%%%%%%%%%%%%%%%%%%%%%%%%%%%%%%%%%%%%%%
\section{Didactic Example}\label{ExampleSection}	%%%%%%%%%%%%%%%%%%%%%%%%%%%%%%
%%%%%%%%%%%%%%%%%%%%%%%%%%%%%%%%%%%%%%%%%%%%%%%%%%%%%%%%%%%%%%%%%%%%%%%%%%%%%%%%%%
In this section, we shall show how the proposed methodology could be applied in a real-world problem related to a large ongoing research project (\href{https://grins.it/}{https://grins.it/}; \href{https://dse.unibo.it/en/university-outreach/next-generation-eu}{https://dse.unibo.it/en/university-outreach/next-generation-eu}). Let us assume that the DM wants to rank Italian Regions according to three main macro-criteria, being Circular Economy $(g_{\mathbf{1}})$, Innovation-Driven Development $(g_{\mathbf{2}})$ and Smart Specialization Strategies $(g_{\mathbf{3}})$. Each of these macro-criteria is then articulated in several elementary criteria which are shown in
Table \ref{alldata}. All elementary criteria have an increasing direction of preference ($+$) except ``Urban Waste Generation" which has a decreasing direction of preference $(-)$. Before applying our proposal, the data are normalized using the procedure described in \cite{GrecoEtAl2019a}. Normalized data are shown in Table \ref{data_normalized}.\\
\begin{sidewaystable} %\setlength{\tabcolsep}{0.9pt} 
\footnotesize
\caption{Regions' performances on elementary criteria divided in three Macro-Criteria: Circular Economy ($g_1$), Innovation-Driven Development ($g_2$), Smart Specialization ($g_3$)}
\label{alldata}
\begin{center} 
\begin{tabular}{m{0.01\textwidth}m{0.06\textwidth}|m{0.06\textwidth}m{0.08\textwidth}m{0.06\textwidth}|m{0.09\textwidth}m{0.09\textwidth}m{0.06\textwidth}|m{0.06\textwidth}m{0.09\textwidth}m{0.09\textwidth}}
\toprule
&&\multicolumn{3}{c}{\textbf{Circular Economy - $g_{(\mathbf{1})}$}} & \multicolumn{3}{c}{\textbf{Innovation-Driven Development - $g_{(\mathbf{2})}$}}&\multicolumn{3}{c}{\textbf{Smart Specialization - $g_{(\mathbf{3})}$}}\\
\cmidrule(lr){3-11}
\textbf{ID}&\textbf{Region}&Urban Waste Generation (kg/capita) (2017) $(g_{(\mathbf{1,1})})$ &Differentiated waste collection (\%) (2017) $(g_{(\mathbf{1,2})})$ &Renewable electricity (percentage) (2017) $(g_{(\mathbf{1,3})})$ &Research and Development Personnel (full-time equivalents per thousand inhabitants) (2017) $(g_{(\mathbf{2,1})})$ &Patents registered at the European Patent Office (EPO) (number per million inhabitants) (2012) $(g_{(\mathbf{2,2})})$ &Number of companies with ISO 14001 certification (2017) 
$(g_{(\mathbf{2,3})})$ &Percentage of Smart specialized companies (2018) $(g_{(\mathbf{3,1})})$ &Percentage of value added by Smart specialised companies to the total value added in the region (2018) $(g_{(\mathbf{3,2})})$ &Number of Smart specialized companies with high intensity of investments in social and environmental responsibility (2018) $(g_{(\mathbf{3,3})})$ \\
\midrule
1&Piedmont&470.69&59.25&35.5&7.04&90.97&1563&32.9&65.5&1040\\
2&Aosta Valley&582.58&61.14&243.5&3.11&51.09&102&26.1&49.5&31\\
%3&Liguria&531.68&48.81&7.3&5.02&56.67&672&29.2&49.2&522\\
\textbf{3}&\textbf{Liguria}&\textbf{531.68}&\textbf{48.81}&\textbf{7.3}&\textbf{5.02}&\textbf{56.67}&\textbf{672}&\textbf{29.2}&\textbf{49.2}&\textbf{522}\\
4&Lombardy&467.25&69.61&21.7&6.57&92.35&3581&32.8&62.4&2637\\
5&Trentino-South Tyrol&487.24&71.58&116.7&6.35&88.17&524&30.7&55&248\\
%6&Veneto&475.88&73.65&21.3&6.62&100.96&2105&35.6&61&1065\\
\textbf{6}&\textbf{Veneto}&\textbf{475.88}&\textbf{73.65}&\textbf{21.3}&\textbf{6.62}&\textbf{100.96}&\textbf{2105}&\textbf{35.6}&\textbf{61}&\textbf{1065}\\
%7&Friuli-Venezia-Giulia&484.11&65.48&23.3&6.73&216.43&585&34.9&63.9&250\\
\textbf{7}&\textbf{Friuli-Venezia Giulia}&\textbf{484.11}&\textbf{65.48}&\textbf{23.3}&\textbf{6.73}&\textbf{216.43}&\textbf{585}&\textbf{34.9}&\textbf{63.9}&\textbf{250}\\
8&Emilia-Romagna&642.54&63.83&19.2&9.49&131.52&1898&33.5&64.7&999\\
9&Tuscany&600&53.88&39.2&6.06&63.97&1403&30.5&60.2&1138\\
10&Umbria&508.39&61.69&37.2&4.4&33.04&395&33.2&51.9&230\\
%11&Marche&532.27&63.25&27&5.29&58.18&586&32.6&58.2&331\\
\textbf{11}&\textbf{Marche}&\textbf{532.27}&\textbf{63.25}&\textbf{27}&\textbf{5.29}&\textbf{58.18}&\textbf{586}&\textbf{32.6}&\textbf{58.2}&\textbf{331}\\
12&Lazio&503.97&45.52&13.2&6.51&23&1430&30.7&67&1046\\
13&Abruzzo&452.52&55.99&44.6&3.46&19.2&488&32.9&53.1&258\\
%14&Molise&376.96&30.72&84.4&3.37&2.93&101&34.5&49.6&71\\
\textbf{14}&\textbf{Molise}&\textbf{376.96}&\textbf{30.72}&\textbf{84.4}&\textbf{3.37}&\textbf{2.93}&\textbf{101}&\textbf{34.5}&\textbf{49.6}&\textbf{71}\\
15&Campania&439.06&52.76&26.4&3.28&9.64&1158&32.7&54.4&1244\\
16&Apulia&462.6&40.44&52.5&2.36&9.43&815&35.5&54.7&742\\
17&Basilicata&345.17&45.29&90.1&2.35&10.29&238&36.1&61.9&152\\
18&Calabria&394.61&39.67&72.6&1.8&9.13&325&38.2&50.4&487\\
19&Sicily&456.01&21.69&25.1&1.84&4.33&780&33.5&50.2&792\\
20&Sardinia&438.29&63.05&36&2.36&5.64&300&31.1&45.5&309\\
\bottomrule
\end{tabular}
\end{center}
\end{sidewaystable}
\begin{sidewaystable} %\setlength{\tabcolsep}{0.9pt} 
\footnotesize
\caption{Normalized performances. $(+)$ and $(-)$ denote elementary criteria having a positive or negative direction of preference, respectively.}
\label{data_normalized}
\begin{center}
\begin{tabular}{m{0.01\textwidth}m{0.06\textwidth}|m{0.06\textwidth}m{0.08\textwidth}m{0.06\textwidth}|m{0.09\textwidth}m{0.09\textwidth}m{0.06\textwidth}|m{0.06\textwidth}m{0.09\textwidth}m{0.09\textwidth}}
\toprule
&&\multicolumn{3}{c}{\textbf{Circular Economy - $g_1$}} & \multicolumn{3}{c}{\textbf{Innovation-Driven Development - $g_2$}}&\multicolumn{3}{c}{\textbf{Smart Specialization - $g_3$}}\\
\cmidrule(lr){3-11}
\textbf{ID}&\textbf{Region}&Urban Waste Generation (kg/capita) (2017) (-) $(g_{(\mathbf{1,1})})$ &Differentiated waste collection (\%) (2017) (+) $(g_{(\mathbf{1,2})})$ &Renewable electricity (percentage) (2017) (+) $(g_{(\mathbf{1,3})})$ &Research and Development Personnel (full-time equivalents per thousand inhabitants) (2017) (+) $(g_{(\mathbf{2,1})})$ &Patents registered at the European Patent Office (EPO) (number per million inhabitants) (2012) (+) $(g_{(\mathbf{2,2})})$ &Number of companies with ISO 14001 certification (2017) (+) $(g_{(\mathbf{2,3})})$ &Percentage of Smart specialized companies (2018) (+) $(g_{(\mathbf{3,1})})$ &Percentage of value added by Smart specialised companies to the total value added in the region (2018) (+) $(g_{(\mathbf{3,2})})$ &Number of Smart specialized companies with high intensity of investments in social and environmental responsibility (2018)(+) $(g_{(\mathbf{3,3})})$ \\
\midrule1&Piedmont&0.53&0.56&0.45&0.68&0.62&0.62&0.5&0.74&0.6\\
2&Aosta Valley&0.26&0.58&1&0.37&0.49&0.33&0.07&0.32&0.32\\
\textbf{3}&\textbf{Liguria}&\textbf{0.38}&\textbf{0.43}&\textbf{0.36}&\textbf{0.52}&\textbf{0.51}&\textbf{0.44}&\textbf{0.27}&\textbf{0.31}&\textbf{0.46}\\
4&Lombardy&0.54&0.69&0.4&0.65&0.62&1&0.5&0.66&1\\
5&Trentino-South Tyrol&0.49&0.71&0.71&0.63&0.61&0.41&0.36&0.46&0.38\\
\textbf{6}&\textbf{Veneto}&\textbf{0.52}&\textbf{0.74}&\textbf{0.4}&\textbf{0.65}&\textbf{0.65}&\textbf{0.73}&\textbf{0.67}&\textbf{0.62}&\textbf{0.61}\\
%7&Friuli-Venezia-Giulia&0.5&0.64&0.41&0.66&1&0.43&0.63&0.7&0.38\\
\textbf{7}&\textbf{Friuli-Venezia Giulia}&\textbf{0.5}&\textbf{0.64}&\textbf{0.41}&\textbf{0.66}&\textbf{1}&\textbf{0.43}&\textbf{0.63}&\textbf{0.7}&\textbf{0.38}\\
8&Emilia-Romagna&0.12&0.62&0.4&0.88&0.74&0.69&0.54&0.72&0.59\\
9&Tuscany&0.22&0.49&0.46&0.61&0.53&0.59&0.35&0.6&0.63\\
10&Umbria&0.44&0.59&0.45&0.48&0.43&0.39&0.52&0.38&0.37\\
%11&Marche&0.38&0.61&0.42&0.55&0.51&0.43&0.48&0.55&0.4\\
\textbf{11}&\textbf{Marche}&\textbf{0.38}&\textbf{0.61}&\textbf{0.42}&\textbf{0.55}&\textbf{0.51}&\textbf{0.43}&\textbf{0.48}&\textbf{0.55}&\textbf{0.4}\\
12&Lazio&0.45&0.39&0.38&0.64&0.4&0.6&0.36&0.78&0.6\\
13&Abruzzo&0.57&0.52&0.48&0.4&0.39&0.41&0.5&0.41&0.38\\
%14&Molise&0.75&0.21&0.6&0.4&0.34&0.33&0.6&0.32&0.33\\
\textbf{14}&\textbf{Molise}&\textbf{0.75}&\textbf{0.21}&\textbf{0.6}&\textbf{0.4}&\textbf{0.34}&\textbf{0.33}&\textbf{0.6}&\textbf{0.32}&\textbf{0.33}\\
15&Campania&0.6&0.48&0.42&0.39&0.36&0.54&0.49&0.45&0.66\\
16&Apulia&0.55&0.33&0.5&0.31&0.36&0.47&0.67&0.45&0.52\\
17&Basilicata&0.82&0.39&0.62&0.31&0.36&0.36&0.71&0.64&0.35\\
18&Calabria&0.71&0.32&0.57&0.27&0.36&0.37&0.84&0.34&0.45\\
19&Sicily&0.56&0.1&0.41&0.27&0.34&0.47&0.54&0.34&0.53\\
20&Sardinia&0.6&0.61&0.45&0.32&0.35&0.37&0.39&0.21&0.4\\
\bottomrule
\end{tabular}
\end{center}
\end{sidewaystable}
In Section \ref{ExamplePreciseInformation}, we shall describe the application of the new framework in case the DM, for example a policy maker, wants to evaluate the regions on Smart Specialization only and they are able to give precise information. In Section \ref{ExampleImpreciseInformation}, instead, we shall show the results of the application of the method in case the DM wants to rank the regions considering the whole set of criteria and they are not able to provide precise information both at global and partial levels. 

%%%%%%%%%%%%%%%%%%%%%%%%%%%%%%
\subsection{Basic Framework}\label{ExamplePreciseInformation}
%%%%%%%%%%%%%%%%%%%%%%%%%%%%%%
At first, let us assume the DM wants to evaluate the Regions only according to the Smart Specialization macro-criterion and that they can express their preferences about 5 regions only, namely $A^R = \{$Veneto, Friuli-Venezia Giulia, Marche, Liguria, Molise$\}=\{a_5^{R},a_4^{R},a_3^{R},a_2^{R},a_1^{R}\}$ articulated as shown in Table \ref{DM_preference_simple}
\begin{table}[h!]
    \centering
    \caption{DM's preferences}
    \label{DM_preference_simple}
    \begin{tabular}{ccc}
    \toprule
         ID&Region&Level\\
         \midrule
         $6$&Veneto&$\nu(a_5^R) = 65$\\[1mm]
         &&$e_{(\mathbf{3},4)}=8$\\[1mm]
         $7$&Friuli-Venezia Giulia&$\nu(a^R_4) = 56$\\[1mm]
         &&$e_{(\mathbf{3},3)} = 6$\\[1mm]
         $11$&Marche&$\nu(a^R_3) = 49$\\[1mm]
         &&$e_{(\mathbf{3},2)} = 5$\\[1mm]
         $3$&Liguria&$\nu(a^R_2) = 43$\\[1mm]
         &&$e_{(\mathbf{3},1)}=5$\\[1mm]
         $14$&Molise&$\nu(a^R_1) = 37$\\[1mm]
	&&$e_{(\mathbf{3},0)} = 36$\\
        \bottomrule
    \end{tabular}
\end{table}	.
To avoid any confusion, we denoted by $a_1^{R},\ldots,a_5^{R}$ the reference regions and we added in the first column the ID of the same alternative in Table \ref{alldata}. Moreover, in all tables regarding the whole set of regions, we put in bold the data corresponding to the reference regions.\\
Following the description in Section \ref{originalmodel}, according to the provided information, \textit{Molise} is the worst among the five reference regions, while \textit{Liguria} is the best among them. Moreover, 36 $(e_{(\mathbf{3},0))})$ blank cards are included between the ``fictitious zero alternative" and \textit{Molise}, 5 $(e_{(\mathbf{3},1)})$ blank cards between \textit{Molise} and \textit{Liguria} and so on.
Solving the LP problem (\ref{regressionmodel}) and considering the Weighted Sum as preference model we get the parameters shown in Table \ref{w_sum_par}. 
\begin{table}
\footnotesize
\caption{Basic Framework - Parameters' values obtained solving the LP problem (\ref{regressionmodel}) assuming the weighted sum as preference model}
\label{w_sum_par}
\setlength{\tabcolsep}{3pt}
\begin{center}
\begin{tabular}{ccccccccccccccc}
\toprule $w_{(\mathbf{3,1})}$&$w_{(\mathbf{3,2})}$&$w_{(\mathbf{3,3})}$&$k_{\mathbf{3}}$&$\sigma^{-}_{(\mathbf{3},5)}$&$\sigma^{+}_{(\mathbf{3},5)}$&$\sigma^{-}_{(\mathbf{3},4)}$&$\sigma^{+}_{(\mathbf{3},4)}$&$\sigma^{-}_{(\mathbf{3},3)}$&$\sigma^{+}_{(\mathbf{3},3)}$&$\sigma^{-}_{(\mathbf{3},2)}$&$\sigma^{+}_{(\mathbf{3},2)}$&$\sigma^{-}_{(\mathbf{3},1)}$&$\sigma^{+}_{(\mathbf{3},1)}$&$\overline{\sigma}$\\ \midrule %0.111&0.459&0.551&0.011&0&0&0&0&0.002&0&0&0.036&0&0&0.039\\ 
0.099&0.41&0.492&0.01&0&0&0&0&0.002&0&0&0.032&0&0&0.035\\
\bottomrule \end{tabular}
 \end{center}
 \end{table}

As one can see, the LP problem is not able to capture the DM's preferences without any estimation errors since $\overline{\sigma}=0.035$. For such a reason, we try to represent the same preferences using the 2-additive Choquet integral. Solving the same LP problem, we get that the Choquet integral is able to represent the DM's preferences without any errors since $\overline{\sigma} = 0$. At this point, we solve the LP problem (\ref{CompatibleMaximizingk}) to maximize the value of the blank card obtaining the M\"{o}bius parameters shown in Table \ref{ch_par}. 
\begin{table} 
\caption{Basic Framework - Parameters’ values obtained solving the LP problem (\ref{CompatibleMaximizingk}) assuming the 2-additive Choquet integral as preference model}
\label{ch_par}
\begin{center}
\begin{tabular}{cccccccc}
\toprule $m(\{g_1\})$&$m(\{g_2\})$&$m(\{g_3\})$&$m(\{g_1,g_2\})$&$m(\{g_1,g_3\})$&$m(\{g_2,g_3\})$&$k$&$\overline{\sigma}$\\ \midrule 0.0953&0.3561&0.7061&0.0712&-0.0953&-0.1334&0.0095&0\\ \bottomrule 
\end{tabular}
\end{center}
\end{table}

    As described in Section \ref{ROR}, to get more robust recommendations on the considered problem, we can compute the possible and necessary preference relations which are shown in Tables \ref{pos_preference} and \ref{nec_preference}, respectively. Referring to Table \ref{pos_preference} entry 1 means that the region in the row is possibly preferred to the region in the column, while entry 0 means that the region in the row is not possibly preferred to the region in the column (analogous interpretations can be given to the entries in Table \ref{nec_preference}). Looking at the two preference relations, it is possible to state that:
    \begin{itemize}        
        \item Lombardy ($ID = 4$) is necessarily preferred to all regions and, therefore, it can be considered the best among them;
        \item Piedmont ($ID = 1$), can be considered the third region since it is necessarily preferred to all regions apart from Lombardy and Lazio ($ID = 12$);   Emilia-Romagna ($ID = 8$) is necessarily preferred to all regions apart from Piedmont, Lombardy and Lazio, while Lazio is necessarily preferred to all regions apart from Piedmont, Lombardy, Emilia-Romagna, Veneto ($ID = 8$) and Campania ($ID = 15$);
        \item Aosta Valley ($ID = 2$), is the worst region since it is not possibly preferred to any region (apart from itself);
        \item Molise ($ID = 14$) is possibly preferred to Aosta Valley and Sardinia ($ID = 20$), only, and, analogously, Sardinia is possibly preferred to Aosta Valley and Molise only.
    \end{itemize}
    \begin{table}
\footnotesize
\caption{Basic Framework - Possible preference relation between regions on Smart Specialization}
\label{pos_preference}
\begin{center}
\begin{tabular}{ccccccccccccccccccccc}
\toprule ID&1&2&3&4&5&6&7&8&9&10&11&12&13&14&15&16&17&18&19&20\\
\midrule
1&1&1&1&0&1&1&1&1&1&1&1&1&1&1&1&1&1&1&1&1\\
2&0&1&0&0&0&0&0&0&0&0&0&0&0&0&0&0&0&0&0&0\\
\textbf{3}&\textbf{0}&\textbf{1}&\textbf{1}&\textbf{0}&\textbf{1}&\textbf{0}&\textbf{0}&\textbf{0}&\textbf{0}&\textbf{1}&\textbf{0}&\textbf{0}&\textbf{1}&\textbf{1}&\textbf{0}&\textbf{0}&\textbf{0}&\textbf{0}&\textbf{0}&\textbf{1}\\
4&1&1&1&1&1&1&1&1&1&1&1&1&1&1&1&1&1&1&1&1\\
5&0&1&0&0&1&0&0&0&0&1&0&0&1&1&0&0&0&0&0&1\\
\textbf{6}&\textbf{0}&\textbf{1}&\textbf{1}&\textbf{0}&\textbf{1}&\textbf{1}&\textbf{1}&\textbf{0}&\textbf{1}&\textbf{1}&\textbf{1}&\textbf{1}&\textbf{1}&\textbf{1}&\textbf{1}&\textbf{1}&\textbf{1}&\textbf{1}&\textbf{1}&\textbf{1}\\
\textbf{7}&\textbf{0}&\textbf{1}&\textbf{1}&\textbf{0}&\textbf{1}&\textbf{0}&\textbf{1}&\textbf{0}&\textbf{0}&\textbf{1}&\textbf{1}&\textbf{0}&\textbf{1}&\textbf{1}&\textbf{0}&\textbf{1}&\textbf{1}&\textbf{1}&\textbf{1}&\textbf{1}\\
8&0&1&1&0&1&1&1&1&1&1&1&1&1&1&1&1&1&1&1&1\\
9&0&1&1&0&1&0&1&0&1&1&1&0&1&1&1&1&1&1&1&1\\
10&0&1&0&0&0&0&0&0&0&1&0&0&0&1&0&0&0&0&0&1\\
\textbf{11}&\textbf{0}&\textbf{1}&\textbf{1}&\textbf{0}&\textbf{1}&\textbf{0}&\textbf{0}&\textbf{0}&\textbf{0}&\textbf{1}&\textbf{1}&\textbf{0}&\textbf{1}&\textbf{1}&\textbf{0}&\textbf{0}&\textbf{0}&\textbf{1}&\textbf{0}&\textbf{1}\\
12&1&1&1&0&1&1&1&1&1&1&1&1&1&1&1&1&1&1&1&1\\
13&0&1&0&0&1&0&0&0&0&1&0&0&1&1&0&0&0&0&0&1\\
\textbf{14}&\textbf{0}&\textbf{1}&\textbf{0}&\textbf{0}&\textbf{0}&\textbf{0}&\textbf{0}&\textbf{0}&\textbf{0}&\textbf{0}&\textbf{0}&\textbf{0}&\textbf{0}&\textbf{1}&\textbf{0}&\textbf{0}&\textbf{0}&\textbf{0}&\textbf{0}&\textbf{1}\\
15&0&1&1&0&1&0&1&0&1&1&1&1&1&1&1&1&1&1&1&1\\
16&0&1&1&0&1&0&0&0&0&1&1&0&1&1&0&1&1&1&1&1\\
17&0&1&1&0&1&0&0&0&0&1&1&0&1&1&0&0&1&1&1&1\\
18&0&1&1&0&1&0&0&0&0&1&1&0&1&1&0&0&0&1&0&1\\
19&0&1&1&0&1&0&0&0&0&1&1&0&1&1&0&0&1&1&1&1\\
20&0&1&0&0&0&0&0&0&0&0&0&0&0&1&0&0&0&0&0&1\\
\bottomrule
\end{tabular}
\end{center}
\end{table}
    \begin{table} \setlength{\tabcolsep}{6pt} \renewcommand\arraystretch{1}
\footnotesize
\caption{Basic Framework - Necessary preference relation between regions on Smart Specialization}
\label{nec_preference}
\begin{center} \begin{tabular}{ccccccccccccccccccccc}
\toprule ID&1&2&3&4&5&6&7&8&9&10&11&12&13&14&15&16&17&18&19&20\\
\midrule
1&1&1&1&0&1&1&1&1&1&1&1&0&1&1&1&1&1&1&1&1\\
2&0&1&0&0&0&0&0&0&0&0&0&0&0&0&0&0&0&0&0&0\\
\textbf{3}&\textbf{0}&\textbf{1}&\textbf{1}&\textbf{0}&\textbf{1}&\textbf{0}&\textbf{0}&\textbf{0}&\textbf{0}&\textbf{1}&\textbf{0}&\textbf{0}&\textbf{1}&\textbf{1}&\textbf{0}&\textbf{0}&\textbf{0}&\textbf{0}&\textbf{0}&\textbf{1}\\
4&1&1&1&1&1&1&1&1&1&1&1&1&1&1&1&1&1&1&1&1\\
5&0&1&0&0&1&0&0&0&0&1&0&0&0&1&0&0&0&0&0&1\\
\textbf{6}&\textbf{0}&\textbf{1}&\textbf{1}&\textbf{0}&\textbf{1}&\textbf{1}&\textbf{1}&\textbf{0}&\textbf{1}&\textbf{1}&\textbf{1}&\textbf{0}&\textbf{1}&\textbf{1}&\textbf{1}&\textbf{1}&\textbf{1}&\textbf{1}&\textbf{1}&\textbf{1}\\
\textbf{7}&\textbf{0}&\textbf{1}&\textbf{1}&\textbf{0}&\textbf{1}&\textbf{0}&\textbf{1}&\textbf{0}&\textbf{0}&\textbf{1}&\textbf{1}&\textbf{0}&\textbf{1}&\textbf{1}&\textbf{0}&\textbf{1}&\textbf{1}&\textbf{1}&\textbf{1}&\textbf{1}\\
8&0&1&1&0&1&1&1&1&1&1&1&0&1&1&1&1&1&1&1&1\\
9&0&1&1&0&1&0&1&0&1&1&1&0&1&1&0&1&1&1&1&1\\
10&0&1&0&0&0&0&0&0&0&1&0&0&0&1&0&0&0&0&0&1\\
\textbf{11}&\textbf{0}&\textbf{1}&\textbf{1}&\textbf{0}&\textbf{1}&\textbf{0}&\textbf{0}&\textbf{0}&\textbf{0}&\textbf{1}&\textbf{1}&\textbf{0}&\textbf{1}&\textbf{1}&\textbf{0}&\textbf{0}&\textbf{0}&\textbf{0}&\textbf{0}&\textbf{1}\\
12&0&1&1&0&1&0&1&0&1&1&1&1&1&1&0&1&1&1&1&1\\
13&0&1&0&0&0&0&0&0&0&1&0&0&1&1&0&0&0&0&0&1\\
\textbf{14}&\textbf{0}&\textbf{1}&\textbf{0}&\textbf{0}&\textbf{0}&\textbf{0}&\textbf{0}&\textbf{0}&\textbf{0}&\textbf{0}&\textbf{0}&\textbf{0}&\textbf{0}&\textbf{1}&\textbf{0}&\textbf{0}&\textbf{0}&\textbf{0}&\textbf{0}&\textbf{0}\\
15&0&1&1&0&1&0&1&0&0&1&1&0&1&1&1&1&1&1&1&1\\
16&0&1&1&0&1&0&0&0&0&1&1&0&1&1&0&1&1&1&1&1\\
17&0&1&1&0&1&0&0&0&0&1&1&0&1&1&0&0&1&1&0&1\\
18&0&1&1&0&1&0&0&0&0&1&0&0&1&1&0&0&0&1&0&1\\
19&0&1&1&0&1&0&0&0&0&1&1&0&1&1&0&0&0&1&1&1\\
20&0&1&0&0&0&0&0&0&0&0&0&0&0&0&0&0&0&0&0&1\\
\bottomrule
\end{tabular}
\end{center}
\end{table}

    Since there is often no clear preference between many pairs of alternatives, we chose to use the SMAA as described in Section \ref{SMAASection}. The SMAA application is based on a sampling of $1$ million compatible measures computing, therefore, the RAIs and the PWIs shown in Tables \ref{RAI_simple} and \ref{PWI_simple}, respectively. 
    \begin{table} 
\footnotesize
\caption{Basic Framework - Rank Acceptability Indices (RAIs) on Smart Specialization}
\label{RAI_simple}
\setlength{\tabcolsep}{1pt} \renewcommand\arraystretch{1} 
\begin{center} \begin{tabular}{>{\centering}m{0.06\textwidth}|cccccccccccccccccccc}
\toprule Rank ID&1&2&3&4&5&6&7&8&9&10&11&12&13&14&15&16&17&18&19&20\\
\midrule
1&0&79.75&20.25&0&0&0&0&0&0&0&0&0&0&0&0&0&0&0&0&0\\
2&0&0&0&0&0&0&0&0&0&0&0&0&0&0&0&0&0&0&0&100\\
%3&0&0&0&0&0&0&0&0&0&0&0&0&0&100&0&0&0&0&0&0\\
\textbf{3}&\textbf{0}&\textbf{0}&\textbf{0}&\textbf{0}&\textbf{0}&\textbf{0}&\textbf{0}&\textbf{0}&\textbf{0}&\textbf{0}&\textbf{0}&\textbf{0}&\textbf{0}&\textbf{100}&\textbf{0}&\textbf{0}&\textbf{0}&\textbf{0}&\textbf{0}&\textbf{0}\\
4&100&0&0&0&0&0&0&0&0&0&0&0&0&0&0&0&0&0&0&0\\
5&0&0&0&0&0&0&0&0&0&0&0&0&0&0&44.94&55.06&0&0&0&0\\
%6&0&0&0&32.44&67.56&0&0&0&0&0&0&0&0&0&0&0&0&0&0&0\\
\textbf{6}&\textbf{0}&\textbf{0}&\textbf{0}&\textbf{32.44}&\textbf{67.56}&\textbf{0}&\textbf{0}&\textbf{0}&\textbf{0}&\textbf{0}&\textbf{0}&\textbf{0}&\textbf{0}&\textbf{0}&\textbf{0}&\textbf{0}&\textbf{0}&\textbf{0}&\textbf{0}&\textbf{0}\\
%7&0&0&0&0&0&0&0&100&0&0&0&0&0&0&0&0&0&0&0&0\\
\textbf{7}&\textbf{0}&\textbf{0}&\textbf{0}&\textbf{0}&\textbf{0}&\textbf{0}&\textbf{0}&\textbf{100}&\textbf{0}&\textbf{0}&\textbf{0}&\textbf{0}&\textbf{0}&\textbf{0}&\textbf{0}&\textbf{0}&\textbf{0}&\textbf{0}&\textbf{0}&\textbf{0}\\
8&0&0&50.96&49.04&0&0&0&0&0&0&0&0&0&0&0&0&0&0&0&0\\
9&0&0&0&0&0&10.05&89.95&0&0&0&0&0&0&0&0&0&0&0&0&0\\
10&0&0&0&0&0&0&0&0&0&0&0&0&0&0&0&0&100&0&0&0\\
%11&0&0&0&0&0&0&0&0&0&0&0&73.73&26.27&0&0&0&0&0&0&0\\
\textbf{11}&\textbf{0}&\textbf{0}&\textbf{0}&\textbf{0}&\textbf{0}&\textbf{0}&\textbf{0}&\textbf{0}&\textbf{0}&\textbf{0}&\textbf{0}&\textbf{73.73}&\textbf{26.27}&\textbf{0}&\textbf{0}&\textbf{0}&\textbf{0}&\textbf{0}&\textbf{0}&\textbf{0}\\
12&0&20.25&28.79&18.51&19.04&13.4&0&0&0&0&0&0&0&0&0&0&0&0&0&0\\
13&0&0&0&0&0&0&0&0&0&0&0&0&0&0&55.06&44.94&0&0&0&0\\
%14&0&0&0&0&0&0&0&0&0&0&0&0&0&0&0&0&0&27.87&72.13&0\\
\textbf{14}&\textbf{0}&\textbf{0}&\textbf{0}&\textbf{0}&\textbf{0}&\textbf{0}&\textbf{0}&\textbf{0}&\textbf{0}&\textbf{0}&\textbf{0}&\textbf{0}&\textbf{0}&\textbf{0}&\textbf{0}&\textbf{0}&\textbf{0}&\textbf{27.87}&\textbf{72.13}&\textbf{0}\\
15&0&0&0&0&13.4&76.55&10.05&0&0&0&0&0&0&0&0&0&0&0&0&0\\
16&0&0&0&0&0&0&0&0&100&0&0&0&0&0&0&0&0&0&0&0\\
17&0&0&0&0&0&0&0&0&0&84.3&15.7&0&0&0&0&0&0&0&0&0\\
18&0&0&0&0&0&0&0&0&0&0&0&26.27&73.73&0&0&0&0&0&0&0\\
19&0&0&0&0&0&0&0&0&0&15.7&84.3&0&0&0&0&0&0&0&0&0\\
20&0&0&0&0&0&0&0&0&0&0&0&0&0&0&0&0&0&72.13&27.87&0\\
\bottomrule
\end{tabular}
\end{center}
\end{table}
    \begin{table}
\caption{Basic Framework - Pairwise Winning Indices (PWIs) on Smart Specialization}
\label{PWI_simple}
\footnotesize
\setlength{\tabcolsep}{1pt} \renewcommand\arraystretch{1} 
\begin{center} 
\begin{tabular}{c|cccccccccccccccccccc}
\toprule ID&1&2&3&4&5&6&7&8&9&10&11&12&13&14&15&16&17&18&19&20\\
\midrule
1&0&100&100&0&100&100&100&100&100&100&100&79.75&100&100&100&100&100&100&100&100\\
2&0&0&0&0&0&0&0&0&0&0&0&0&0&0&0&0&0&0&0&0\\
%3&0&100&0&0&100&0&0&0&0&100&0&0&100&100&0&0&0&0&0&100\\
\textbf{3}&\textbf{0}&\textbf{100}&\textbf{0}&\textbf{0}&\textbf{100}&\textbf{0}&\textbf{0}&\textbf{0}&\textbf{0}&\textbf{100}&\textbf{0}&\textbf{0}&\textbf{100}&\textbf{100}&\textbf{0}&\textbf{0}&\textbf{0}&\textbf{0}&\textbf{0}&\textbf{100}\\
4&100&100&100&0&100&100&100&100&100&100&100&100&100&100&100&100&100&100&100&100\\
5&0&100&0&0&0&0&0&0&0&100&0&0&44.94&100&0&0&0&0&0&100\\
%6&0&100&100&0&100&0&100&0&100&100&100&32.44&100&100&100&100&100&100&100&100\\
\textbf{6}&\textbf{0}&\textbf{100}&\textbf{100}&\textbf{0}&\textbf{100}&\textbf{0}&\textbf{100}&\textbf{0}&\textbf{100}&\textbf{100}&\textbf{100}&\textbf{32.44}&\textbf{100}&\textbf{100}&\textbf{100}&\textbf{100}&\textbf{100}&\textbf{100}&\textbf{100}&\textbf{100}\\
%7&0&100&100&0&100&0&0&0&0&100&100&0&100&100&0&100&100&100&100&100\\
\textbf{7}&\textbf{0}&\textbf{100}&\textbf{100}&\textbf{0}&\textbf{100}&\textbf{0}&\textbf{0}&\textbf{0}&\textbf{0}&\textbf{100}&\textbf{100}&\textbf{0}&\textbf{100}&\textbf{100}&\textbf{0}&\textbf{100}&\textbf{100}&\textbf{100}&\textbf{100}&\textbf{100}\\
8&0&100&100&0&100&100&100&0&100&100&100&50.96&100&100&100&100&100&100&100&100\\
9&0&100&100&0&100&0&100&0&0&100&100&0&100&100&10.05&100&100&100&100&100\\
10&0&100&0&0&0&0&0&0&0&0&0&0&0&100&0&0&0&0&0&100\\
%11&0&100&100&0&100&0&0&0&0&100&0&0&100&100&0&0&0&73.73&0&100\\
\textbf{11}&\textbf{0}&\textbf{100}&\textbf{100}&\textbf{0}&\textbf{100}&\textbf{0}&\textbf{0}&\textbf{0}&\textbf{0}&\textbf{100}&\textbf{0}&\textbf{0}&\textbf{100}&\textbf{100}&\textbf{0}&\textbf{0}&\textbf{0}&\textbf{73.73}&\textbf{0}&\textbf{100}\\
12&20.25&100&100&0&100&67.56&100&49.04&100&100&100&0&100&100&86.6&100&100&100&100&100\\
13&0&100&0&0&55.06&0&0&0&0&100&0&0&0&100&0&0&0&0&0&100\\
%14&0&100&0&0&0&0&0&0&0&0&0&0&0&0&0&0&0&0&0&27.87\\
\textbf{14}&\textbf{0}&\textbf{100}&\textbf{0}&\textbf{0}&\textbf{0}&\textbf{0}&\textbf{0}&\textbf{0}&\textbf{0}&\textbf{0}&\textbf{0}&\textbf{0}&\textbf{0}&\textbf{0}&\textbf{0}&\textbf{0}&\textbf{0}&\textbf{0}&\textbf{0}&\textbf{27.87}\\
15&0&100&100&0&100&0&100&0&89.95&100&100&13.4&100&100&0&100&100&100&100&100\\
16&0&100&100&0&100&0&0&0&0&100&100&0&100&100&0&0&100&100&100&100\\
17&0&100&100&0&100&0&0&0&0&100&100&0&100&100&0&0&0&100&84.3&100\\
18&0&100&100&0&100&0&0&0&0&100&26.27&0&100&100&0&0&0&0&0&100\\
19&0&100&100&0&100&0&0&0&0&100&100&0&100&100&0&0&15.7&100&0&100\\
20&0&100&0&0&0&0&0&0&0&0&0&0&0&72.13&0&0&0&0&0&0\\
\bottomrule
\end{tabular}
\end{center}
\end{table}

\noindent Analyzing these indices, it is possible to notice that:
    \begin{itemize}
        \item Aosta Valley ($ID = 2$) is always ranked last ($b_{\mathbf{3}}^{20}(a_{2})=100\%$) while Lombardy ($ID = 4$) is always ranked first ($b_{\mathbf{3}}^{1}(a_{4})=100\%$);
        \item Second to last and third to last positions are taken by Molise ($ID = 14$) and  Sardinia ($ID = 20$). In particular, Molise is most frequently in the second to last position ($b_{\mathbf{3}}^{19}(a_{14})=72.13\%$), while,  Sardinia is in the third to last position with the same frequency ($b_{\mathbf{3}}^{18}(a_{20})=72.13\%$);
        \item Piedmont ($ID = 1$) is always in the second and third positions ($b_{\mathbf{3}}^{2}(a_{1})=79.75\%$ and $b_{\mathbf{3}}^{3}(a_{1})=20.25\%$), Emilia-Romagna ($ID = 8$) fills the third or the fourth position in the whole ranking ($b_{\mathbf{3}}^{3}(a_{8})=50.96\%$ and $b_{\mathbf{3}}^{4}(a_{8})=49.04\%$), while Veneto ($ID=6$) is always in fourth or fifth position ($b_{\mathbf{3}}^{4}(a_{6})=32.44\%$ and $b_{\mathbf{3}}^{5}(a_{6})=67.56\%$); 
        %Lazio ($ID = 12$) occupies a position between the second and the sixth with similar frequencies; 
        \item Tuscany ($ID = 7$), Trentino-South Tyrol ($ID = 16$), Emilia Romagna ($ID = 3$) and Basilicata ($ID = 10$) are always set on ranks $8$, $9$, $14$ and $17$, respectively. 
    \end{itemize}
    Considering the PWIs, the data in Table \ref{PWI_simple} show that in most of the cases, there is an evident preference for one alternative over the other since $p_{\mathbf{r}}(a,b)=100\%$ and, therefore, $p_{\mathbf{r}}(b,a)=0\%$. The only doubtful cases are the following: 
    \begin{itemize}
    \item Trentino-South Tyrol ($ID=5$) is preferred to Abruzzo ($ID=13$) with a frequency of $44.94\%$, while the vice versa is true for $55.06\%$ of the cases;
    \item Emilia-Romagna ($ID=8$) and Lazio ($ID=12$) are preferred one to the other with similar frequencies since $p_{\mathbf{3}}(8,12)=50.96\%$, while $p_{\mathbf{3}}(12,8)=49.04\%$. 
    \end{itemize}
    In all other cases, one region is preferred to the other with a frequency not lower than 67.56\%. \\
    
    \begin{table}[!h]
\centering
\footnotesize
\caption{Basic Framework - Expected Ranking and ranking obtained using the approximation of the barycenter shown in Table \ref{bar_par}}
\label{ER_BAR_simple}
\begin{subtable}{0.5\textwidth}
\centering
\caption{Expected Ranking}
\label{ER_simple}
\begin{tabular}{cccc}
\toprule
Rank&Region&ID&$ER(\cdot)$\\
\midrule
1&Lombardy&4&100\\
2&Piedmont&1&220.25\\
3&Emilia-Romagna&8&349.04\\
4&Lazio&12&376.55\\
\textbf{5}&\textbf{Veneto}&\textbf{6}&\textbf{467.56}\\
6&Campania&15&596.64\\
7&Tuscany&9&689.95\\
\textbf{8}&\textbf{Friuli-Venezia Giulia}&\textbf{7}&\textbf{800}\\
9&Apulia&16&900\\
10&Basilicata&17&1015.7\\
11&Sicily&19&1084.3\\
\textbf{12}&\textbf{Marche}&\textbf{11}&\textbf{1226.27}\\
13&Calabria&18&1273.73\\
\textbf{14}&\textbf{Liguria}&\textbf{3}&\textbf{1400}\\
15&Abruzzo&13&1544.94\\
16&Trentino-South Tyrol&5&1555.06\\
17&Umbria&10&1700\\
18&Sardinia&20&1827.87\\
\textbf{19}&\textbf{Molise}&\textbf{14}&\textbf{1872.13}\\
20&Aosta Valley&2&2000\\
\bottomrule
\end{tabular}
\end{subtable}% <---- don't forget this %
\begin{subtable}{0.5\textwidth}
\centering
\caption{Ranking based on barycenter parameters}
\label{BAR_simple_rank}
\begin{tabular}{cccc}
\toprule 
Rank&Region&ID&$U(a_j)$\\
\midrule
1&Lombardy&4&0.88213\\
2&Piedmont&1&0.63429\\
3&Emilia-Romagna&8&0.62726\\
4&Lazio&12&0.62696\\
\textbf{5}&\textbf{Veneto}&\textbf{6}&\textbf{0.61956}\\
6&Campania&15&0.60514\\
7&Tuscany&9&0.58144\\
\textbf{8}&\textbf{Friuli-Venezia Giulia}&\textbf{7}&\textbf{0.53377}\\
9&Apulia&16&0.51954\\
10&Basilicata&17&0.51014\\
11&Sicily&19&0.49586\\
\textbf{12}&\textbf{Marche}&\textbf{11}&\textbf{0.46705}\\
13&Calabria&18&0.46239\\
\textbf{14}&\textbf{Liguria}&\textbf{3}&\textbf{0.40986}\\
15&Abruzzo&13&0.40613\\
16&Trentino-South Tyrol&5&0.40601\\
17&Umbria&10&0.39053\\
18&Sardinia&20&0.36098\\
\textbf{19}&\textbf{Molise}&\textbf{14}&\textbf{0.35267}\\
20&Aosta Valley&2&0.27889\\
\bottomrule
\end{tabular}		
\end{subtable}
\end{table}
    Let us observe that the barycenter shown in Table \ref{bar_par} is computed as the average, component by component, of the sampled measures. 
    \begin{table}
	\footnotesize
	\caption{Basic Framework - Parameters of the approximated barycenter}
	\label{bar_par}
	\setlength{\tabcolsep}{6pt}
	\renewcommand\arraystretch{1}
	\begin{center}
		\begin{tabular}{cccccc}
			\toprule 
			$m(\{g_{(\mathbf{3,1})}\})$&$m(\{g_{(\mathbf{3,2})}\})$&$m(\{g_{(\mathbf{3,3})}\})$&$m(\{g_{(\mathbf{3,1})},g_{(\mathbf{3,2})}\})$&$m(\{g_{(\mathbf{3,1})},g_{(\mathbf{3,3})}\})$&$m(\{g_{(\mathbf{3,2})},g_{(\mathbf{3,3})}\})$\\
			%\midrule
			0.0918&0.3540&0.7306&0.0768&-0.0100&-0.2432\\
			\bottomrule
		\end{tabular}
	\end{center}
\end{table}

    The first observation is that the two rankings are the same. On the one hand, considering the top of the ranking, Lombardy is in the first position, followed by Piedmont, while Emilia-Romagna and Lazio, which had similar frequencies to be preferred one to the other, are third and fourth, respectively. On the other hand, looking at the bottom of the same ranking, the last three places are taken by Sardinia, Molise, and Aosta Valley.\\
    Commenting on the barycenter shown in Table \ref{bar_par}, taking into account the weights assigned to single criteria without considering interactions, one can observe that Number of Smart specialized companies with high intensity of investments in social and environmental responsibility (2018) ($g_{\mathbf{(3,3)}}$) has a greater weight than Percentage of value added by Smart specialized companies to the total value added in the region (2018) $(g_{(\mathbf{3,2})})$ that, in turn, has a greater weight  than Percentage of Smart specialized companies (2018) $(g_{(\mathbf{3,1})})$. Moreover, there is a positive interaction between $g_{(\mathbf{3,1})}$ and $g_{(\mathbf{3,2})}$ $\left(m\left(\{g_{(\mathbf{3,1})},g_{(\mathbf{3,2})}\}\right)=0.0768\right)$, while there is a negative interaction between $g_{(\mathbf{3,1})}$ and $g_{(\mathbf{3,3})}$ $\left(m\left(\{g_{(\mathbf{3,1})},g_{(\mathbf{3,3})}\}\right)=-0.0100\right)$ as well as between $g_{(\mathbf{3,2})}$ and $g_{(\mathbf{3,3})}$ $\left(m\left(\{g_{(\mathbf{3,2})},g_{(\mathbf{3,3})}\}\right)=-0.2432\right)$.
    %\input{ER_simple}
    %\input{BAR_simple}

    %%%%%%%%%%%%%%%%%%%%%%%%%%%%%%%%%%%%%%%%%%%%%%%%%%%%%%%%%%%%%%%%%%%%%%%%%
    \subsection{Imprecise Information}\label{ExampleImpreciseInformation}%%%%
    %%%%%%%%%%%%%%%%%%%%%%%%%%%%%%%%%%%%%%%%%%%%%%%%%%%%%%%%%%%%%%%%%%%%%%%%%
    Let's suppose in this case that the DM provides their preferences on the same 5 regions above considered both globally ($g_0$) and on the single macro-criteria (Circular Economy - $g_1$, Innovation-Driven Development - $g_2$, Smart Specialization - $g_3$). Moreover, for each macro-criterion, the DM provides interval or imprecise information (see Section \ref{ImprecisionSection}). Let us assume that the preference is articulated as shown in Table \ref{DM_preference_general} and let us comment on it.
    \begin{table}[!h]
\scriptsize
    \centering
    \caption{Imprecise Information - DM's preference at global and macro-criteria level}
    \label{DM_preference_general}
    \begin{subtable}{0.9\textwidth}
    \caption{Global level - $g_{(\mathbf{0})}$}
    \label{MissingImpreciseInformationGlobal}
    \centering
    \begin{tabular}{>{\centering}p{0.1\columnwidth}>{\centering}p{0.3\columnwidth}>{\centering\arraybackslash}m{0.5\columnwidth}}
    \toprule
    ID&Region&Level\\
    \midrule
    6 &Veneto&\multirow{2}{*}{$e_{\mathbf{(0},4)}\in\left[e_{\mathbf{(0},4)}^L,e_{\mathbf{(0},4)}^U\right] = \left[1,7\right]$}\\[1.5mm]
    &&\\
    7 &Friuli-Venezia Giulia &  \multirow{2}{*}{$e_{\mathbf{(0},3)}\in\left[e_{\mathbf{(0},3)}^L,\;e_{\mathbf{(0},3)}^U\right] = \left[7,?\right]$}\\[1.5mm]
    & & \\
    11 & Marche & \multirow{2}{*}{$e_{\mathbf{(0},2)}\in\left[e_{\mathbf{(0},2)}^L,\;e_{\mathbf{(0},2)}^U\right] = \left[1,6\right]$}\\[1.5mm]
    & &  \\
    3 & Liguria & \multirow{2}{*}{$e_{\mathbf{(0},1)}\in\left[e_{\mathbf{(0},1)}^L,\;e_{\mathbf{(0},1)}^U\right] = \left[?,5\right]$}\\[1.5mm]
    &&  \\
    14 & Molise & \multirow{2}{*}{$e_{\mathbf{(0},0)}\in\left[e_{\mathbf{(0},0)}^L,\;e_{\mathbf{(0},0)}^U\right] = \left[40,50\right]$}\\
    & & \\
    %\bottomrule
    \end{tabular}
    \end{subtable}
    \begin{subtable}{0.9\textwidth}
    \caption{Circular Economy - $g_{(\mathbf{1})}$}
    \label{MissingImpreciseInformationL1}
    \centering
    \begin{tabular}{>{\centering}p{0.1\columnwidth}>{\centering}p{0.3\columnwidth}>{\centering\arraybackslash}p{0.5\columnwidth}}
    \toprule
    ID&Region&Level\\
    \midrule
    6 & Veneto & \multirow{2}{*}{$e_{\mathbf{(1},4)}\in\left[e_{\mathbf{(1},4)}^{L},e_{\mathbf{(1},4)}^{U}\right]=\left[1,6\right]$}\\[1.5mm]
    && \\
    7 & Friuli-Venezia Giulia & \multirow{2}{*}{$e_{\mathbf{(1},3)}\in\left[e_{\mathbf{(1},3)}^{L},e_{\mathbf{(1},3)}^{U}\right]=\left[?,4\right]$}\\[1.5mm]
    & & \\
    11 & Marche & \multirow{2}{*}{$e_{\mathbf{(1},2)}\in\left[e_{\mathbf{(1},2)}^{L},e_{\mathbf{(1},2)}^{U}\right]=\left[1,3\right]$}\\[1.5mm]
    && \\
    14 & Molise & \multirow{2}{*}{$e_{\mathbf{(1},1)}\in\left[e_{\mathbf{(1},1)}^{L},e_{\mathbf{(1},1)}^{U}\right]=\left[?,5\right]$}\\[1.5mm]
    && \\
    3 &Liguria& \multirow{2}{*}{$e_{\mathbf{(1},0)}\in\left[_{\mathbf{(1},0)}^{L},e_{\mathbf{(1},0)}^{U}\right]=\left[9,16\right]$}\\[1.5mm]
    && \\
    %\bottomrule
    \end{tabular}
    \end{subtable}
    \begin{subtable}{0.9\textwidth}
    \caption{Innovation-Driven - $g_{(\mathbf{2})}$}
    \label{MissingImpreciseInformationL2}
    \centering
    \begin{tabular}{>{\centering}p{0.1\columnwidth}>{\centering}p{0.3\columnwidth}>{\centering\arraybackslash}p{0.5\columnwidth}}
    \toprule
    ID&Region&Level\\
    \midrule
    7 & Friuli-Venezia Giulia & \multirow{2}{*}{$e_{\mathbf{(2},4)}\in\left[e_{\mathbf{(2},4)}^{L},e_{\mathbf{(2},4)}^{U}\right]=\left[1,5\right]$}\\[1.5mm]
    && \\
    6 & Veneto &  \multirow{2}{*}{$e_{\mathbf{(2},3)}\in\left[e_{\mathbf{(2},3)}^{L},e_{\mathbf{(2},3)}^{U}\right]=\left[2,7\right]$}\\[1.5mm]
    && \\
    11 & Marche & \multirow{2}{*}{$e_{\mathbf{(2},2)}\in\left[e_{\mathbf{(2},2)}^{L},e_{\mathbf{(2},2)}^{U}\right]=\left[?,5\right]$}\\[1.5mm]
    && \\
    3 & Liguria & \multirow{2}{*}{$e_{\mathbf{(2},1)}\in\left[e_{\mathbf{(2},1)}^{L},e_{\mathbf{(2},1)}^{U}\right]=\left[1,6\right]$}\\[1.5mm]
    && \\
    14 & Molise & \multirow{2}{*}{$e_{\mathbf{(2},0)}\in\left[e_{\mathbf{(2},0)}^{L},e_{\mathbf{(2},0)}^{U}\right]=\left[10,20\right]$}\\[1.5mm]
    && \\
    %\bottomrule
    \end{tabular}
    \end{subtable}
    \begin{subtable}{0.9\textwidth}
    \caption{Smart Specialization - $g_{(\mathbf{3})}$}
    \label{MissingImpreciseInformationL3}
    \centering
    \begin{tabular}{>{\centering}p{0.1\columnwidth}>{\centering}p{0.3\columnwidth}>{\centering\arraybackslash}p{0.5\columnwidth}}
    \toprule
    ID&Region&Level\\
    \midrule
    6 & Veneto & \multirow{2}{*}{$e_{\mathbf{(3},4)}\in\left[e_{\mathbf{(3},4)}^{L},e_{\mathbf{(3},4)}^{U}\right]=\left[1,6\right]$}\\[1.5mm]
    && \\
    7 & Friuli-Venezia Giulia & \multirow{2}{*}{$e_{\mathbf{(3},3)}\in\left[e_{\mathbf{(3},3)}^{L},e_{\mathbf{(3},3)}^{U}\right]=\left[3,?\right]$}\\[1.5mm]
    && \\
    11 & Marche & \multirow{2}{*}{$e_{\mathbf{(3},2)}\in\left[e_{\mathbf{(3},2)}^{L},e_{\mathbf{(3},2)}^{U}\right]=\left[2,5\right]$}\\[1.5mm]
    && \\
    3 & Liguria & \multirow{2}{*}{$e_{\mathbf{(3},1)}\in\left[e_{\mathbf{(3},1)}^{L},e_{\mathbf{(3},1)}^{U}\right]=\left[?,7\right]$}\\[1.5mm]
    && \\
    14 & Molise & \multirow{2}{*}{$e_{\mathbf{(3},0)}\in\left[e_{\mathbf{(3},0)}^{L},e_{\mathbf{(3},0)}^{U}\right]=\left[10,18\right]$}\\[1.5mm]
    && \\
   \bottomrule
    \end{tabular}
    \end{subtable}
\end{table}

    Considering the information in Table \ref{MissingImpreciseInformationGlobal}, the DM retains that, at the global level, Molise is the least preferred region, followed by Liguria, Marche, Friuli-Venezia Giulia and, finally, Veneto, being the most preferred. On the one hand, the DM provides imprecise information regarding the number of blank cards between the zero fictitious level and Molise $\left(\left[e_{\mathbf{(0,0)}}^L,\;e_{\mathbf{(0,0)}}^R\right] = \left[40,50\right]\right)$, the number of blank cards between Liguria and Marche $\left(\left[e_{\mathbf{(0,2)}}^L,\;e_{\mathbf{(0,2)}}^R\right] = \left[1,6\right]\right)$ and the number of blank cards between Friuli-Venezia Giulia and Veneto $\left(\left[e_{\mathbf{(0,4)}}^L,e_{\mathbf{(0,4)}}^U\right] = \left[1,7\right]\right)$. On the other hand, the DM gives missing information on the number of blank cards between Molise and Liguria since $e_{\mathbf{(0,1)}}^L=?$ and for the number of blank cards between Marche and Friuli-Venezia Giulia since $e_{\mathbf{(0,3)}}^U=?.$ The information contained in Tables \ref{MissingImpreciseInformationL1}-\ref{MissingImpreciseInformationL3} can be interpreted analogously. One thing that should be underlined is that using the MCHP within the proposed framework (see Section \ref{RORSMAAMCHPMissingInformationSection}) the DM can provide different information on the considered macro-criteria. For example, while the five reference regions are ordered in the same way at the global level (see Table \ref{MissingImpreciseInformationGlobal}) and on Smart Specialization (see Table \ref{MissingImpreciseInformationL3}), the same does not happen for the other two macro-criteria. Indeed, while Liguria is preferred to Molise at the global level, the vice versa is true on Circular Economy (see Table \ref{MissingImpreciseInformationL1}). Analogously, while Veneto is preferred to Friuli-Venezia Giulia at the global level, the vice versa is true on Innovation-Driven Development (see Table \ref{MissingImpreciseInformationL2}). This permits the DM to give more detailed information not only at the global but also at the partial one.

    Solving the LP problem (\ref{ImpreciseMCHP}) and assuming the weighted sum as preference model, we find $\overline{\sigma}=0$. Therefore, this time, the weighted sum can represent the preferences given by the DM. Solving the LP problem (\ref{EpsilonMCHPImprecise}) maximizing the blank cards' values, we get the weights of the elementary criteria shown in Table \ref{ElementaryCriteriaWeightsMissingInformation} as well as the number of blank cards and the value of one blank card (globally and on all macro-criteria) shown in Table \ref{BlankCardsMissingInformation}. 

    \begin{table}[h!]
\caption{Imprecise Information - Elementary criteria weights obtained solving LP problem (\ref{EpsilonMCHPImprecise}) and assuming the weighted sum as preference model}
\label{ElementaryCriteriaWeightsMissingInformation}
\begin{center}
\begin{tabular}{ccccccccc}
\toprule
$w_{(\mathbf{1},1)}$&$w_{(\mathbf{1},2)}$&$w_{(\mathbf{1},3)}$&$w_{(\mathbf{2},1)}$&$w_{(\mathbf{2},2)}$&$w_{(\mathbf{2},3)}$&$w_{(\mathbf{3},1)}$&$w_{(\mathbf{3},2)}$&$w_{(\mathbf{3},3)}$\\
\midrule
%0.1118&0.1382&0&0.4459&0.0547&0.0119&0.0376&0.1845&0.1725\\
0.097&0.119&0&0.385&0.047&0.01&0.033&0.159&0.149\\
\bottomrule
\end{tabular}
\end{center}
\end{table}

    \begin{table}
\caption{Imprecise Information - Number of blank cards and value of a single blank card both at the global and partial levels}
\label{BlankCardsMissingInformation}
\begin{center}
\begin{tabular}{cc|cc|cc|cc}
\toprule
\multicolumn{2}{c}{$g_{(\mathbf{0})}$}&\multicolumn{2}{c}{$g_{(\mathbf{1})}$}&\multicolumn{2}{c}{$g_{(\mathbf{2})}$}&\multicolumn{2}{c}{$g_{(\mathbf{3})}$}\\
\midrule
$e_{(\mathbf{0},4)}$&1.1604&$e_{(\mathbf{1},4)}$&1.2139&$e_{(\mathbf{2},4)}$&1&$e_{(\mathbf{3},4)}$&2.7161\\
$e_{(\mathbf{0},3)}$&10.2139&$e_{(\mathbf{1},3)}$&1.2674&$e_{(\mathbf{2},3)}$&4.9015&$e_{(\mathbf{3},3)}$&3\\
$e_{(\mathbf{0},2)}$&6&$e_{(\mathbf{1},2)}$&1&$e_{(\mathbf{2},2)}$&0&$e_{(\mathbf{3},2)}$&4.8209\\
$e_{(\mathbf{0},1)}$&5&$e_{(\mathbf{1},1)}$&0.3637&$e_{(\mathbf{2},1)}$&6&$e_{(\mathbf{3},1)}$&0\\
$e_{(\mathbf{0},0)}$&40&$e_{(\mathbf{1},0)}$&13.0545&$e_{(\mathbf{2},0)}$&19.3095&$e_{(\mathbf{3},0)}$&18\\
%$k_{\mathbf{0}}$&0.011&$k_{\mathbf{1}}$&0.0073&$k_{\mathbf{2}}$&0.0098&$k_{\mathbf{3}}$&0.0073\\
$k_{\mathbf{0}}$&0.0095&$k_{\mathbf{1}}$&0.0063&$k_{\mathbf{2}}$&0.0085&$k_{\mathbf{3}}$&0.0063\\
\bottomrule
\end{tabular}
\end{center}
\end{table}
    
    As one can see, the number of blank cards to be included between two successive levels is perfectly compatible with the DM's preferences. For example, while the DM stated that $e_{(\mathbf{0},4)}\in\left[1,7\right]$ (see Table \ref{MissingImpreciseInformationGlobal}), the number of blank cards obtained solving the LP problem is $e_{(\mathbf{0},4)}=1.1604$ (see Table \ref{BlankCardsMissingInformation}). Analogously, while the DM stated that $e_{(\mathbf{3},3)}\in\left[3,?\right]$ (see Table \ref{MissingImpreciseInformationL3}), the value obtained solving the previously mentioned LP problem is $e_{(\mathbf{3},3)}=3$ (see Table \ref{BlankCardsMissingInformation}). Considering, instead, the weights obtained for the considered elementary criteria, it seems that Renewable Electricity (percentage) (2017) has not any importance since $w_{(\mathbf{1},3)}=0$, while, the most important elementary criterion is Research and Development Personnel (full-time equivalents per thousand inhabitants) (2017) ($w_{(\mathbf{2},1)}=0.385$), followed by Percentage of value added by Smart Specialized companies to the total value added in the region (2018) $(w_{(\mathbf{3},2)}=0.159)$ and Number of Smart Specialized companies with high intensity of investments in social and environmental responsibility (2018) $(w_{(\mathbf{3},3)}=0.149)$.

To get more robust recommendations on the problem at hand, we computed the ROR and SMAA methodologies. To save space, we do not report here all the corresponding tables that are provided, instead, as supplementary material. Let us comment on the results obtained by SMAA. 

%Applying the ROR, one can observe that the necessary and possible preference relations shown in the previous section are different from those obtained now. Indeed, while those preference relations were based on SMART Specialization information only, here, the two preference relations are obtained considering information both at the global and partial levels. For example, while Lombardia ($ID=4$) was necessarily preferred to all other regions on SMART specialization, now, Lombardia is not necessarily preferred to three regions being Veneto ($ID=5$), Campania ($ID=6$) and Friuli-Venezia-Giulia ($ID=8$). Considering instead, the possible preference relation in most of the cases one region is globally possibly preferred to the other and vice-versa. 

Considering the RAIs at the global and partial levels, the following can be observed: 
\begin{itemize}
\item At the global level, Lombardy ($ID=4$) is always at the first position $(b_{\mathbf{0}}^{1}(a_{4})=100\%)$, while, Sicily ($ID=19$) is always in the last rank position $(b_{\mathbf{0}}^{20}(a_{19})=100\%)$; 
\item On Circular Economy, Veneto ($ID=6$) is always first $(b_{\mathbf{1}}^{1}(a_{6})=100\%)$, while Sicily is always at the bottom of the ranking $(b_{\mathbf{1}}^{20}(a_{19})=100\%)$;
\item On Innovation-Driven Development, Emilia-Romagna ($ID=8$) is robustly in the first position $(b_{\mathbf{2}}^{1}(a_{8})=100\%)$, while Calabria ($ID=18$) is always in the last rank position $(b_{\mathbf{2}}^{20}(a_{18})=100\%)$; 
\item On Smart Specialization, Lombardy is in the first rank position in all considered cases $(b_{\mathbf{3}}^{1}(a_{4})=100\%)$, while the last position is taken by Aosta Valley ($ID=2$) and Sardinia ($ID=20$) with frequencies $90.19\%$ and $9.81\%$, respectively.
\end{itemize}
%To better clarify the importance of taking into account preference information not only at the global level but also at the partial one and, analogously, to have more detailed recommendations at the partial level, let us comment on the results obtained by SMAA. \\

Even if the previous analysis permits to appreciate the recommendations on the considered regions at the partial level, let us underline even more this aspect considering the data in Table \ref{RAI_MCHP_comparison}. 
\begin{table}
\footnotesize
\caption{Imprecise Information - Comparison between RAIs computed at global and partial levels by the MCHP: Focus on Lombardy ($ID=4$) and Trentino-South Tyrol ($ID=5$)}
\label{RAI_MCHP_comparison}
\setlength{\tabcolsep}{2pt}
\begin{center}
\begin{tabular}{>{\centering}m{0.06\textwidth}|ccccccccccccccccccccc}
\toprule
Rank ID&&1&2&3&4&5&6&7&8&9&10&11&12&13&14&15&16&17&18&19&20\\
\midrule
\multirow{4}{*}{4}&$g_{(\mathbf{0})}$&100&0&0&0&0&0&0&0&0&0&0&0&0&0&0&0&0&0&0&0\\
&$g_{(\mathbf{1})}$&0&61.07&38.93&0&0&0&0&0&0&0&0&0&0&0&0&0&0&0&0&0\\
&$g_{(\mathbf{2})}$&0&0&0&0&100&0&0&0&0&0&0&0&0&0&0&0&0&0&0&0\\
&$g_{(\mathbf{3})}$&100&0&0&0&0&0&0&0&0&0&0&0&0&0&0&0&0&0&0&0\\
\midrule
\multirow{4}{*}{5}&$g_{(\mathbf{0})}$&0&0&0&0&0&0&100&0&0&0&0&0&0&0&0&0&0&0&0&0\\
&$g_{(\mathbf{1})}$&0&38.93&61.07&0&0&0&0&0&0&0&0&0&0&0&0&0&0&0&0&0\\
&$g_{(\mathbf{2})}$&0&0&0&0&0&100&0&0&0&0&0&0&0&0&0&0&0&0&0&0\\
&$g_{(\mathbf{3})}$&0&0&0&0&0&0&0&0&0&0&0&0&22.39&77.61&0&0&0&0&0&0\\
\bottomrule
\end{tabular}
\end{center}
\end{table} 
\noindent In the table, we report the RAIs at the global and partial levels for two regions, that is, Lombardy ($ID=4$) and Trentino-South Tyrol ($ID=5$). As one can see, while at the global level and on Smart Specialization, Lombardy is always in the first rank position, the same cannot be said for Circular Economy and Innovation-Driven Development. On the one hand, considering Circular Economy, Lombardy is in the second or third positions with frequencies of 61.07\% and 38.93\%, respectively, while, on the other hand, it is always fifth on Innovation-Driven Development.\\
The different recommendations provided by the introduced framework are even more evident considering Trentino-South Tyrol. Indeed, while at the global level, it is always in the seventh position, on Circular Economy it is between the second and the third positions ($b_{\mathbf{1}}^{2}(a_{5})=38.93\%)$ and $(b_{\mathbf{1}}^{3}(a_{5})=61.07\%$), it is always in the sixth position on Innovation-Driven Development, and, finally, it is between positions 13th and 14th on Smart Specialization with frequencies of 22.39\% and 77.61\%, respectively.   

To summarize the results of the RAIs, let us show in Table \ref{ExpectedRankingTableMissingInformation} the expected rankings both at the global and partial levels.

\begin{table}
\centering
\footnotesize
\caption{Imprecise Information - Expected rankings of the considered regions both at the global and partial levels}
\label{ExpectedRankingTableMissingInformation}
\begin{subtable}{0.5\textwidth}
\centering
\caption{Global level - $g_{(\mathbf{0})}$}
\label{}
\begin{tabular}{cccc}
\toprule
Rank&Region&ID&$ER(\cdot)$\\
\midrule
%			1 & Lombardia & 4 & -100 \\ 
%			2 & Emilia-Romagna & 8 & -200 \\ 
%			3 & Piemonte & 1 & -312.38 \\ 
%			\textbf{4} & \textbf{Veneto} & \textbf{6} & \textbf{-387.62} \\ 
%			\textbf{5} & \textbf{Friuli-Venezia-Giulia} & \textbf{7} & \textbf{-500} \\ 
%			6 & Lazio & 12 & -600 \\ 
%			7 & Trentino-Alto Adige & 5 & -700 \\ 
%			8 & Tuscany & 9 & -800 \\ 
%			\textbf{9} & \textbf{Marche} & \textbf{11} & \textbf{-900} \\ 
%			10 & Campania & 15 & -1000 \\ 
%			11 & Umbria & 10 & -1100 \\ 
%			12 & Basilicata & 17 & -1217.35 \\ 
%			\textbf{13} & \textbf{Liguria} & \textbf{3} & \textbf{-1282.65} \\ 
%			14 & Abruzzo & 13 & -1400 \\ 
%			15 & Puglia & 16 & -1500 \\ 
%			\textbf{16} & \textbf{Molise} & \textbf{14} & \textbf{-1600} \\ 
%			17 & Sardegna & 20 & -1700 \\ 
%			18 & Calabria & 18 & -1808.29 \\ 
%			19 & Valle d'Aosta & 2 & -1891.71 \\ 
%			20 & Sicilia & 19 & -2000 \\ 
1 & Lombardy & 4 & 100 \\ 
2 & Emilia-Romagna & 8 & 200 \\ 
3 & Piedmont & 1 & 345.76 \\ 
\textbf{4} & \textbf{Veneto} & \textbf{6} & \textbf{354.24} \\ 
\textbf{5} & \textbf{Friuli-Venezia Giulia} & \textbf{7} & \textbf{500} \\ 
6 & Lazio & 12 & 600 \\ 
7 & Trentino-South Tyrol & 5 & 700 \\ 
8 & Tuscany & 9 & 800 \\ 
\textbf{9} & \textbf{Marche} & \textbf{11} & \textbf{900} \\ 
10 & Campania & 15 & 1000 \\ 
11 & Umbria & 10 & 1100 \\ 
12 & Basilicata & 17 & 1226.97 \\ 
\textbf{13} & \textbf{Liguria} & \textbf{3} & \textbf{1273.03} \\ 
14 & Abruzzo & 13 & 1400 \\ 
15 & Apulia & 16 & 1500 \\ 
\textbf{16} & \textbf{Molise} & \textbf{14} & \textbf{1608.87} \\ 
17 & Sardinia & 20 & 1691.13 \\ 
18 & Calabria & 18 & 1834.82 \\ 
19 & Aosta Valley & 2 & 1865.18 \\ 
20 & Sicily & 19 & 2000 \\
\bottomrule
\end{tabular}
\end{subtable}% <---- don't forget this %
\begin{subtable}{0.5\textwidth}
\centering
\caption{Circular economy - $g_{(\mathbf{1})}$}
\label{}
\begin{tabular}{cccc}
\toprule 
Rank&Region&ID&$ER(\cdot)$\\
\midrule
%			\textbf{1} & \textbf{Veneto} & \textbf{6} & \textbf{-100} \\ 
%			2 & Lombardia & 4 & -217.8 \\ 
%			3 & Trentino-Alto Adige & 5 & -282.2 \\ 
%			4 & Sardegna & 20 & -400 \\ 
%			5 & Basilicata & 17 & -500 \\ 
%			\textbf{6} & \textbf{Friuli-Venezia-Giulia} & \textbf{7} & \textbf{-600} \\ 
%			7 & Piemonte & 1 & -700 \\ 
%			8 & Abruzzo & 13 & -800 \\ 
%			9 & Campania & 15 & -900 \\ 
%			10 & Umbria & 10 & -1000 \\ 
%			\textbf{11} & \textbf{Marche} & \textbf{11} & \textbf{-1100} \\ 
%			12 & Calabria & 18 & -1200 \\ 
%			13 & Molise & 14 & -1322.68 \\ 
%			%14 & Valle d'Aosta & 2 & -1377.32 \\ 
%			%15 & Puglia & 16 & -1500 \\ 
%			14 & Puglia & 16 & -1500 \\ 
%			15 & Valle d'Aosta & 2 & -1377.32 \\ 
%			16 & Lazio & 12 & -1600 \\ 
%			\textbf{17} & \textbf{Liguria} & \textbf{3} & \textbf{-1700} \\ 
%			18 & Emilia-Romagna & 8 & -1800 \\ 
%			19 & Tuscany & 9 & -1900 \\ 
%			20 & Sicilia & 19 & -2000 \\ 
\textbf{1} & \textbf{Veneto} & \textbf{6} & \textbf{100} \\ 
2 & Lombardy & 4 & 238.93 \\ 
3 & Trentino-South Tyrol & 5 & 261.07 \\ 
4 & Sardinia & 20 & 400 \\ 
5 & Basilicata & 17 & 500 \\ 
\textbf{6} & \textbf{Friuli-Venezia Giulia} & \textbf{7} & \textbf{600} \\ 
7 & Piedmont & 1 & 700 \\ 
8 & Abruzzo & 13 & 800 \\ 
9 & Campania & 15 & 900 \\ 
10 & Umbria & 10 & 1000 \\ 
\textbf{11} & \textbf{Marche} & \textbf{11} & \textbf{1100} \\ 
12 & Calabria & 18 & 1200 \\ 
13 & Aosta Valley & 2 & 1348.36 \\ 
\textbf{14} & \textbf{Molise} & \textbf{14} & \textbf{1351.64} \\ 
15 & Apulia & 16 & 1500 \\ 
16 & Lazio & 12 & 1600 \\ 
\textbf{17} & \textbf{Liguria} & \textbf{3} & \textbf{1700} \\ 
18 & Emilia-Romagna & 8 & 1800 \\ 
19 & Tuscany & 9 & 1900 \\ 
20 & Sicily & 19 & 2000 \\ 
\bottomrule
\end{tabular}		
\end{subtable}
\begin{subtable}{0.5\textwidth}
\centering
\caption{Innovation-Driven - $g_{(\mathbf{2})}$}
\label{}
\begin{tabular}{cccc}
\toprule 
Rank&Region&ID&$ER(\cdot)$\\
\midrule
%			1 & Emilia-Romagna & 8 & -100 \\ 
%			\textbf{2} & \textbf{Friuli-Venezia-Giulia} & \textbf{7} & \textbf{-200} \\ 
%			3 & Piemonte & 1 & -300 \\ 
%			\textbf{4} & \textbf{Veneto} & \textbf{6} & \textbf{-400} \\ 
%			5 & Lombardia & 4 & -500 \\ 
%			6 & Trentino-Alto Adige & 5 & -600 \\ 
%			7 & Lazio & 12 & -700 \\ 
%			8 & Tuscany & 9 & -800 \\ 
%			\textbf{9} & \textbf{Marche} & \textbf{11} & \textbf{-900} \\ 
%			\textbf{10} & \textbf{Liguria} & \textbf{3} & \textbf{-1000} \\ 
%			11 & Umbria & 10 & -1100 \\ 
%			12 & Abruzzo & 13 & -1200 \\ 
%			\textbf{13} & \textbf{Molise} & \textbf{14} & \textbf{-1307.94} \\ 
%			14 & Campania & 15 & -1392.59 \\ 
%			15 & Valle d'Aosta & 2 & -1499.47 \\ 
%			16 & Puglia & 16 & -1600 \\ 
%			17 & Basilicata & 17 & -1700 \\ 
%			18 & Sardegna & 20 & -1800 \\ 
%			19 & Sicilia & 19 & -1900 \\ 
%			20 & Calabria & 18 & -2000 \\ 
1 & Emilia-Romagna & 8 & 100 \\ 
\textbf{2} & \textbf{Friuli-Venezia Giulia} & \textbf{7} & \textbf{200} \\ 
3 & Piedmont & 1 & 300 \\ 
\textbf{4} & \textbf{Veneto} & \textbf{6} & \textbf{400} \\ 
5 & Lombardy & 4 & 500 \\ 
6 & Trentino-South Tyrol & 5 & 600 \\ 
7 & Lazio & 12 & 700 \\ 
8 & Tuscany & 9 & 800 \\ 
\textbf{9} & \textbf{Marche} & \textbf{11} & \textbf{900} \\ 
\textbf{10} & \textbf{Liguria} & \textbf{3} & \textbf{1000} \\ 
11 & Umbria & 10 & 1100 \\ 
12 & Abruzzo & 13 & 1200 \\ 
\textbf{13} & \textbf{Molise} & \textbf{14} & \textbf{1310.59} \\ 
14 & Campania & 15 & 1389.44 \\ 
15 & Aosta Valley & 2 & 1499.97 \\ 
16 & Apulia & 16 & 1600 \\ 
17 & Basilicata & 17 & 1700 \\ 
18 & Sardinia & 20 & 1800 \\ 
19 & Sicily & 19 & 1900 \\ 
20 & Calabria & 18 & 2000 \\
\bottomrule
\end{tabular}		
\end{subtable}%
\begin{subtable}{0.5\textwidth}
\centering
\caption{Smart Specialization - $g_{(\mathbf{3})}$}
\label{}
\begin{tabular}{cccc}
\toprule 
Rank&Region&ID&$ER(\cdot)$\\
\midrule
%			1 & Lombardia & 4 & -100 \\ 
%			2 & Lazio & 12 & -200 \\ 
%			3 & Piemonte & 1 & -300 \\ 
%			4 & Emilia-Romagna & 8 & -400 \\ 
%			\textbf{5} & \textbf{Veneto} & \textbf{6} & \textbf{-500} \\ 
%			6 & Tuscany & 9 & -600 \\ 
%			\textbf{7} & \textbf{Friuli-Venezia-Giulia} & \textbf{7} & \textbf{-700} \\ 
%			8 & Campania & 15 & -800 \\ 
%			9 & Basilicata & 17 & -900 \\ 
%			10 & Puglia & 16 & -1000 \\ 
%			\textbf{11} & \textbf{Marche} & \textbf{11} & \textbf{-1100} \\ 
%			12 & Sicilia & 19 & -1200 \\ 
%			13 & Calabria & 18 & -1337.46 \\ 
%			14 & Trentino-Alto Adige & 5 & -1362.54 \\ 
%			15 & Abruzzo & 13 & -1500 \\ 
%			16 & Umbria & 10 & -1600 \\ 
%			\textbf{17} & \textbf{Liguria} & \textbf{3} & \textbf{-1700} \\ 
%			\textbf{18} & \textbf{Molise} & \textbf{14} & \textbf{-1800} \\ 
%			19 & Sardegna & 20 & -1916.24 \\ 
%			20 & Valle d'Aosta & 2 & -1983.76 \\ 
1 & Lombardy & 4 & 100 \\ 
2 & Lazio & 12 & 200 \\ 
3 & Piedmont & 1 & 300 \\ 
4 & Emilia-Romagna & 8 & 400 \\ 
\textbf{5} & \textbf{Veneto} & \textbf{6} & \textbf{500} \\ 
6 & Tuscany & 9 & 600 \\ 
\textbf{7} & \textbf{Friuli-Venezia Giulia} & \textbf{7} & \textbf{700} \\ 
8 & Campania & 15 & 800 \\ 
9 & Basilicata & 17 & 900 \\ 
10 & Apulia & 16 & 1000 \\ 
\textbf{11} & \textbf{Marche} & \textbf{11} & \textbf{1100} \\ 
12 & Sicily & 19 & 1200 \\ 
13 & Calabria & 18 & 1322.39 \\ 
14 & Trentino-South Tyrol & 5 & 1377.61 \\ 
15 & Abruzzo & 13 & 1500 \\ 
16 & Umbria & 10 & 1600 \\ 
\textbf{17} & \textbf{Liguria} & \textbf{3} & \textbf{1700} \\ 
\textbf{18} & \textbf{Molise} & \textbf{14} & \textbf{1800} \\ 
19 & Sardinia & 20 & 1909.81 \\ 
20 & Aosta Valley & 2 & 1990.19 \\
\bottomrule
\end{tabular}		
\end{subtable}
\end{table}

Looking at the data in the tables, one has the confirmation of the goodness of Lombardy being in the first rank-position at the global level and on Smart Specialization, while, it is second-ranked on Circular Economy and fifth-ranked on Innovation-Driven Development. Analogously, looking at the bottom of the rankings, Sicily performs badly at the global level and on Circular Economy being placed in the last position but also on Innovation-Driven Development since it is in the second to last position. However, it is not so badly ranked on Smart Specialization since it is placed in the 12th rank position. 

Regarding some other regions, a floating behavior can be observed. For example, Trentino-South Tyrol is seventh at the global level, while, it is third on Circular Economy, sixth on Innovation-Driven Development and fourteenth on Smart Specialization. This means that, while Circular Economy can be considered a strong point of the region, Smart Specialization can be considered a weak point deserving improvement. Analogously, Sardinia, is seventeenth at the global level, while it is ranked fourth on Circular economy, eighteenth on Innovation-Driven Development and nineteenth on Smart Specialization.\\
This information can be used by policymakers and political governments to evaluate which are the weak and strong points of each region and, consequently, build politics aiming to improve the weak points by pushing on the strong points.

%%%%%%%%%%%%%%%%%%%%%%%%%%%%%%%%%%%%%%%%%%%%%%%%%%%     
\section{Conclusions}\label{ConclusionSection}%%%%%
%%%%%%%%%%%%%%%%%%%%%%%%%%%%%%%%%%%%%%%%%%%%%%%%%%%

Preference parameter elicitation procedures are fundamental to Multiple Criteria Decision Aiding (MCDA) and must take into account two main requirements: on the one hand, the preference information requested from the Decision Maker (DM) must be as simple and easy as possible, and on the other hand, it must be as rich and precise as possible. In this perspective, we considered the newly introduced Deck of cards based Ordinal Regression (DOR) \citep{barbati2023new} that permits to collect preference information not only in terms of ranking of reference alternatives but also of intensity of preferences. This can be done through the Deck of Cards Method (DCM) that, in general, is perceived as a relatively intuitive and straightforward method and, as discussed at length in this paper, also with other well-known MCDA procedures such as AHP, BWM and MACBETH. In this context, we took into consideration some issues that are particularly relevant in real-world applications: robustness concerns; imprecision, incompleteness and ill determination of the preference information; and hierarchy of considered criteria. With this aim, Robust Ordinal Regression, Stochastic Multicriteria Acceptability Analysis, and Multiple Criteria Hierarchy Process were extended to DOR. Observe that the proposed approach puts together two main research streams in MCDA, that is, on the one hand, the scaling procedures related to DCM, AHP and MACBETH, and, on the other hand, the ordinal regression with its extensions such as Robust Ordinal Regression. The resulting methodology is very flexible and versatile. Indeed, according to the requirements and the previous experiences of the DM,  one can select 1) the most appropriate scaling procedure such as DCM, AHP, BWM, MACBETH or other analogous methodologies acquiring and processing DM's pairwise preferences comparisons between reference alternatives;
2) the most adequate formal  model of a value function, e.g. weighted sum, piecewise additive value function, Choquet integral, assigning an overall evaluation to each alternative;
3) the most satisfactory form to represent preferences on the set of alternatives taking into account imprecision and indetermination in the preference parameters, e.g. possible and necessary preferences, pairwise preferences probability or probability of attaining a given ranking. This so rich and ductile methodology can be successfully applied to complex real-world decision problems characterized by (i) many heterogeneous criteria structured in a hierarchical way; (ii) vague and approximate preference information; (iii) plurality of experts, stakeholders, policymakers and, in general, diversified actors participating to the decision process. We have shown the potential of the proposed methodology in a didactic example related to the ranking of Italian Regions with respect to Circular Economy, Innovation-Driven Development and Smart Specialization Strategies. The example is based on a large ongoing national project financed by the European Union - NextGenerationEU where Italian Regions are evaluated on several elementary criteria so that those considered in our example are just a small subset. Let us underline that the methodology can be applied across various contexts such as environmental or human development and potentially in the construction of composite indicators in any domain (for an updated state-of-the-art survey on composite indicators see \citealt{GrecoEtAl2019}). 

For the future, we plan to apply the proposed methodology to some real-world decision problems aiming, on the one hand, to test its effectiveness, reliability and validity and, on the other hand, to develop some customized procedures and protocols for applications in specific domains. Related to this, we would like to test the advantages of the proposed methodology in some decision experiments that could unveil spaces for further improvements. A further domain for future research is the extension of DOR to sorting decision problems in which the alternatives, rather than ordered from best to worst, have to be assigned to predefined ordered classes. We would like also to explore a possible extension of DOR to outranking methods such as ELECTRE \citep{FigueiraEtAl2013,GovindanJepsen2016} and PROMETHEE \citep{BehzadianEtAl2010,BransVincke1985}.

%In this paper, we consider  the recently introduced application of the Deck of Cards Method (DCM) to ordinal regression proposing an extension to Robust Ordinal Regression and Stochastic Multiattribute Acceptability Analysis. In Multiple Criteria Decision Aiding context, the proposed methodology permits to assign a value to each alternative evaluated on a set of criteria hierarchically structured. The Decision Maker can provide precise or imprecise information at different levels of the hierarchy of criteria using the classical DCM framework. This information is therefore used to infer a value function compatible with it. The compatible value function can be a simple weighted sum, a piecewise linear value function, a general monotonic value function, or a Choquet integral. To provide robust recommendations to the Decision Maker, we consider the Robust Ordinal Regression and the Stochastic Multicriteria Acceptability Analysis because, even if in different ways, both of them take into account the whole set of models compatible with the preference information provided by the Decision Maker. The applicability of the proposed methodology is shown by a didactic example in which Italian regions are evaluated on criteria representing Circular Economy, Innovation Driven Development and Smart Specialization Strategies.
%
%
%
%
\section*{Acknowledgments}
\noindent This study was funded by the European Union -  NextGenerationEU, in the framework of the GRINS -Growing Resilient, INclusive and Sustainable project (GRINS PE00000018  CUP E63C22002120006). The views and opinions expressed are solely those of the authors and do not necessarily reflect those of the European Union, nor can the European Union be held responsible for them. \bibliographystyle{agsm}
\bibliography{Full_bibliography}
%\printbibliography
%%%%%%%%%%%%%%%%%%%%%%%%%%%%%%%%%%%%%%%%%%%%%%%%%%%%%%%%%%%%%%%%%
\appendix
%\section*{Appendix A}
\section{Different preference models}\label{PreferenceModels}%%%%
%%%%%%%%%%%%%%%%%%%%%%%%%%%%%%%%%%%%%%%%%%%%%%%%%%%%%%%%%%%%%%%%%
The procedure described in Sections \ref{originalmodel} - \ref{ImpreciseInformation}, proposing to apply the DCM to build a value function compatible with a few DM's preferences, is independent of the type of the assumed preference model $U$. For this reason, in the following, we shall recall four different types of value functions that can be used under the proposed framework to build a compatible value function, namely, a weighted sum, a piecewise linear value function \citep{JacquetLagrezeSiskos1982}, a general additive value function \citep{GrecoMousseauSlowinski2008} and the Choquet integral value function \citep{Choquet1953,Grabisch1996}. The choice of the considered value function will determine the constraints composing the $E^{Model}$ set mentioned at first in the LP problem (\ref{regressionmodel}) included in Section \ref{originalmodel} and, therefore, present in the subsequent mathematical problems. The four preference functions are described in the following.
\begin{description}  
\item[Weighted sum:] denoting by $w_1, \ldots, w_n$ the weights of criteria $g_1, \ldots, g_n$, $E^{Model}$ is the following set of constraints: 
$$
\left.
\begin{array}{l}
U(a)=\displaystyle\sum_{i=1}^{n}w_i\cdot g_i(a)\;\mbox{for all}\; a\in A,\\[1mm]
w_i\geqslant 0,\;\mbox{for all}\;i=1,\ldots,n, \\[1mm]
\displaystyle\sum_{i=1}^{n}w_i=1.
\end{array}
\right\}E^{Model}_{WS} 
$$
Under the MCHP framework, the value assigned by $U$ to alternative $a$ on macro-criterion $g_{\mathbf{r}}$ is $U_{\mathbf{r}}\left(a\right)=\displaystyle\sum_{\mathbf{t}\in E(g_{\mathbf{r}})}w_{\mathbf{t}}\cdot g_{\mathbf{t}}(a)$ where $w_{\mathbf{t}}$ is the weight attached to the elementary criterion $g_{\mathbf{t}}$ and they are such that $w_{\mathbf{t}}\geqslant 0$ for all $\mathbf{t}\in EL$. In this case, the parameters defining the value function $U$ are the weights $w_{\mathbf{t}}$, $\mathbf{t}\in EL$;
\item[Piecewise linear value function:] in this case, the value function $U$ is such that $U(a)=\displaystyle\sum_{i=1}^{n} u_i(g_i(a))$ for all $a\in A$ and each $u_i:A\rightarrow\rea$ is a piecewise linear value function for each $g_i\in G$. Assuming that all criteria are expressed on a quantitative scale and $X_i=\{x_i^{1},\ldots,x_i^{m_i}\}\subseteq \mathbb{R}$ is the set of possible values that can be obtained on criterion $g_i\in G$, each function $u_i$ is defined by $\gamma_i$ breakpoints $y_i^1, \ldots, y_i^{\gamma_i}\in\mathbb{R}$ such that $y_i^1<\cdots<y_i^{\gamma_i}$ and 
$$
\left[x_i^{1},x_{i}^{m_i}\right]=\left[y_{i}^{1},y_{i}^{2}\right]\cup\cdots\cup\left[y_{i}^{\gamma_i-1},y_{i}^{\gamma_i}\right].
$$
Therefore, the marginal value function $u_i$ will be defined by means of $u_i(y_i^1),\ldots,u_i(y_i^{\gamma_i})$ only so that, if $g_i(a)\in[y_i^{q-1},y^{q}_{i}]$, $q=2,\ldots,\gamma_i$, we have 
$$
u_i(g_i(a))=\frac{y_{i}^q-g_i(a)}{y_i^q-y_i^{q-1}}\cdot u_i\left(y_i^{q-1}\right)+\frac{g_i(a)-y_i^{q-1}}{y_i^{q}-y_i^{q-1}}\cdot u_i\left(y_i^{q}\right)
$$
or, equivalently, 
$$
u_i(g_i(a))=\frac{g_i(a)-y_i^{q-1}}{y_i^{q}-y_i^{q-1}}\cdot \left(u_i\left(y_i^{q}\right)-u_i\left(y_i^{q-1}\right)\right)+u_i(y_i^{q-1}).
$$
Consequently, $E^{Model}$ is replaced by the following set of constraints:
$$
\left.
\begin{array}{l}
U(a)=\displaystyle\sum_{i=1}^{n} u_i(g_i(a))\;\mbox{for all}\; a\in A,\\[1mm]
u_i(y_i^{q-1}) \leqslant u_i(y_i^{q}),\; \mbox{for all}\;i= 1, \ldots, n, \;\mbox{and for all}\;q=2,\ldots,\gamma_i,\\[2mm]
u_i(x_i^1) = 0, \; \mbox{for all}\; i = 1, \ldots, n,\\[1mm]
\displaystyle\sum_{i=1}^{n} u_i\left(x_i^{m_i}\right) = 1.
\end{array}
\right\}E^{Model}_{PS} 
$$
Under the MCHP framework, the value of $a$ w.r.t. macro-criterion $g_{\mathbf{r}}$ is equal to $U_{\mathbf{r}}\left(a\right)=\displaystyle\sum_{\mathbf{t}\in E\left(g_{\mathbf{r}}\right)}u_{\mathbf{t}}\left(g_{\mathbf{t}}(a)\right)$. In this case, the set of criteria $G$ has to be replaced by the set of elementary criteria $\{g_{\mathbf{t}},\; \mathbf{t}\in EL\}$;
\item[General additive value function:] Denoting by $X_i=\{x_i^{1},\ldots,x_i^{m_i}\}\subseteq \mathbb{R}$ the set of possible values that can be obtained on criterion $g_i\in G$, the marginal value function $u_i$ depends on $u_i(x_i^{1}),\ldots,u_i(x_i^{m_i})$ only. $E^{Model}$ is replaced by the following set of constraints:
$$
\left.
\begin{array}{l}
U(a)=\displaystyle\sum_{i=1}^{n} u_i(g_i(a))\;\mbox{for all}\; a\in A,\\[1mm]
u_i(x_i^{f-1}) \leqslant u_i(x_i^{f}),\; \mbox{for all}\;i= 1, \ldots, n\; \mbox{and for all}\; f = 2,\ldots,m_i,\\[2mm]
u_i(x_i^1) = 0, \; \mbox{for all}\; i = 1, \ldots, n,\\[2mm]
\displaystyle\sum_{i=1}^{n} u_i\left(x_i^{m_i}\right) = 1.
\end{array}
\right\}E^{Model}_{GA}
$$
Under the MCHP framework the definition of $U_{\mathbf{r}}$ is the same considered when the preference model is a piecewise linear value function. The only difference w.r.t. the previous case is that for each marginal value function $u_{\mathbf{t}}$, $\mathbf{t}\in EL$, the breakpoints coincide with the evaluations of the alternatives on that elementary criterion, that is, $y_{\mathbf{t}}^{1}=x_{\mathbf{t}}^1$, $y_{\mathbf{t}}^{2}=x_{\mathbf{t}}^2$, $\ldots$, $y_{\mathbf{t}}^{\gamma_{\mathbf{t}}}=x_{\mathbf{t}}^{m_{\mathbf{t}}}$;
\item[Choquet integral value function:] as known in literature, a value function $U$ can be written in an additive way ($U(a) = \displaystyle \sum_{g_i \in G} u_i(g_i(a))$) iff the criteria from $G$ are mutually preferentially independent \citep{Wakker1989}. However, in real-world applications, the criteria present generally a certain degree of interaction. In particular, given $g_i, g_{i'}\in G$, we say that $g_i$ and $g_{i'}$ are positively (negatively) interacting if the importance given to them together is greater (lower) than the sum of their importance when considered alone. To take into account these interactions, non additive integrals are used in literature \citep{GrabischLabreuche2016} and the most known is the Choquet integral. It is based on a capacity being a set function $\mu:2^G\rightarrow [0,1]$ such that the following constraints are satisfied:
\begin{description}  
    \item[1a)] $\mu(R)\leqslant \mu(S)$ for all $R\subseteq S \subseteq G$ (monotonicity),
    \item[2a)] $\mu(\emptyset) = 0$ and $\mu(G) = 1$ (normalization).
\end{description}
Known $\mu$, the Choquet integral of $(g_1(a), \ldots, g_n(a))$ is computed as follows:
\begin{equation*}
\label{Choquetcapacity}
Ch_{\mu}(a) = \displaystyle\sum_{i = 1}^{n} [g_{(i)}(a) - g_{(i-1)}(a)]\cdot \mu(A_{(i)})    
\end{equation*}
where $_{(\cdot)}$ is a permutation of the indices of criteria in $G$ such that $0 = g_{(0)}(a) \leqslant g_{(1)}(a)\leqslant \cdots \leqslant g_{(n)}(a)$ and $A_{(i)} = \{g_{i^{'}}\in G:\; g_{i^{'}}(a)\geqslant g_{(i)}(a)\}$.\\
Since the use of the Choquet integral implies the definition of $2^{|G|}-2$ values (one for each subset of $G$ different from $\emptyset$ and $G$), in real-world applications the M\"{o}bius transformation of $\mu$ \citep{Rota1964} and $k$-additive capacities \citep{Grabisch1997} are used:
\begin{itemize}
    \item the M\"{o}bius transformation of a capacity $\mu$ is a set function $m:2^G\rightarrow \mathbb{R}$ such that $\mu(T) = \displaystyle \sum_{S\subseteq T}m(S)$,
    \item a capacity $\mu$ is said $k$-additive if its M\"{o}bius transformation is such that $m(T) = 0$ for all $T\subseteq G,|T|\geqslant k$. 
\end{itemize}
Because, as stated in \citep{Grabisch1997} and demonstrated by several real-world applications (see, e.g., \citep{grabisch2002subjective,berrah2007towards,angilella2018robust}), as well as in combination with other methodologies such as multi-objective optimization (see, e.g., \citep{branke2016using}), 2-additive capacities represent a useful compromise between a fully additive but simplistic model (a weighted sum, implying independence between criteria) and a fully general but difficult-to-handle Choquet integral model (which poses challenging elicitation issues). Therefore, in the following, we consider the Choquet integral in terms of a 2-additive capacity formulated as follows:
%Because 2-additive capacities are a useful compromise between a fully additive model (a weighted sum, implying independence) and a fully general Choquet integral model (raising difficult elicitation issues) \citep{MayagBouyssou2020}, in the following, let us recall the formulation of the Choquet integral in terms of a 2-additive capacity:
$$
U(a)=\displaystyle\sum_{g_i \in G} m(\{g_i\})\cdot g_i(a) + \displaystyle\sum_{\{g_i,g_{i'}\} \subseteq G} m(\{g_i,g_{i'}\})\cdot min\{g_i(a), g_{i'}(a)\}.
$$
$E^{Model}$ is therefore replaced by the following set of constraints:
$$
\left.
\begin{array}{l}
U(a)=\displaystyle\sum_{g_i \in G} m(\{g_i\})\cdot g_i(a) + \displaystyle\sum_{\{g_i,g_{i'}\} \subseteq G} m(\{g_i,g_{i'}\})\cdot min\{g_i(a), g_{i'}(a)\}\;\mbox{for all}\; a\in A,\\[1mm]
m(\{ \emptyset \}) = 0, \\[1mm]
\displaystyle\sum_{g_i\in G} m(\{ g_i\}) + \displaystyle \sum_{\{g_i, g_{i'}\} \subseteq G} m(\{g_i, g_{i'} \}) = 1,\\[1mm]
m(\{g_i \}) +\displaystyle \sum_{g_i \in T}^{} m(\{g_i, g_{i'}\} )\geqslant 0,\; \mbox{for all}\; g_i \in G,\; \mbox{and}\; T\subseteq G
\end{array}
\right\}E^{Model}_{Choquet}
$$
\noindent where the last three constraints are the equivalent of normalization (\textbf{2a)}) and monotonicity (\textbf{1a)}) constraints when a $2$-additive capacity is used.
Under the MCHP framework, the value of an alternative $a\in A$ w.r.t. macro-criterion $g_{\mathbf{r}}$ is given by 
$$
U_{\mathbf{r}}(a)=\displaystyle\sum_{\mathbf{t}\in E(g_{\mathbf{r}})}m\left(\{g_{\mathbf{t}}\}\right)\cdot g_{\mathbf{t}}\left(a\right)+\sum_{\{\mathbf{t}_1,\mathbf{t}_2\}\subseteq E(g_{\mathbf{r}})}m\left(\{g_{\mathbf{t}_1},g_{\mathbf{t}_2}\}\right)\cdot\min\{g_{\mathbf{t}_1}\left(a\right),g_{\mathbf{t}_2}\left(a\right)\}.
$$
\end{description}
Let us conclude this section observing that all mathematical programming problems discussed in the previous sections are linear because the considered preference models are affine in their parameters, that is, $U_{\alpha\mathbf{\pi}+(1-\alpha)\mathbf{\pi}'}(a)=\alpha\cdot U_{\mathbf{\pi}}(a) +(1-\alpha)\cdot U_{\mathbf{\pi}'}(a)$ for all $\alpha\in[0,1]$. In particular,
\begin{itemize}
\item if $U$ is a weighted sum, then  $\mathbf{\pi}=\left[w_1,\ldots,w_n\right]$,
\item if $U$ is a piecewise linear value function, then 
$\mathbf{\pi}=\left[u_i\left(y_i^{q}\right)\right]_{\substack{i=1,\ldots,n \\ q=2,\ldots,\gamma_i}}$,
\item if $U$ is a general additive value function, then 
$\mathbf{\pi}=\left[u_i\left(x_i^{f}\right)\right]_{\substack{i=1,\ldots,n \\ f=2,\ldots,m_i}}$,
\item If $U$ is a 2-additive Choquet integral, then, $\mathbf{\pi}=\left[\left[m\left(\{g_i\}\right)\right]_{g_i\in G},\left[m\left(\{g_i,g_{i'}\}\right)\right]_{\{g_i,g_{i'}\}\subseteq G}\right]$.
\end{itemize}

\section{Extended formulation of the LP problems to be solved in the example 2.2}
\label{sec:appendixB}
\begin{itemize} 
\item In example \ref{didacticexample2}:
\begin{center}
\begin{equation}\label{ExampleGA}
\begin{array}{l}
\;\;\;\;\overline{\sigma} = \min \left \{ \sigma^{+}(a) + \sigma^{-}(a) + \sigma^{+}(b) + \sigma^{-} (b) \right \},\;\;\mbox{subject to},\\[3mm]
\left.
\begin{array}{l}
\left.
\begin{array}{l}
U(a)= u_1(0.3) + u_2(0.7) \\[1,5mm]
U(b) = u_1(0.4) + u_2(0.6)\\[1,5mm]
U(c) = u_1(0.8) + u_2(1)\\[1,5mm]
u_1(0.3) \leqslant u_1(0.4) \leqslant u_1(0.8)\\[1,5mm]
u_2(0.6) \leqslant u_2(0.7) \leqslant u_2(1)\\[1,5mm]
u_1(0.3) = 0\\[1,5mm]
u_2(0.6) = 0\\[1,5mm]
u_1(0.8) + u_2(1) = 1\\[1,5mm]
\end{array}
\right \} E^{Model}_{GA}\\[1,5mm]
\;\;U(a) \geqslant U(b)\\[1,5mm]
\;\;U(a) - \sigma^{+}(a) + \sigma^{-}(a) = k\cdot 100\\[1,5mm]
\;\;U(b) - \sigma^{+}(b) + \sigma^{-}(b) = k\cdot 70\\[1,5mm]
\;\;k \geqslant 0\\[1,5mm]
\;\;\sigma^{+}(a),\; \sigma^{-}(a), \; \sigma^{+}(b), \;\sigma^{-}(b) \geqslant 0.
\end{array}
\right \} E^{DM}
\end{array}
\end{equation}
\end{center}
\begin{center}
    \begin{equation}
    \label{secondstepexample}
        \begin{array}{l}
        \;\;\ k^{*} = \max k,\; \text{subject to},\\[1,5mm]
        \left.
        \begin{array}{l}
             E^{DM}\\[1,5mm]
             \sigma^{+}(a)+\sigma^{-}(a)+\sigma^{+}(b)+\sigma^{-}(b) \leqslant \overline{\sigma} + \eta(\overline{\sigma}).
        \end{array}
        \right \} E^{DM^{'}}
        \end{array}
    \end{equation}
\end{center}
\end{itemize}
\end{document}